\documentstyle[amssymb,amsfonts]{amsart}

\newenvironment{proof}{\noindent {\bf Proof} }{\endprf\par}
\def \endprf{\hfill  {\vrule height6pt width6pt depth0pt}\medskip}
\def\emph#1{{\it #1}}
\def\textbf#1{{\bf #1}}

\newcommand{\bea}{\begin{eqnarray}}
\newcommand{\eea}{\end{eqnarray}}
\def\beaa{\begin{eqnarray*}}
\def\eeaa{\end{eqnarray*}}

\def\dcal{{\cal D}_1}
\def\dcall{{\cal D}_2}
\def\dcalll{\,^\star{\cal D}_1}
\def\dcallll{\,^\star{\cal D}_2}

\def\Dcal{{\cal D}}

\def\MM{{\cal M}}
\def\II{{\cal I}}
\def\ba{\begin{array}}
\def\ea{\end{array}}
\def\be#1{\begin{equation} \label{#1}}
\def \eeq{\end{equation}}

\newcommand{\nn}{\nonumber}

\def\rr{{\bf R}}

\def\nn{\nonumber}

\def\Gd{A}
\def\Bd{\und{A}}

\def\Lie{{\cal L}}
\def\tr{\mbox{tr}}
\def\err{\mbox{Err}}

\def\err{\mbox{Err}_{\ep}}
\def\Err{\mbox{Err}_{0}}

\def\gg{{\bf g}}

\newcommand\und{\underline}

\def\RR{{\cal R}}

\renewcommand{\div}{\mbox{div }}
\newcommand{\curl}{\mbox{curl }}

\def\a{\alpha}
\def\b{\beta}
\def\ga{\gamma}
\def\Ga{\Gamma}
\def\de{\delta}
\def\De{\Delta}
\def\ep{\epsilon}

\def\la{\lambda}
\def\La{\Lambda}
\def\si{\sigma}
\def\Si{\Sigma}
\def\om{\omega}
\def\Om{\Omega}

\def\rhoc{\check{\rho}}
\def\sic{\check{\si}}
\def\bboc{\check{\bb}}
\def\Roc{\check{R}}
\def\muc{\tilde{\mu}}
\def\th{\theta}
\def\ze{\zeta}
\def\nab{\nabla}

\def\bb{\underline{\b}}
\def\lap{\Delta}
\def\cga{\overset\circ{\ga}}

\def\hot{\widehat{\otimes}}
\newcommand{\trchb}{\tr \chib}
\newcommand{\chih}{\hat{\chi}}
\newcommand{\chib}{\underline{\chi}}

\newcommand{\etab}{\underline{\eta}}
\newcommand{\chibh}{\underline{\hat{\chi}}\,}
\def\f14{\frac{1}{4}}
\def\f12{{\frac{1}{2}}}

\def\c{\cdot}

\newcommand{\les}{\lesssim}

\newcommand{\HH}{{\mathcal H}}
\newcommand{\BB}{{\mathcal B}}
\newcommand{\PP}{{\mathcal P}}

\newcommand{\NN}{{\mathcal N}}

\def\Lb{\underline{L}}

\def\pr{\partial}

\def\chih{\hat{\chi}}
\def\trch{\mbox{tr}\chi}
\def\trchav{\overline{\trch}}
\newcommand{\ddd}{\nab}

\newcommand{\nnab}{\overline{\nab}}
\newcommand{\dddd}{\overline{\ddd}}

\begin{document}
\theoremstyle{plain}
  \newtheorem{theorem}[subsection]{Theorem}
  \newtheorem{conjecture}[subsection]{Conjecture}
  \newtheorem{proposition}[subsection]{Proposition}
  \newtheorem{lemma}[subsection]{Lemma}
  \newtheorem{corollary}[subsection]{Corollary}

\theoremstyle{remark}
  \newtheorem{remark}[subsection]{Remark}
  \newtheorem{remarks}[subsection]{Remarks}

\theoremstyle{definition}
  \newtheorem{definition}[subsection]{Definition}

\include{psfig}
\title[rough Einstein metrics ]{Causal geometry of Einstein-Vacuum 
spacetimes with finite curvature flux}
\author{Sergiu Klainerman}
\address{Department of Mathematics, Princeton University,
 Princeton NJ 08544}
\email{ seri@@math.princeton.edu}

\author{Igor Rodnianski}
\address{Department of Mathematics, Princeton University, 
Princeton NJ 08544}
\email{ irod@@math.princeton.edu}
\subjclass{35J10\newline\newline
The first author is partially supported by NSF grant 
DMS-0070696. The second author is a Clay Prize Fellow and is partially 
supported by NSF grant DMS-01007791
}
\vspace{-0.3in}
\begin{abstract}One of the  central difficulties
of settling the $L^2$-bounded curvature conjecture
for the Einstein -Vacuum equations is to be able to control
the causal structure of  spacetimes with such limited regularity.
 In this paper we 
show how to circumvent this difficulty by showing that the geometry
of null hypersurfaces of Enstein-Vacuum spacetimes
can be controlled in terms of initial data and the total
curvature flux through the hypersurface.
\end{abstract}
\maketitle
\section{Introduction}
We consider the Einstein Vacuum equations,
\be{eq:I1}
\rr_{\a\b}=0
\end{equation}
where  $\rr_{\a\b}$
denotes  the  Ricci curvature tensor  of  a four dimensional Lorentzian space time  $(\MM,\,  \gg)$.  
The fundamental problem of
the subject is to 
 study the long term regularity and asymptotic 
properties
 of the Cauchy developments of general, asymptotically flat,  
initial data sets\footnote{Though the results presented here  are local in nature and  apply
equally well to  more general data.} $(\Si, g, k)$. In so far as local regularity is concerned it 
is natural to ask what are the minimal local regularity properties of
the initial data which guarantee the existence and 
uniqueness of local developments. 

The first step 
in this   direction  was taken by 
 Y. C. Bruhat \cite{Bruhat} who proved 
a local in time existence result under the assumption\footnote{These are local Sobolev spaces on $\Si$;
more precise conditions can be given in terms of  weighted Sobolev spaces at infinity.} that 
$g\in H^s_{\mbox{loc}}(\Si)$, $k\in H^{s-1}_{\mbox{loc}}(\Si)$ with $k\ge 4$. The result, which
was  was
later improved to $s>5/2$ in \cite{HKM}, 
depended only on  energy estimates and Sobolev type inequalities
and did not require the study of null hypersurfaces.
It was conjectured\footnote{In \cite{Kl:2000} the issue was also raised 
 of going all the way to the {\sl critical exponent} $s_c=3/2$. We shall argue
below that this may not be possible,   see also 
footnote \ref{foot:H2sharp} below. } in 
\cite{Kl:2000} that one can  significantly improve  the result to $s=2$
 which corresponds to initial data sets with bounded 
$L^2$ curvature, result which would be 
  particularly satisfying in view of the naturalness
 of the norms involved. The conjecture
was
 motivated by the progress made earlier 
on semilinear type equations such as Wave Maps and Yang Mills equations.

At that time the conjecture was made  it seemed however completely out of reach. 
It was clear that in order
to improve  the exponent $s>5/2$ one had to abandon the naive use of
 Sobolev inequalities of  the classical
argument and rely instead on Strichartz and
bilinear type estimates.
 The problem was  that one needed  to extend these estimates to wave operators 
on very rough background metrics. The first results below $5/2$ are  due to Bahouri-Chemin \cite{Ba-Ch1},
 Tataru \cite{Ta},  \cite{KR:Duke}  and
relied, indeed, on proving Strichartz type estimates on such backgrounds. To do this
 they had to rely on adequate notions of
approximate fundamental
 solutions( or vectorfields in the case of  \cite{KR:Duke}) for the corresponding wave operators
based on an adapted version of the classical    geometric optics  construction. This  construction depends
heavily on the regularity  properties of null hypersurfaces associated to these backgrounds and  
  applies to general type 
of quasilinear wave equations including the  reduced Einstein vacuum equations in wave coordinates.  
In \cite{KR:Annals} we were able to reach, for the particular  case of  Einstein vacuum equations, any
exponent
$s>2$. A similar result was also  obtained in \cite{Sm-Ta} for
 the general class of quasilinear
wave equations mentioned above.

The case $s=2$ is far more difficult. First of all such a result
 cannot hold for general quasilinear wave equations,
see \cite{Lind}. As the experience with semilinear wave 
equations demonstrates, see
 discussion in \cite{Kl:2000} , to prove such a
result we need the following ingredients:
{\em\begin{enumerate}
\item Provide  a system of coordinates relative to which \eqref{eq:I1}
verifies an appropriate version of the null condition.
\item Prove  appropriate bilinear estimates for solutions 
to homogeneous wave equations, of the type $\square_\gg \phi=0$, on
 a fixed Einstein Vacuum  background( endowed with
the coordinate system indicated in 1.) 
with bounded $L^2$ curvature.
\end{enumerate}
}
To prove bilinear estimates we need a good notion of approximate solutions;
this requires  the third ingredient, typical to quasilinear equations,

\quad 3. {\em Make sense of null hypersurfaces, on Einstein vacuum  backgrounds with only
  $\smallskip   \qquad  L^2$ bounds on their curvature tensor,  and provide
appropriate estimates on }  
$\smallskip   \qquad$ {\em their geometry.  }

In this paper we make a first step towards settling the bounded  $L^2-$ curvature conjecture by 
providing the main ideas and techniques needed to  deal with the last ingredient  mentioned above.
More precisely we shall consider an outgoing  null
 hypersurface\footnote{  In this paper, for clarity,  we only consider 
 spherical geodesic foliations of a fixed null hypersurface. Nevertheless
our methods can be easily extended to deal with other foliations and 
other type of null hypersurfaces. They can also be extended to families
of null hypersurfaces as it would be in fact needed to tackle the 
$L^2$ curvature conjecture.
 } $\HH$, initiating
on a compact $2$ surface $S_0\subset\Si$ diffeomorphic to ${\bf S^2}$,   given by the level hypersurfaces 
of an optical function $u$, i.e. solution to the Eikonal equation 
\be{eq:I3}
\gg^{\a\b}\pr_\a u\pr_\b u=0.
\end{equation}
Let $L=-\gg^{\a\b}\pr_a u\, \pr_b $ be the corresponding null generator vectorfield and $s$ 
its affine parameter, i.e. $L(s)=1,\quad s|_{S_0}=0$. The level surfaces $S_s$ of $s$
generates the geodesic foliation on $\HH$. We shall denote
by $\nab$ the covariant differentiation
on $S_s$ and by $\ddd_L$ the projection to $S_s$ of the
covariant derivative with respect to $L$, see section 2. We also denote by $r$ the
function on $\HH$ defined by  $r=r(s)=\sqrt{(4\pi)^{-1}|S_s|}   $, with $|S_s|$ the area of $S_s$. Let
$\HH_t$ be the portion of $\HH$ between $s=0$ and $s=t$ and, for simplicity, assume $\HH=\HH_1$.
  We introduce the  total curvature flux\footnote{The justification for this quantity can
be found in \cite{Chr-Kl} or \cite{Kl-Nic} in connection to the Bell Robinson tensor and the Bianchi
identities verified by the curvature tensor $\rr$.} along
$\HH$ to be the integral, see precise definition in section 2, 
\be{eq:I3'}
\RR_0=\bigg(\|\a\|_{L^2(\HH)}^2+\|\b\|_{L^2(\HH)}^2+
\|\rho\|_{L^2(\HH)}^2+\|\si\|_{L^2(\HH)}^2+\|\bb\|_{L^2(\HH)}^2\bigg)^\f12
\end{equation} 
with $\a,\b,\rho,\si,\bb$ null   components of the curvature tensor $\, \rr$.
More precisely, for any vectorfields $X,Y$  on $\HH$ tangent to the $s$ foliation,
\be{eq:I3''}
\a(X,Y)=\rr(X,L,Y,L), \quad \b(X)=\f12 \rr(X,L,\Lb,L), \quad 
\rho=\frac 1 4 \rr(L,\Lb,L,\Lb)\ldots
\end{equation}
 see  definition \ref{def:nullcurv} for the other components $\si,\bb$.
 We shall assume that  the total curvature flux $\RR_0$   is finite
and show that the geometry of the null hypersurface $\HH_t$, for sufficiently small
 $t>0$, depends only on  initial conditions on $S_0$, which we denote by $\II_0$
and will make precise later,  and the size of  $\RR_0$. Alternatively
we will take $t=1$ and make $\RR_0$ and $\II_0$ sufficiently
 small\footnote{The smallness of $\II_0$ implies that the metric on $S_0$
is close to that of  the standard sphere. }.

The geometry of $\HH$ depends in particular of 
 the null second fundamental
form
\be{eq:I4}
\chi(X,Y)=<D_{X}L\,,\,Y>
\end{equation}
with $X,Y$ arbitrary vectorfields tangent to the  $s$-foliation. We 
denote by $\trch$ the trace of $\chi$,
i.e. $\trch=\de^{ab}\chi_{ab}$ where $\chi_{ab}$ are the components of $\chi$ relative to an
orthonormal frame $(e_a)_{a=1,2}$ on the leaves of the $s$-foliation. 
 It satisfies the well known Raychadhouri equation,
 see \eqref{eq:trchgeod},
\be{eq:I5}
\frac{d}{ds}\trch+\f12 (\trch)^2=-|\chih|^2
\end{equation}
with $\chih_{ab}=\chi_{ab}-\f12\trch\de_{ab}$ the traceless part of $\chi$.

Even in Minkowski space one cannot adequately control the geometry of null
hypersurfaces without a uniform bound on, at least, $\trch$. Indeed one can 
show  by an explicit calculation, see  section \ref{sect:blow-up}, that
if $L^\infty$ norm of $\trch$, on an initial 2-surface,  is allowed to become arbitrarily large 
caustics can occur instantaneously and therefore all $L^p$ norms, $p>2$ must
 become infinite\footnote{\label{foot:H2sharp}This  fact leads  us to conjecture that one cannot
control solutions of the Einstein equations with initial data rougher than $H^2$, i.e.
in  some $H^s$, $s<2$.  Thus the $H^2$ conjecture appears to be a natural endpoint result
!}.  To avoid caustics we  are thus  led to ask whether we can bound the $L^\infty$ norm of
$\trch$ simply in terms of the initial data norm $\II_0$ and total curvature flux $\RR_0$.
At first glance this seems impossible. Indeed, by integrating  \eqref{eq:I5}, in order to
control $\|\trch\|_{L^\infty}$ we need to control uniformly  the integrals
$\int_{\Ga}|\chih|^2$, along all  null  geodesic generators $\Ga$ of $\HH$.
 On the other hand  $\chih$ verifies a transport equation of the form,
\be{eq:I6}
\frac{d}{ds}\chih+\f12 \trch\c \chih=-\a
\end{equation}
where $\a$, the  null  curvature component defined  by \eqref{eq:I3''}, is only $L^2$ integrable
on $\HH$.   Thus, unless there is a miraculous gain of a spatial  derivative(along the
surfaces of the $s$- foliation) of the  integrals of $\a$ along the null geodesic generators $\Ga$, we have no chance to 
control the $L^\infty$ norm of $\trch$.  Fortunately this cancellation occurs and it is best
seen by observing that $\chih$ verifies a Codazzi type equation of the form,
\be{eq:I7}
\div\chih =-\b+\f12 \nab\trch+\f12  \trch\c\ze-\ze\c\chih
\end{equation}
 with $\b$ as defined   above and where $\ze$ denotes the torsion of the $s$ foliation, see  
\eqref{eq:torsion}.
As it was noted  and made use of  in \cite{Chr-Kl},
$\div\chih$ defines an elliptic system  on the leaves of the $s$ foliation
and therefore we expect\footnote{Assuming that the remaining quadratic terms in \eqref{eq:I7}
can also be controlled.} that $\chih$ behaves like $\Dcal^{-1}(-\b+\f12 \nab\trch)$ with $\Dcal^{-1}$
a pseudodifferential operator of order $-1$. We are thus  led to   control, uniformly, 
the integrals
$$\int_\Ga | \Dcal^{-1}(-\b+\f12 \nab\trch)|^2$$
along the null generators $\Ga$ 
or,  as there is no chance of an additional cancellation between the two terms, the integrals
\be{eq:I7'}
I_1=\int_\Ga | \Dcal^{-1}\b|^2, \quad I_2=\int_\Ga|\Dcal^{-1}\c\nab\trch|^2.
\end{equation}
 {\bf  Consider  first $I_1$}. Can we estimate it in terms of $\|\b\|_{L^2(\HH)}$ ?.
In our previous work\footnote{Similar arguments were also used in \cite{Sm-Ta}.}
 \cite{KR:Annals} we had interpreted $\b$ as corresponding to the term
 $\pr^2\gg + (\pr\gg)^2$, where $\gg$ is a space-time metric. In wave coordinates
   the Einstein metric $\gg$ solves a nonlinear wave equation which can 
 be written schematically in the form,
 $$
 \gg^{\mu\nu} \pr^2_{\mu\nu} \gg = (\pr \gg)^2.
 $$ 
 Thus, 
 $$
 \int_\Ga | \Dcal^{-1}\b|^2\approx \int_\Ga |\pr\phi|^2,
 $$
 where, in a first approximation, we can think of $\phi$ as an $H^2$ solution of a linear wave 
 equation in Minkwoski space, $\Box \phi =0$. 
 Unfortunately, it is known (see \cite{KM}), that the  integral 
 $ \int_\Ga |\pr\phi|^2$ can become arbitrary large, independent of length 
 of the null segment $\Ga$ and the size of data in $H^2$.
 This phenomenon, in fact,  was used as a counterexample to the Strichartz estimate
\be{eq:Str-s}
 \|\pr\phi\|_{L^2_t L^\infty_x} \les \|\phi_0\|_{H^s}
 \end{equation}
 for $s=2$. The Strichartz estimate \eqref{eq:Str-s}, however, holds true 
 for $s>2$. This suggests that an adapted version of it,  for  {\sl 
rough backgrounds}, 
 could be used to
control the  local geometry of null
hypersurfaces, provided that the
regularity assumptions 
 are consistent with $\gg\in H^s$,
$s>2$. This strategy was, in fact,
implemented
 in \cite{KR:Annals} and \cite{Sm-Ta}.
 
 Another way of looking at the question of whether $\int_\Ga | \Dcal^{-1}\b|^2$
 can be bounded by $\|\b\|_{L^2({\cal H})}$  is to 
 interpret it as a restriction problem.
 This makes sense dimensionally and it corresponds to a trace type theorem
of the type 
\be{eq:I8}\|U\|_{L^2(\Ga)}\les \|\nab U\|_{L^2(\HH)}.
\end{equation}
applied to the tensor  $U=\Dcal^{-1}\b$. 
 Unfortunately it is well known
that this type of trace theorems, just as the 
 sharp Sobolev embedding into $L^\infty$, fail to be true. 
  We might  be able
to overcome this difficulty if we could write  $\b=\ddd_LQ$ with $Q$ a tensor
verifying\footnote{We denote by $\nnab=(\nab, \ddd_L)$ 
 all first derivatives on $\HH$ and by $\nnab^2 Q$ all second derivatives. } 
$\|\dddd^2
Q\|_{L^2(\HH)}+\| Q\|_{L^2(\HH)}\les
\RR_0$. Indeed such a result holds true in flat space\footnote{That is 
the estimate $\big(\int_\Ga|\ddd_LQ|^2\big)^\f12\les
 \|\nnab^2Q\|_{L^2(\HH)}+\|Q\|_{L^2(\HH)}$
holds true in flat space.}. The breakthrough,  which allows us to make
use of an appropriately adapted version of the   sharp  trace theorem
mentioned above, is provided by the Bianchi identities,
\be{eq:I9}
D_{[\muc}\rr_{\a\b]\ga\de}=0
\end{equation}
which, expressed   relative to our null pair $L,\,\Lb$, takes the form 
of a system of first order equations connecting the $\ddd_L$ derivatives 
of the null  components $\a, \b, \rho,\si\ldots$ to their  spatial derivatives $\nab$. 
In particular we have, ignoring the quadratic and higher order terms\footnote{It turns
out that the ignored  terms are not so easy to treat, we shall need to make
various renormalizations. In particular  we shall need to renormalize
the Bianchi identities \eqref{eq:I10} by introducing
new quantities $(\rhoc,\sic)$, see \eqref{eq:renormalizedcurv}. },  see
\eqref{eq:bianchro},\eqref{eq:bianchsi} for precise formulas,
\be{eq:I10}
\div\b=L(\rho)+\ldots,\quad
\curl\b=-L(\si)+\ldots
\end{equation}
Once more this is an elliptic Hodge system in $\b$ and we can write, formally
$$
\Dcal^{-1}\b=\Dcal^{-2}L(\rho,-\si)+\ldots=\ddd_L\big(\Dcal^{-2}(\rho,-\si)\big)
+[\Dcal^{-2},\ddd_L](\rho,-\si)+\ldots
$$
with $\Dcal^{-2} $ a pseudodifferential operator of order $-2$.
Thus, ignoring for the moment the  commutator $[\Dcal^{-2},\ddd_L](\rho,-\si)$ 
 and the other
error terms, we have indeed
\be{eq:I11}
\Dcal^{-1}\b=\ddd_LQ+\ldots,\qquad Q=\Dcal^{-2}(\rho,-\si)
\end{equation}
and one can check, with the help of the null Bianchi identities and $L^2$ elliptic 
estimates,\footnote{Ignoring the multitude of error terms generated in the process.
In reality the error terms require delicate arguments.} that
$\|\dddd^2Q\|_{L^2(\HH)}+\|Q\|_{L^2(\HH)}$ can indeed  be bounded by $\RR_0$. 

This circle of ideas seem to take care of the integral  $I_1$ in \eqref{eq:I7'}.
Unfortunately we hit another  serious obstacle with $I_2$. The problem is that 
the operator $\Dcal^{-1}\c\nab $ is a nonlocal operator of order zero and therefore does not
map $L^\infty$ into $L^\infty$. To overcome this difficulty we are forced
to try to prove a stronger estimate for $\trch$. The idea is to try to prove
the boundedness of $\trch$ not only in $L^\infty$ but rather in a Besov space
of type $B^1_{2,1}(S_s)$ which  both imbeds in $L^\infty(S_s)$ and is stable relative
to operators of order zero. This simple, unavoidable fact,
 forces us to work with spaces defined by Littlewood-Paley projections
and therefore adds a lot of technical baggage to this work.
Moreover the standard LP- theory, based on Fourier transform,
would require the use of local  coordinates on $\HH$. The
most natural coordinates would be those transported by
the hamiltonian flow generated by $L$. Yet these {\sl transported coordinates}
 lose $1/2$ derivative, relative to what one expects in view of the
finite flux assumption( see proposition \ref{prop:weak-regular-cond}), and thus  we  
prefer   to rely on an invariant,  geometric,  version of the LP-theory.  We develop
such a theory, together with
an appropriate  paradifferential calculus,   in \cite{KR2} by following   the heat flow
approach of \cite{S}. An informal introduction to the sharp trace type
theorems in  Besov  spaces needed in our work  and the geometric 
LP theory on which they rely is given in \cite{KR3}.

  Yet another difficulty
is due to the presence of the torsion element $\ze$ in the Codazzi  equation \eqref{eq:I7}.
It turns out that the estimates for $\ze$ are at least as complex as those for $\chi$.
Moreover they  depend on another parameter of the foliation, the extrinsic second
fundamental form $\chib$ for 
which we also need to provide delicate estimates.
We shall provide a more complete informal description
of the  main circle of ideas of how to treat
$\chi, \ze, \chib$  in section 4 after we  make a full discussion
of the null structure and Bianchi identities.

In addition to the difficulties already mentioned we have  to confront another 
unexpected hurdle;  to arrive at
the crucial formula \eqref{eq:I11} one needs to estimate 
the  commutators between the $\ddd_L$ derivative and 
the pseudodifferential operator $\Dcal^{-2}$ applied to the curvature components 
$(\rho,-\si)$. It turns out that far from being lower order
these commutators have the same level of differentiability
as the principal terms and require, themselves, a nontrivial
 renormalization\footnote{We deal with 
them by an infinite sequence of renormalizations, see lemma \ref{le:decompbandmore}. }. 
 This difficulty, typical to end point results, is also a major technical complication
in the proof of the  sharp trace type results on which this work relies. 

 We are now ready to state
our main result in a preliminary  simplified version.
A full version of the
theorem can be found in section 3.

{\bf Main Theorem.}{\em(First version)\quad Consider   an outgoing  null hypersurface $\HH=\HH_1$,
initiating on a compact $2$ surface $S_0\subset\Si$ diffeomorphic to ${\bf S^2}$ and
foliated by the geodesic foliation associated to the affine parameter $s$
with $s|_{S_0}=0$. Assume that both the initial data quantity 
$\II_0$  and  the total curvature flux $\RR_0$  are 
sufficiently small. Then\footnote{Observe that $\frac 2 r$ is the value of $\trch$
 for an  outgoing null cone in Minkowski space.}  ,
\be{eq:I12}
\|\trch-\frac{2}{r}\|_{L^\infty(\HH)}\les \II_0+\RR_0
\end{equation}
Additional estimates hold true for $\chih$, $\ze$ and $\chib$.}

The proof of the theorem requires not only the material of this paper
but also two other separate  papers. In \cite{KR2}  we discuss a geometric version
of Littlewood-Paley theory  needed in our proof. In \cite{KR3} we use
 this geometric LP theory 
to prove the trace type theorems to which we have alluded above as 
well as   other technical results.  
In this paper we rely  on the  results of \cite{KR2}--\cite{KR3}   to prove 
our main theorem. 

In section 2 we discuss the main  geometric quantities associated 
to null hypersurfaces $\HH$ , discuss the null structure and null
Bianchi identities as well as various renormalizations needed later.
 We refer the reader to the relevant  material in \cite{Chr-Kl}
and \cite{Kl-Nic}.

In section 3 we formulate a precise version 
of the main theorem and 
give a more elaborate description
of the main ideas in its proof. We also formulate 
the simplest bootstrap  assumption {\bf BA1} and discuss
some important simple consequences of it.

In section 4 we introduce our other bootstrap assumptions {\bf BA2}--{\bf BA4}
  and discuss some of their basic consequences.  We prove
that the surfaces of our foliation satisfy the {\bf WS} condition
 which allows us to use  the results of \cite{KR2}. We  also derive estimates
 for the Gauss curvature $K$. We show that the $L^2$ norm  of $K$ on the null
hypersurface $\HH$  is bounded by the total curvature  flux. We also
 give uniform bounds   on the the $L^2$ norm of some negative fractional derivative 
of $K$ restricted to the leaves of the foliation. This important  latter result
requires a commutator argument whose proof is postponed to the Appendix.
  We also provide 
proofs for our main $L^2$-elliptic
estimates for  Hodge systems.

In section 5 we describe our main  sharp, Besov, bilinear trace theorems
and  Besov product  results needed in the proof
of the main theorem. Most results stated in this section
require complicated arguments  based on the geometric 
LP-theory developed in \cite{KR2}--\cite{KR3}.

In section 6 we describe the structure and  the estimates  of the main 
commutator terms encountered in the  proof of the main  theorem. 

In section 7 we rely on  the results of previous sections to complete 
the  proof  of the main theorem.

In the Appendix we recall the definition
of fractional powers of the Laplace-Beltrami operator 
on the leaves of the  geodesic foliation of $\HH$, based
on the heat flow approach developed in \cite{KR2}, and provide the proof of 
the commutator estimate needed in section 4. 

\section{Null Hypersurfaces} We recall in this section
the basic geometric concepts associated to null hypersurface
and refer the reader to \cite{Chr-Kl}
and \cite{Kl-Nic} for a complete derivation 
of the null structure and null Bianchi
identities.

We consider a null hypersurface
$\HH$  of  an Einstein vacuum space-time $(\MM, g)$. We assume
 that $\HH$ is foliated by the level surfaces  $S_v$ of a
given function $v$ with $S_0$ fixed on an initial space-like hypersurface 
$\Si_0$. Let $L$ be the null
geodesic generator of
$\HH$,
 i.e. $$<L,L>=0,\qquad D_LL=0$$
with $<\,,\,>= <\,,\,>_g$ denoting the metric  of $\MM$.
At any point $P\in S_v\subset \HH$  we denote by 
$\Lb$ the null vector conjugate to $L$ relative
to the $S_v$ foliation, i.e. 
$$<L\,,\,\Lb>=-2\,,\qquad <\Lb, X>=0 \quad \mbox{for all } \, X\in T_p(S_v).$$
We shall say that $L,\Lb$ form the canonical null  pair associated
to the foliation.
\begin{remark}\label{re:constants}
Observe that the null pair is uniquely defined up
to a  function constant along  the null geodesic generators of $\HH$.
\end{remark}
We shall denote by $\ga$ the induced metric on $S_v$, by $\nab$ the induced
covariant derivative and $K$ the Gauss curvature. An arbitrary orthonormal frame
on $S_v$ will be denoted by $(e_a)_{a=1,2}$. Clearly,
$$<e_a\,,\, L>=<e_a\,,\,\Lb>=0, \quad <e_a,e_b>=\de_{ab}.$$
A null pair together with an orthonormal frame $(e_a)_{a=1,2}$
as above is called a null 
frame associated to the foliation.
We recall, see \cite{Chr-Kl}, \cite{Kl-Nic}, the  definitions of the following basic geometric
quantities:
\begin{definition} The null second fundamental forms $\chi,\chib$
 of the foliation $S_v$ are given by 
\be{eq:nullsforms}
\chi_{ab}=<D_aL\,,\, e_b>,\qquad \chib_{ab}=<D_a\Lb\,,\, e_b>
\end{equation} 
The torsion is given by,
\be{eq:torsion}
\ze_a=\f12 <D_a L\,,\,\Lb>
\end{equation}
\end{definition}
\begin{remark} \label{framechangeS}In view of the remark \ref{re:constants}
we have the  freedom to replace $L$, on the original surface $S_0$,
by $L'=\om L$ where $\om $ is an arbitrary scalar function
on $S$. Therefore, choosing  $\Lb'=\om^{-1}\Lb$ we have the transformation
formulae,
\be{eq:framechangeonS} \chi'=\om \chi,\qquad \chib'=\om^{-1}\chib,\qquad \ze'=\ze-\nab\log\om
\end{equation}   
\end{remark}
\begin{definition}\label{def:nullcurv}
The null components of the curvature tensor
$R$ of the space-time metric $g$ are given
by:
\beaa
\a_{ab}&=&R(L,e_a,L, e_b)\,,\qquad \b_a=\f12 R(e_a, L,\Lb, L) ,\\ \rho&=&\frac{1}{4}R(\Lb, L, \Lb,
L)\,,
\qquad
\si=\frac{1}{4}\, ^{\star} R(\Lb, L, \Lb, L)\\
\bb_a&=&\f12R(e_a,\Lb,\Lb, L)\,,\qquad \underline{\a}_{ab}=R(\Lb, e_a,\Lb, e_b)
\eeaa
where $^\star R$ denotes the Hodge dual of $R$. The null decomposition 
of $^\star R$ can be related to that of $R$ according to the formulas, see\cite{Chr-Kl} :
\beaa
\a(^\star R)&=&-^\star\a(R), \quad \b(^\star R)=-^\star \b(R), \quad \rho(^\star R)=\si(R)\\
\si(^\star R)&=&-\rho(R), \quad \bb(^\star R)=-^\star \bb(R), \quad \underline{\a}(^\star
R)=^\star\underline{\a}(R)
\eeaa
\end{definition}
Observe that all tensors defined above are $S_v$-tangent, or simply
S-tangent when there is no danger of confusion,
i.e., tangent to the leaves of the $v$- foliation.

\begin{definition} The particular case of the foliation
defined by the level surfaces of the affine parameter $s$
of the null geodesic generator $L$ of $\HH$, i.e.,
$$L(s)=1, $$
is called the geodesic(or background)  foliation of $\HH$. Given
any other foliation $S_v$ we denote by 
\be{eq:nulllapse}
\frac{dv}{ds}=\Om_v=\Om
\end{equation}
the null lapse of the $v-$ foliation.
\end{definition}
To compare two arbitrary foliations it suffices to compare
any one of them with the geodesic foliation $S_s$. Let
$L,\Lb, e_a$ denote a null frame  at $P\in\HH$ for the geodesic foliation.
Then the vectors 
 $L,\Lb',e_a'$ with
\be{eq:framechange}
e_a'=e_a-\Om^{-1}e_a(v)L ,\qquad \Lb'=\Lb-2 \Om^{-1} e_{a}(v) e_{a} - 
2\Om^{-2}|\nab v|^{2} L.
\end{equation}
form a null  frame for  the  $v$ foliation at $P$.    Indeed, $e_a'(v)
=e_a(v)-\Om^{-1}e_a(v)\frac{dv}{ds}=0$. Taking into account that $L$ is orthogonal
to all tangent vectors to $\HH$ we  easily check that
$<e_a', e_b'>=<e_a,e_b>=\de_{ab}$ and $<\Lb', \Lb'>=<\Lb', e_a'>=0$ while $<L,\Lb'>=-2$.
\begin{proposition} The null second fundamental form
$\chi$ and curvature component $\a$ are intrinsic to the null hypersurface $\HH$,
i.e. they do not change by passing from one foliation to another.
\be{eq:invariance}
\chi'_{ab}=\chi_{ab},\qquad \a_{ab}'=\a_{ab}
\end{equation}
The torsion $\ze$ verifies the following transformation formula
\be{eq:invariancezeta}
 \ze_a'=\ze_a -
\Om^{-1} \chi_{ab}\, e_{b}(v)
\end{equation}
The curvature component $\b$ verifies,
\be{eq:transfbeta}
\b'_a=\b_a-\Om^{-1}e_b(v)\,\a_{ab}
\end{equation}
Similarly $\rho'-\rho,\si'-\si $ can be expressed as 
linear combinations of $\a,\b $ with
coefficients depending on $\Om^{-1}\nab v$
 and  $\bb'-\bb $ as  a linear combination
of  $\a, \b,\rho,\si$.

\label{prop:invariance}
\end{proposition}

\begin{proof}:\quad
Indeed  in view of the relations \eqref{eq:framechange} between the frames $\Lb', e_1', e_2'$
and $\Lb, e_1, e_2$ and the fact that $L$
is null geodesic,
\beaa
\chi'_{ab}&=&< D_{\,e'_a} L,e_{b'}>=<D_{e_a}L+\Om^{-1}e_a(v) D_LL\,,\,
e_b+\Om^{-1}e_a(v)L>=\chi_{ab}\\
2\ze'_a&=&<D_{e_a'}L, \Lb'>=<D_{e_a}L,\Lb'>=<D_a L\,,\, \Lb>- 
2\Om^{-2} e_{b}(v) <D_{a}L, e_{b}>\\ &=&2\ze_a -2\Om^{-1} \chi_{ab} e_{b}(v) 
\eeaa
as desired. 
To check the invariance of the null components of the curvature tensor
we write with the help of \eqref{eq:framechange}, 
\beaa
\a'_{ab}&=&R(L, e_a', L, e_b')=R(L,e_a, L, e_b)=\a_{ab}\\
 \b'_a&=&\f12 
R(e_a', L, \Lb', L)=\f12 R(e_a, L, \Lb-2\Om^{-1} e_b(v)e_b,  L)=\b_a-\Om^{-1} e_b(v)\a_{ab}\\
\rho'&=&\frac{1}{4}R(\Lb'-2\Om^{-1}e_a(v)e_a, L, \Lb' -2\Om^{-1}e_b(v)e_b,  L)\\
&=&\rho -2\Om^{-1}e_a(v)\b_a+\Om^{-2}E_a(v) e_b(v)\a_{ab}\\
\si'&=&\frac{1}{4}^\star R(\Lb'-2\Om^{-1}e_a(v)e_a, L, \Lb' -2\Om^{-1}e_b(v)e_b,  L)\\
&=&\si +2\Om^{-1}e_a(v)^\star\b_a-\Om^{-2}E_a(v) e_b(v)^\star\a_{ab}
\eeaa
The last equality holds in view of the formulas, see definition \ref{def:nullcurv}, relating the
null decomposition of
$R$ to that of its Hodge dual $^\star R$. We observe also that 
 the Hodge dual is invariant relative to \eqref{eq:framechange}, i.e. 
$\in_{L\Lb' e_1'e_2'}=\in_{L\Lb e_1 e_2}$. The transformation formula for $\bb$ can be derived in the same
manner.
\end{proof}
 We define now   the other   connection  coefficient 
associated to an arbitrary  $v$-foliation:
\be{eq:etab}
\etab_a=\f12<e_a\,,\, D_L\Lb>
\end{equation}
\begin{proposition} The Ricci coefficient $\etab$
associated to an arbitrary  $v$-foliation verifies:
\be{eq:etabzeta}
\etab_a=-\ze_a+e_a(\log\Om)=-\ze_a+e_a(\log\Om)
\end{equation}
\label{prop:etabzeta}
\end{proposition}
\begin{proof}:\quad
Consider the commutator $[L\,,\,e_a]$. We claim
that
$$[L\,,\,e_a]=-e_a(\log\Om) L+X$$
where $X$ is a vectorfield tangent to $S_v$.
Indeed,
this follows easily from the fact that 
$$[L\,,\,e_a](v)=Le_a(v)-e_aL(v)=-e_a (\Om)$$
Therefore,
\beaa
2\etab_a&=&<e_a\,,\, D_L\Lb>=-<D_Le_a\,,\, \Lb>=-<D_{e_a}L\,,\,\Lb>-<[L\,,\,e_a]\,,\,\Lb>\\
&=&-2\ze_a+2e_a(\log\Om)
\eeaa
as desired.
\end{proof}
\begin{definition}
The S-tensors $\chi,\ze, \chib, \etab$  form
the connection coefficients of the $v-$ foliation.
The following Ricci equations hold true:
 \bea
D_a\Lb&=&\chib_{ab} e_b+\ze_a \Lb\,,\qquad D_aL=\chi_{ab}e_b-\ze_a L\\
 D_L\Lb&=&2\,\etab_b\, e_b\,,\qquad\qquad D_L L=0\\
D_ae_b&=&\nab_a e_b+\f12 \chi_{ab} \Lb+\f12\chib_{ab}L
\eea
Also,
\be{eq:projection}D_Le_a=\ddd_Le_a+\etab_a L
\end{equation}\end{definition}
\begin{definition}
In general, given any tensor $\pi$ tangent
to the leaves of the $v$-foliation,   $\ddd_L \pi $ denotes the projection of $D_L \pi$ on $S_v$.
Clearly, $\ddd_L\ga=0$. A frame $(e_a)_{a=1,2}$ is said to be Fermi propagated if $\nab_L e_a=0$.
\end{definition}
\begin{definition}
  Given  an S-tangent one  form  $  F  $  we denote by $^\star   F  $ its
Hodge dual  $^\star  F  _a=\in_{ab}  F  _b$. Similarly if 
if $  F  $ is an S-tangent  symmetric traceless 2-tensor we define the left and right Hodge duals
$^\star   F  _{ab}=\in_{ac}  F  _{cb}$, $  F  ^\star_{ab}=  F  _a^c\in_{cb}$ and observe 
that $  F  ^\star=-^\star  F  $. For both one forms and symmetric traceless S-tensors $  F  $
we have $^\star(^\star  F  )=-  F  $. Given two S-tangent 
traceless symmetric(or one forms) tensors $  F  ,
\eta$ we denote by $  F  \wedge\eta$ the scalar product  between $  F  $ and the Hodge dual of
$\eta$. Thus for one forms,  $(  F  \wedge\eta) =(  F  \c^\star\eta)=\in_{ab}  F  _a  F  _b$
and for symmetric traceless tensors,
 $(  F  \wedge\eta) =(  F  \c^\star\eta)=\in_{ab}  F  _{ac}  F  _{bc}$.
\end{definition}
\begin{definition}
Given an S-tangent symmetric tensor $  F  $ we denote by $\tr  F  =\de^{ab}  F  _{ab}$ its trace
and by $\hat{  F  } $ its traceless part, i.e.
$\hat{  F  }_{ab}=  F  _{ab}-\f12(\tr  F  )\de_{ab}$
\end{definition}

\subsection{Structure  equations for null-foliations}
We shall use the general formulation of  [Chr-Kl].
With  the notation used on page 147 of [Chr-Kl] the  full parameters of the 
$v$ foliation verify the following:
$$
H=\chi,\, \bar{H}=\chib, \,Y=0, \,\bar{Z}=\etab,\, V=\ze, \Om=0
$$
The structure equations of the $v$-foliation are( in
view of the formula
on pages 168--169 of [Chr-Kl]):
\bea
L(\trch)&=&-|\chih|^2-\f12(\trch)^2\label{eq:transportchi}\\
\ddd_L\chih&=&-\trch\c \chih-\a\\
\div\chih& =&-\b+\f12 \nab\trch+\f12  \trch\ze-\ze\c\chih\\
L(\trchb)&=&2\div\etab  +2\rho-   \f12\trch\c \trchb-\chih\c\chibh+2|\etab|^2\label{eq:transportchih}\\
\ddd_L\chibh&=&\nab\hot\etab-\f12\trch\c\chibh -\f12 \trchb\c\chih
+\etab\hot\etab\label{eq:transportchibh}\\
\ddd_L\ze&=&-\b+\chi\c(-\ze+\etab)\label{eq:transportze}\\
\curl\etab&=&\f12\chih\wedge\chibh-\si\\
K&=&-\frac{1}{4}\trch\trchb +\f12\chibh\c\chih-\rho\\
\div\chibh&=&\f12\nab\trchb-\f12\trchb\ze+\ze\c\chibh+\underline{\b}
\eea
Recall also that $\etab$ and $\ze$
are related by the equation 
\be{eq:etabze2}
\etab=-\ze+\nab\log\Om
\end{equation}
where the null lapse  $\Om=\frac{dv}{ds}$ is  a free parameter. 
\subsection{ Null Bianchi equations}
In view of the formulas on page 161 of [Chr-Kl]
the Bianchi equations for $\b, \rho,\si$ of the $v$-foliation
are:
\bea
\ddd_L\b+2\trch\b&=&\div\a+(2\ze+\etab)\a\label{eq:bianchb}\\
L(\rho)+\frac{3}{2}\trch\rho&=&\div\b -\f12 \chibh \c\a+\ze\c\b +2\etab\c\b\label{eq:bianchro}\\
L(\si)+\frac{3}{2}\trch\si&=&-\curl \b +\f12 \chibh \wedge\a-\ze\wedge\b
-2\etab\wedge\b\label{eq:bianchsi}\\
\ddd_L\bb+\trch \bb &=&-\nab\rho+(\nab \si)^\star+2\chibh\c \b -3(\etab\c\rho-\etab^\star
\si)\label{eq:transportbb}
\eea
\subsection{Renormalized Null Bianchi}
The presence of  $\chibh$ on the right hand side \eqref{eq:bianchro}, 
\eqref{eq:bianchsi}
\eqref{eq:transportbb} creates serious difficulties in the applications
we have in mind. In what follows we remove this difficulty
by introducing the renormalized null curvature components 
\be{eq:renormalizedcurv}
\rhoc=\rho-\f12 \chih\c\chibh,\qquad  \sic=\si-\f12 \chih\wedge\chibh,
\qquad\bboc=\bb+2\chibh\ze
\end{equation}
 Using the transport equations \eqref{eq:transportchih} and then \eqref{eq:transportchih}
we notice that,
\beaa
-\chibh\c\a&=&\chibh\c(\ddd_L\chih +\trch\c\chih)\\
&=&L(\chih\c\chibh)+\trch\c\chih\c\chibh-
\chih\c(\nab\hot\etab-\f12\trch\chibh-\f12\trchb\c\chih+\etab\hot\etab)
\eeaa
Thus,
\beaa
L(\rho-\f12 \chih\c\chibh)+\frac{3}{2}\trch\c\rho&=&\div \b +\ze\c\b+2\etab\c \b+
\f12\trch\c\chih\c\chibh\\
&-&\f12\chih\c(\nab\hot\etab-\f12\trch\chibh-\f12\trchb\c\chih+\etab\hot\etab)
\eeaa
or
\bea 
L(\rho-\f12 \chih\c\chibh)+\frac{3}{2}\trch\c(\rho-\f12 \chih\c\chibh)&=&\div\b+(\ze+2\etab)\c\b
\label{eq:transportrho'}\\
&-&\f12\chih\c(\nab\hot\etab-\f12\trchb\c\chih+\etab\hot\etab).\nn
\eea
A similar equation holds 
for $\si$,
\bea \label{eq:transportsi'}
L(\si-\f12 \chih\wedge\chibh)+\frac{3}{2}\trch\c(\si-\f12\chih\wedge\chibh)
&=&-\curl\b -(\ze+2\etab)\wedge\b\\
&-&\f12\chih\wedge(\nab\hot\etab+\etab\hot\etab).\nn
\eea
We shall make a similar modification
for the transport equation verified by $\bb$. 
The idea is to eliminate $2\chibh\c \b$ from
the right hand side of \eqref{eq:transportbb}
with the help of \eqref{eq:transportze} and \eqref{eq:transportchibh}.
Indeed,
\beaa
2\chibh\c \b&=&2\chibh\c (-\ddd_L\ze+\chi\c(-\ze+\etab)\\
&=&-2L(\chibh\c\ze)+2\chi\c\chibh\c(-\ze+\etab)+2(\ddd_L\chibh)\c \ze\\
&=&-2L(\chibh\c\ze)+2(\nab\hot\etab-\f12\trch\c\chibh -\f12 \trchb\c\chih
+\etab\hot\etab)\c \ze\\&+&2\chi\c\chibh\c(-\ze+\etab)
\eeaa
Hence,
\bea
\ddd_L(\bb+2\chibh\c\ze) &=&-\nab\rho+(\nab\si)^\star  +2 (\nab\hot\etab)\c\ze   
-3(\etab\c\rho-\etab^\star
\si)-\trch \bb\label{eq:transportmodifiedbb}\\
&+&2\ze\c(-\f12\trch\c\chibh -\f12 \trchb\c\chih
+\etab\hot\etab)+2\chi\c\chibh\c(-\ze+\etab)\nn
\eea
\subsection{Commutation formulas}
\begin{proposition}
Consider an arbitrary k-covariant,
 S-tangent vectorfield $F_{\underline{a}}=F_{a_1\ldots a_k}$
Then,
\bea
\label{eq:commutationLnab}
\ddd_L\nab_b F_{\underline{a}}-\nab_b\ddd_L F_{\underline{a}}& =&-\chi_{bc}\nab_c
F_{\underline{a}}+(\ze_b+\etab_b)\ddd_LF_{\underline{a}}\\
&+&\sum_i(\chi_{a_i b}\,\etab_c-\chi_{bc}\,\etab_{a_i}-\in_{a_ic}\,^\star
\b_b)F_{a_1\ldots c\ldots a_k}\nn
\eea
In particular for scalars f,
\be{eq:commscalarf}
\ddd_L\nab_b f-\nab_b\ddd_L f=-\chi_{bc}\nab_c
f+(\ze_b+\etab_b)\ddd_L f
\end{equation}
Also, for a one form $F$,
\be{eq:commutationdiv}
L(\div F)-\div(\ddd_LF)=-\chi\c \nab F +(\ze+\etab)\c\ddd_L F -\f12 \trch\c\etab\c f
 -\chih\c\etab\c F -\b\c F
\end{equation}
Finally, for scalars $f$,
\bea
L(\De f)-\De(Lf)&=&-2\chi\c\nab^2 f+2(\ze+\etab)\c\nab L(f)+
(\div \ze +\div\etab+|\ze+\etab|^2)L(f)\nn\\
&+&\bigg(\f12 \trch(\ze+\etab)+\chih\c(\ze-\etab)+\nab\trch+\trch\ze\bigg)\nab f 
\label{eq:commlap}
\eea
\label{prop:commutationL}
\end{proposition}
\begin{remark}
In this paper we shall only make use of the above
commutator formulas for geodesic foliations, in which
case  $\ze+\etab=0$ and therefore the terms involving $\ddd_L$
 vanish.
\end{remark}
\begin{remark}In view of \eqref{eq:commutationdiv} one
would expect
the presence of a term of the form $\b\c\nab f$ on the right hand side
of \eqref{eq:commlap}. This terms cancells   however 
due to the Cdazzi equation for $\chih$.
of 
\end{remark}
\subsection{Mass aspect function}
Differentiating the equation \eqref{eq:transportze}
and using the commutation formula \eqref{eq:commutationLnab}
we infer,
\bea
\ddd_L(\nab_b\ze_a)&=& -\nab_b\b_a +
\nab_b\big(\chi_{ac}(-\ze_c+\etab_c)\big)-\chi_{bc}\nab_c\ze_a\label{eq:derivativeze}\\
&+& (\ze_b+\etab_b)\c\ddd_L\ze_a+\ze_c(\chi_{ab}\etab_c-\chi_{bc}\etab_a-\in_{ac}\,^\star\b_b)\nn
\eea
Contracting we derive,
\beaa
L(\div \ze)&=&-\div\b-(\ze-\etab)\c\div\chi+
\chi\c\nab(\etab-2\ze)+(\ze+\etab)\c\ddd_L\ze\\
&+&\f12\trch\ze\c\etab-\chih\c\ze\c\etab-\ze\c\b\\
&=&-\div\b-(\ze-\etab)\c (\nab\trch+\f12\ze\trch-\b-\ze\c\chih)     \\
&+&
\chi\c\nab(\etab-2\ze)+(\ze+\etab)\c\big(-\b+\chi\c(-\ze+\etab)\big)\\
&+&\f12\trch\ze\c\etab-\chih\c\ze\c\etab-\ze\c\b\\
&=&-\div\b-\ze\c\b-(\ze-\etab)\c\nab\trch+\chi\c\nab(\etab-2\ze)+\trch(-\ze\c\ze+\ze\c\etab+\f12
\etab\c\etab)\\ &+&\chih\c(-2\ze\c\etab+\etab\c\etab)
\eeaa
In other words,
\bea
L(\div \ze)&=&  -\div\b-\ze\c\b-(\ze-\etab)\c\nab\trch+\trch(\f12 \div\etab-\div\ze)+ \chih\c
(\f12\nab\hot\etab -\nab\hot\ze)\nn
\\
&+&\trch(-\ze\c\ze+\ze\c\etab+\f12
\etab\c\etab)
-\chih\c(\ze\hot\etab+\f12\etab\hot\etab)\label{eq:divzeta}
\eea

Combining this with \eqref{eq:transportrho'} and introducing the
mass aspect function 
\be{eq:massaspect} \muc =-\div \ze+\f12 \chih\c \chibh -\rho
\end{equation}
we derive the following:
\beaa
-L(\muc)&=&L(\div\ze)+L(\rho-\f12\chih\c\chibh)\\
&=&  -\div\b-\ze\c\b-(\ze-\etab)\c\nab\trch+ \trch(\f12\div\etab-\div\ze)+ \chih\c
(\f12\nab\hot\etab-\nab\hot\ze) \nn
\\
&+&\trch(-\ze\c\ze+\ze\c\etab+\f12
\etab\c\etab)
-\chih\c(\ze\hot\etab+\f12\etab\hot\etab)\\
&-&\frac{3}{2}\trch(\rho-\f12\chih\c\chibh)+\div\b+(\ze+2\etab)\c\b\\
&-&
\f12\chih\c(\nab\hot\etab-\f12\trchb\chih+\etab\hot\etab)\\
&=&-\chih\c(\nab\hot\ze) +2\etab\c\b -\frac{3}{2}\trch(\rho-\f12\chih\c\chibh)-(\ze-\etab)\c\nab\trch+
\trch(\f12\div\etab-\div\ze)\\
&+&\trch(-\ze\c\ze+\ze\c\etab+\f12
\etab\c\etab)-\chih\c(\etab\hot\etab+\ze\hot\etab)+\frac{1}{4}\tr\chib\chih\c\chih
\eeaa
We have thus proved the following
\begin{proposition}
The mass aspect $\muc $ verifies the following transport equation:
\bea
L(\muc)+\trch\muc&=&\chih\c(\nab\hot\ze)-2\etab\c\b +\trch(\rho-\f12\chih\c\chibh)+
(\ze-\etab)\c\nab\trch\nn\\
&-&\f12\trch\div\etab -\trch(-\ze\c\ze+\ze\c\etab+\f12
\etab\c\etab)\nn\\
&+&\chih\c(\ze\hot\etab+\etab\hot\etab)-\frac{1}{4}\tr\chib\chih\c\chih\label{eq:masstransport}
\eea
\end{proposition}

This transport equation should be viewed
in combination with  the following  Hodge elliptic system
for $\ze$:
\bea
\div\ze&=&-\muc+\f12 \chih\c\chibh-\rho\\
\curl\ze&=& -\curl\etab=\si-\f12 \chih\wedge\chibh
\eea
\subsection{Transport equation for $\nab\trch$} For future reference we also record
the transport equation for $\nab\trch$. This can be easily derived by
differentiating the transport equation \eqref{eq:transportchi} for $\trch$ and then use the commutation
formulas \eqref{eq:commutationLnab}. We derive,
\bea
\ddd_L(\nab\trch)&=& -\frac{3}{2} \trch\c\nab \trch-\chih\c \nab\trch-2\chih\c \nab \chih\\
&-&(\ze+\etab)\big(|\chih|^2+\f12 (\trch)^2\big)\nn
\eea

\subsection{Geodesic  foliation}
The simplest foliation of a null hypersurface $\HH$ is given by the geodesic foliation
$\Omega=1$.  In view of proposition \ref{prop:etabzeta}, see also   \eqref{eq:etabze2},
we infer that 
\be{eq:etabgeodesic}
\etab=-\ze.
\end{equation}
The structure equations of the geodesic foliation
take the form,
\bea
\frac{d}{ds}\tr\chi&=&-\frac 12 (\tr \chi)^2-|\chih|^2 \label{eq:trchgeod}\\
\nab_L\chih&=&-\trch\c \chih-\a\label{eq:chihgeod}\\
\ddd_L(\nab\trch)&=& -\frac{3}{2} \trch\c\nab \trch-\chih\c \nab\trch-2\chih\c \nab
\chih\label{eq:transportnabtrchi}\\
\nab_L\ze&=&-\trch\ze -2\chih\c\ze -\b\label{eq:zegeod}\\
\ddd_L(\nab\trch)&=& -\frac{3}{2} \trch\c\nab \trch-\chih\c \nab\trch-2\chih\c \nab \chih
\label{eq:transportnabtrch}\\
\frac{d}{ds}\trchb&=&-\f12 \trch \trchb-2\div\ze-\chih\c\chibh+2|\ze|^2+2\rho\label{eq:trchbgeod}\\
\nab_L\chibh&=&-\nab\hot \ze-\f12(\trch\c\chibh+\trchb\c\chih)+\ze\hot\ze\label{eq:chibhgeod}\\
 \div \chih&=&\f12\,\nab\trch-\chih\c\ze -\b\label{eq:dvichihgeod}\\
\curl \ze&=&-\f12\chih\wedge\chibh+\si\label{eq:divzegeod}\\
\div\chibh&=&\f12\nab\trchb-\f12\trchb\ze+\ze\c\chibh+\underline{\b} \label{eq:dvichibhgeod}\\
K&=&-\frac{1}{4}\trch\trchb+\f12 \chih\c \chibh -\rho\label{eq:Gausseq}.
\eea

We also have the renormalized  Bianchi identities,
\bea
\nab_L\b&=&\div \a+\ze\c\a\label{eq:bianchibeta}\\
L(\rhoc)+\frac{3}{2}\trch\c\rhoc &=&\div\b-\ze\c\b
+\f12\chih\c(\nab\hot\ze+\f12\trchb\c\chih-\ze\hot\ze).\label{eq:transportrho''}\\
L(\sic)+\frac{3}{2}\trch\c\sic
&=&-\curl\b +\ze\wedge\b
+\f12\chih\wedge(\nab\hot\ze-\ze\hot\ze)\label{eq:transportsi''}\\
\ddd_L(\bboc) &=&-\nab\rho+(\nab \si)^\star -2 (\nab\hot\ze)\c\ze   
+3(\ze\c\rho-\ze^\star
\si)-\trch \bb\label{eq:transportmodifiedbbgeod}\\
&+&2\ze\c(-\f12\trch\c\chibh -\f12 \trchb\c\chih
+\ze\hot\ze)-4\chi\c\chibh\c\ze\nn
\eea
where
\be{eq:rhocsic}
\rhoc=\rho-\f12 \chih\c\chibh, \qquad \sic=\si-\f12\chih\wedge\chibh, \qquad \bboc=\bb+2\chibh\c\ze.
\end{equation}

The transport equation for
the mass aspect function $\muc$ takes the form,
\bea
\frac{d}{ds} \muc+\trch\muc&=&2\ze\c\b+\trch(\rho-\f12\chih\c\chibh)+2\ze\c\nab\trch
+\f12\trch\c \div\ze+\chih\c(\nab\hot\ze)\nn\\
&+&\frac{3}{2}\trch\ze\c\ze-\frac{1}{4}\trchb\chih\c\chih\label{eq:masstransportgeodesic}
\eea
This transport equation should be viewed
in combination with  the following  Hodge elliptic system
for $\ze$:
\bea
\div\ze&=&-\muc+\f12 \chih\c\chibh-\rho=-\muc+\rhoc\label{eq:hodgeze}\\
\curl\ze&=& \si-\f12 \chih\wedge\chibh=\sic\nn
\eea
It turns out that the form  of the transport equation \eqref{eq:masstransportgeodesic}
is not quite convenient in applications. To eliminate the curvature terms on the right
hand side we observe, from the transport equation for $\ze$,
$$\frac{d}{ds}|\ze|^2+2\trch |\ze|^2=-4\chih\c\ze\c\ze-2\ze\c \b.$$
Adding this to \eqref{eq:masstransportgeodesic} and substituting
$\div\ze=-\muc-\rhoc$ we derive,
\bea
\frac{d}{ds} \mu+\frac{3}{2}\trch \mu&=&\chih\c(\nab\hot\ze)+\f12 \trch\rhoc+2\ze\c\nab\trch
\label{eq:newmasstransportgeodesic}\\ &-& 4\chih\c\ze\c\ze + \trch|\ze|^2 -\frac 1 4
\trchb|\chih|^2\nn
\eea
where
\be{eq:newmassaspect}
\mu=\muc+|\ze|^2=-\div\ze+\f12\chih\c\chibh -\rho+|\ze|^2
\end{equation}
\begin{remark} For a Fermi propagated frame $(e_a)_{a=1,2}$, i.e. 
$\nab_L e_a=0$,  we have  $(\nab_L\chi)_{ab}=\frac{d}{ds}\chi_{ab}, (\nab_L\ze)_a=\frac{d}{ds}\ze_a$
etc.,  and  therefore all corresponding structure  equations can be interpreted as differental
equations along null geodesics.
\end{remark}
\begin{remark}
As explained in Remark \ref{framechangeS}
 we have the choice of selecting the null vector $L$ on the initial surface
$S_0$.  Passing from
one fixed  null pair $L,\Lb$ to another $L'=\om L, \Lb'=\om^{-1} \Lb$
 we have $\chih'=\om\chih$,\, $\chibh'=\om^{-1}\chibh$,\,  $\ze'=\ze-\nab\log\om$.
It is easy to check that the combination $\nab\trch+\ze\trch$ remains invariant
 while $\muc\,'=\muc+\lap\log\om$.
\end{remark}
\begin{remark}
  In view of the remark above it is thus possible to choose a null pair $L,\Lb$ relative 
to which  the mass aspect function $\muc$ is constant on $S_0$. Indeed all we have to do is 
solve the equation $\lap\log\om=\muc-\overline{\muc}$, where $\overline{\muc}$
denotes  the average of $\muc$ on $S_0$.
\end{remark}
Finally we record below the following well known:
\begin{lemma} Relative to the null geodesic foliation on  $\HH$, 
\be{eq:derivaverage}
\frac{d}{ds}\int_{S_s} fdA_{s}=\int_{S_s}\big(\frac{d}{ds} f+\trch f\big) dA_s
\end{equation}
 where $dA_s=dA_{\ga_s}$ denotes the volume element on $S_s$, with  $\ga_s=\ga(s)$ the
induced metric on $S_s$. In particular
denoting by $r=r(s)=\sqrt{(4\pi)^{-1}|S_s|}   $, with $|S_s|$ the area of $S_s$, we have
\be{eq:drds}
\frac{dr}{ds}=\frac{r}{2} \overline{\trch}
\end{equation}
\label{le:transportintegrform}
\end{lemma}
\begin{proof}:\quad We parametrize the hypersurface $\HH$ py picking 
arbitrary coordinates $\om^1,\om^2$ on $S_0$ and following the trajectories
of the null generator vectorfield $L$. More precisely let $x^\mu=x^\mu(s,\om^1,\om^2)$
be the flow defined by:
\be{eq:geodflow}
\frac{dx^\mu}{ds}=L^\mu, 
\end{equation}
with $x^\mu(0,\om^1,\om^2)$ the point on $S_0$ of coordinates $(\om^1,\om^2)$
and $s$ the affine parameter of $L$.  We claim that relative to the coordinates $s,
\om^1,\om^2$ on $\HH$  the metric $\ga_s$ verifies,
\be{eq:transportga}
\frac{d}{ds}\ga_{ab}=2\chi_{ab}
\end{equation}

 Indeed  relative to the
coordinates
$s,
\om^1,\om^2$ on $\HH$  we have $L=\frac{\partial}{\partial s}$ and  since
 $[\frac{\partial}{\partial s},\frac{\partial}{\partial \om^a}]=0 $ we infer from
$\nab_L\ga=0$, and $\ga_{ab}=\ga(\frac{\partial}{\partial \om^a}, \frac{\partial}{\partial \om^b})$,
\beaa
0&=&(\nab_L\ga)(\frac{\partial}{\partial \om^a}, \frac{\partial}{\partial \om^b})=
\frac{d}{ds}\ga_{ab}-\ga(\nab_{\frac{\partial}{\partial s}}\frac{\partial}{\partial \om^a}, 
\frac{\partial}{\partial \om^b})-\ga(\frac{\partial}{\partial \om^a}, 
\nab_{\frac{\partial}{\partial s}}\frac{\partial}{\partial \om^b})\\
&=&\frac{d}{ds}\ga_{ab}-\ga(\nab_{\frac{\partial}{\partial \om^a}} L, 
\frac{\partial}{\partial \om^b})-\ga(\frac{\partial}{\partial \om^a}, 
\nab_{\frac{\partial}{\partial \om^b}}L )\\
&=&\frac{d}{ds}\ga_{ab}-2\chi_{ab}
\eeaa
as desired.

Therefore, denoting $|\ga|=\det(\ga_{ab})$, 
$\frac{d}{ds}\log |\ga|=\ga^{ab}\frac{d}{ds}\ga_{ab}=2\trch\quad$
or,
\be{eq:volumega}
\frac{d}{ds}\sqrt{|\ga|}=\trch\sqrt{|\ga|}
\end{equation}
Now,  relative to the coordinates $s, \om^1,\om^2$, 
$\int_{S_s} fd\muc_\ga=\int \int f\sqrt{|\ga|}\,  d\om^1d\om^2$, therefore,
$$\frac{d}{ds}\int_{S_s} fdA_s= 
\int \int \frac{d}{ds}(f\sqrt{|\ga|})\,  d\om^1d\om^2=
\int_{S_s}\big(\frac{d}{ds} f +\trch f\big)dA_s
$$
as desired.
\end{proof}
\begin{definition} We denote by $\ga_t(\om)$ the null geodesic $x^\mu(s,\om^2,\om^2)\subset\HH$  
initiating on $S_0$ at the point of coordinates $\om=(\om^1, \om^2)$ and ending on $S_{t}$.
Given a scalar function $f$ on $\HH$ we denote by 
$\int_{\ga_t(\om)} f$, or simply  $\int_{\ga_t} f$ when no confusion is possible,  the  integral 
$$\int_{\ga_t} f =\int_{\ga_t} f ds=\int_0^tf(s, \om^1, \om^2) ds.$$
\end{definition}
\begin{definition} We denote by $\HH_t$ the portion of the null hypersurface $\HH$
 corresponding to $0\le s\le t$.
Given a scalar function $f$ on $\HH$ we denote by $\int_{\HH_t} f$ the integral,
$$\int_{\HH_t} f=\int_0^t\int_{S_s} fdA_s=\int_0^t\int\int f \sqrt{|\ga|}(s,\om^1,\om^2) ds 
d\om^1 d\om^2$$
\end{definition}
\begin{remark} In view of the above definitions we easily check the co-area  formula,
\be{eq:coareacone}
\int_{\HH_t} f=\int_{S_0}\big(\int_{\ga_t}f\frac{\sqrt{|\ga_s|}}{\sqrt{|\ga_0|}}ds\big)
\end{equation}
\end{remark}

Observe that   $v_s=\frac{\sqrt{|\ga_s|}}{\sqrt{|\ga_0|}}$  verifies $\frac{d}{ds} v_s= \trch v_s$
 and therefore it can be used  as an integrating factor
according  the formula,
\be{eq:integratingfactor}
\ddd_L (v^k U)=v^k\big(\ddd_LU +k\, \trch U\big)
\end{equation}
for arbitrary $S$-tangent tensors $U$. 
 As a consequence of this formula,
 applied to scalar functions f,  we have:
\begin{proposition} The equation
$\frac{d}{ds}f +k\, \trch f=g,$
for scalars $f$, $g$ 
has the solution 
\be{eq:niceintegralform}
f(s,\om)= v_s^{-k}\bigg( f(0,\om)+ \int_0^s 
v_{s'}^{k} \,\,g(s',\om) ds  \bigg)
\end{equation}
\label{prop:formulatranspgeneral}
\end{proposition}
In what follows we shall use lemma \ref{le:transportintegrform}
to derive a transport  equation for $\overline{\trch}-\frac{2}{r}$. 
Using the transport equation  \eqref{eq:trchgeod} for $\trch$ and lemma \ref{le:transportintegrform}
we derive,
\beaa
\frac{d}{ds}\overline{\trch}&=&-\frac{2}{r}\frac{dr}{ds}\overline{\trch}
+\frac{1}{4\pi r^2}\int_{S_s}(\frac{d}{ds}\trch+\trch^2)\\
&=&-(\,\overline{\trch}\,)^2+\frac{1}{4\pi r^2}\int_{S_s}\big(\f12 (\trch)^2-|\chih|^2\big)\\
&=&-(\,\overline{\trch}\,)^2+\f12 \overline{\,(\trch)^2}-\overline{\,|\chih|^2}\\
&=&-(\,\overline{\trch}\,)^2+\f12
\overline{\,(\trch-\overline{\trch}+\overline{\trch})^2}-\overline{\,|\chih|^2}\\
&=&-\f12 (\,\overline{\trch}\,)^2+\f12 \overline{\,(\trch-\overline{\trch})^2}-\overline{\,|\chih|^2}
\eeaa
On the other hand
\beaa
\frac{d}{ds}(2r^{-1})&=&-2r^{-2}\frac{r}{2}\overline{\trch}=-r^{-1}\overline{\trch}
\eeaa
Thus,
\beaa
\frac{d}{ds}\big(\overline{\trch}-\frac{2}{r}\big)&=&-\f12 (\,\overline{\trch}\,)^2+\f12
\overline{\,(\trch-\overline{\trch})^2}-\overline{\,|\chih|^2}
+\frac{1}{r}\,\overline{\trch}\,\\
&=&-\f12\overline{\trch}\c\big(\overline{\trch}-\frac{2}{r}\big)
+\f12
\overline{\,(\trch-\overline{\trch})^2}-\overline{\,|\chih|^2}
\eeaa
or,
$$\frac{d}{ds}\big(\overline{\trch}-\frac{2}{r}\big)+\f12\big(\overline{\trch}-\frac{2}{r}\big)
=\f12\bigg(\big(\trch-\overline{\trch}\big) \big(\overline{\trch}-\frac{2}{r}\big)+
 \overline{\,(\trch-\overline{\trch})^2}\bigg)-\overline{\,|\chih|^2}.
$$
On the other hand,
\beaa
\frac{d}{ds}\big(\trch-\overline{\trch}\big)&=&-\frac{1}{2}(\trch^2-(\,\overline{\trch}\,)^2\big)-
\big(|\chih|^2-\overline{\,|\chih|^2}\big)\\
&=&-\frac{1}{2}(\trch-\,\overline{\trch}\,\big)   (\trch+\,\overline{\trch}\,\big)      -
\big(|\chih|^2-\overline{\,|\chih|^2}\big)
\eeaa
Thus,
$$\frac{d}{ds}\big(\trch-\overline{\trch}\big)+\f12\trch\big(\trch-\overline{\trch}\big)
= -\f12\overline{\trch } \big(\trch-\overline{\trch}\big) -\big(|\chih|^2-\overline{\,|\chih|^2}\big) $$
We infer that,
\begin{proposition} Let $W=\overline{\trch}-\frac{2}{r}$ and $V=\trch-\overline{\trch}$.
Then,
\be{eq:transpaveragetrch}
\frac{d}{ds} W+\f12\trch W=\f12( V\c W + V^2)-\overline{\,|\chih|^2}
\end{equation}
\be{eq:transptrchminusaveragetrch}
\frac{d}{ds} V+\f12\trch V=-\f12 \overline{\trch}\, \, V-\big(|\chih|^2-\overline{\,|\chih|^2}\big)
\end{equation}
\label{prop:travertrch}
\end{proposition}

\subsection{$L^2$ curvature norm associated to a  geodesic foliation}
Given a geodesic foliation on $\HH$    and $F$ an  arbitrary $S$-tensor
on $\HH$, and denoting by  $\HH_{t}$ the portion of $\HH$
 contained in the interval $s\in [0,t]$,   we define its $L^2$ norm to be 
\be{eq:L2normH}
\|F\|_{L^2(\HH_{t})}=\big(\int_0^{t} ds\int_{S_s}  |F|^2 \big)^\f12
\end{equation}
where $ \int_{S_s}  |F|^2$
 denotes the integral with respect to  the
volume element $d\muc_s$ of  $S_s$. If $\HH=\HH_1$
we shall simply write $\|F\|_{L^2}=\|F\|_{L^2(\HH)}$.
\begin{definition}
Consider $\a, \b,\rho,\si, \bb$ to be the null decomposition
of the curvature tensor relative to the geodesic  foliation.
We  consider the norm:
\be{eq:L2curv}
\RR[t]=\bigg(\|\a\|_{L^2(\HH_t)}^2+\|\b\|_{L^2(\HH_t)}^2+
\|\rho\|_{L^2(\HH_t)}^2+\|\si\|_{L^2(\HH_t)}^2+\|\bb\|_{L^2(\HH_t)}^2\bigg)^\f12
\end{equation} 
We say that the  null hypersurface $\HH=\HH_1$ has bounded $L^2$
curvature flux, relative to the geodesic-foliation, if
$\RR_0=\RR[1]<\infty$.
\end{definition}
We shall also make us of the following notations:
\begin{definition} Given an arbitrary $S$-tangent tensor $F$ on $\HH=\HH_1$  we denote 
$$\|\nnab F\|_{L^2}=\|\nab F\|_{L^2}+\|\ddd_L F\|_{L^2}.$$
We also introduce the following  norms ,
\beaa
\NN_1(F)&=&\|F\|_{L^2}+\|\nnab F\|_{L^2}\\
&=&\|F\|_{L^2}+\|\nab F\|_{L^2}+\|\ddd_L F\|_{L^2}\\
\NN_2(F)&=&\|F\|_{L^2}+\|\nnab F\|_{L^2}+\|\nab\nnab F\|_{L^2}\\
&=&\|F\|_{L^2}+\|\nab F\|_{L^2}+\|\nab^2 F\|_{L^2}\\
&+&\|\ddd_L
F\|_{L^2}+\|\nab\c\ddd_L F\|_{L^2}\\
\eeaa
We write $F\in\NN_1(\HH)$ or  $F\in\NN_2(\HH)$, if $\NN_1(F)<\infty$, resp. 
$\NN_2(F)<\infty$.
\label{def:NNnorms}
\end{definition}
\subsection{Hodge systems}
We consider the following Hodge operators acting on  $2$ surface $S$ 
diffeomorphic to  the standard $S^2$ sphere, such as  the leaves $S_v$ of
a given foliation:
\begin{enumerate}
\item The operator $\dcal $ takes any $1$-form $  F  $ into the pairs of
functions $(\div  F  \,,\, \curl  F  )$
\item The operator $\dcall$ takes any $S$ tangent symmetric, traceless
tensor $  F  $ into the $S$ tangent one form $\div  F  $.
\item The operator $\dcalll$ takes the pair of scalar functions 
$(\rho, \si)$ into the $S$-tangent 1-form $-\nab \rho+(\nab \si)^\star$.
\item The operator $\dcallll$ takes 1-forms $F$  on $S$ into  the 2-covariant, symmetric,
traceless tensors $-\f12 \widehat{\Lie_F\ga }$ with $\Lie_F\ga$ the traceless part of the 
Lie derivative of the metric $\ga$ relative to $F$, i.e.
$$\widehat{(\Lie_F\ga)}_{ab}=\nab_b F_a+\nab_a F_b-(\div F)\ga_{ab}.$$

\end{enumerate}
Observe that the kernels of both $\dcal$ and $\dcall$ in $L^2(S)$ are 
trivial and that $\dcalll$, resp. $\dcallll$ are  the $L^2$ adjoints of
 $\dcal$, respectively $\dcall$. The
kernel of $\dcalll$ consists of pairs of  constant functions $(\rho, \si)$ while
that of $\dcallll$ consists 
of the set of all conformal Killing vectorfields on $S$. In particular  the  $L^2$- range 
of $\dcal$ consists of all pairs of functions $\rho, \si$ on $S$ with
vanishing mean curvature. The $L^2$ range of $\dcall$ consists of all
$L^2$ integrable 1-forms on $S$ which are orthogonal
to the Lie algebra of  all conformal Killing vectorfields
on $S$, see proposition \ref{prop:hodgeident}. Accordingly we shall
consider the inverse operators $\dcal^{-1}$ and $\dcall^{-1}$ and
implicitly assume that they are defined on the  $L^2$ subspaces 
identified above.

 Finally we record the following simple
identities,
\bea
\dcalll\c\dcal&=&-\lap+K,\qquad \dcal\c\dcalll=-\lap\label{eq:dcalident}\\
\dcallll\c\dcall&=&-\f12\lap+K,\qquad \dcall\c\dcallll=-\f12(\lap+K)\label{eq:dcallident}
\eea

\subsection{Null Structure equations in Minkowski space}\label{sect:blow-up}
In this section we show how to derive explicit formulas
for the null structure equations in the  flat case of the  Minkowski
 space. These explicit formulas are well known, see for example \cite{F-S}.

We consider a null hypersurface  $\HH$ in Minkowski space.
Let $S_0\subset \Si_0$ an initial 2-surface with
$\Si_t$ the standard space like
hypersurfaces generated by $t$.
Clearly $L=T+N$, $\Lb=T-N$
where $T$ denotes the unit normal to $\Si_t$
while $N$ is the exterior unit normal to $S_t=\Si_t\cap \HH$ on $\Si_t$.
Clearly the affine parameter $s$  of $L$ coincides with $t$.
We also observe that $\ze=0$ and 
$\chi=-\chib=\th$ where $\th$
denotes the null second fundamental form
of  the surfaces $S_t$ with respect to $N$.
Thus the geometry of $\HH$ is entirely
described by  $\chi$, which we interpret here as $2\times 2$
 matrix with eigenvalues $\la$ and $\muc$,  and we have the matrix equation,
\be{eq:matrixchieq}
\frac{d}{ds}\chi=-\chi^2
\end{equation}
We look for solutions of \eqref{eq:matrixchieq}
of the form $\chi(s)=U\c\La(s)\c U^{-1}$ with $U$
independent of $s$ and $\La$ diagonal. Therefore,
$$\frac{d}{ds}\la=-\la^2,\qquad \frac{d}{ds}\mu=-\mu^2.$$
and we infer that,
\be{eq:exactsolutionlamu}\la(s)=\frac{\la(0)}{1+s\la(0)},\qquad \mu(s)=\frac{\mu(0)}{1+s\mu(0)}
\end{equation}
Now let $$a=\trch=\la+\mu,\qquad b=\la\c\mu.$$
Recall  that  $\chi^2-a\c\chi+b\c I=0$, with $I$ the identity
matrix.  On the other hand, since 
$\chih=\chi-\f12 a\c I$, we find
\be{eq:chihmatrix}
|\chih|^2=2(-b+\frac{1}{4} a^2)
\end{equation}
Now, in view of \eqref{eq:exactsolutionlamu},
\beaa
a(s)&=&\frac{\la(0)}{1+s\la(0)}+\frac{\mu(0)}{1+s\mu(0)}\\
&=&\frac{a(0)+2sb(0)}{\big(1+s\la(0)\big)\c\big(1+s\mu(0)\big)}
\eeaa
and 
\beaa
b(s)&=&\frac{\la(0)}{1+s\la(0)}\c\frac{\mu(0)}{1+s\mu(0)}=
\frac{b(0)}{\big(1+s\la(0)\big)\c\big(1+s\mu(0)\big)}\\
\eeaa
Now, in view of \eqref{eq:chihmatrix}, we easily calculate
\beaa
-\f12 |\chih(s)|^2&=&
(b(s)-\frac{1}{4} a^2(s))=\frac{4b(0)\c\De(s)-\big(a(0)+2sb(0)\big)^2}{\De(s)^2}\\
&=&\frac{4b(0)-a(0)^2}{4\De(s)}=-\f12|\chih(0)|^2\c\frac{1}{\De(s)^2}
\eeaa
where $\De(s)=\big(1+s\la(0)\big)\c\big(1+s\mu(0)\big)=1+sa(0)+s^2b(0)$.

Therefore,
$$|\chih(s)|^2=\frac{|\chih(0)|^2}{\De(s)^2}$$
and we easily deduce that,
\be{eq:finalformula1}
   \trch(s)=\frac{\trch(0)+s\big(\f12 \trch(0)^2-|\chih(0)|^2\big)}{\De(s)},\qquad   
\chih(s)=\frac{\chih(0)}{\De(s)}
\end{equation}
It remains to identify the zeroes of $\De(s)$ i.e. 
$$s=-\la(0)^{-1},\qquad s=-\mu(0)^{-1}.$$
Clearly $\la(0)+\mu(0)=\trch(0)$, $\la(0)\c\mu(0)=b(0)=\frac{1}{4}\trch(0)^2-\f12 |\chih(0)|^2$.
Thus, 
\be{eq:geometriclamuformulas}
\la=\frac{\trch(0)-\sqrt{2}|\chih(0)|}{2},\qquad
\mu=\frac{\trch(0)+\sqrt{2}|\chih(0)|}{2}
\end{equation}
Thus,
for $s>0$, singularities only occur if  $\la<0$
i.e. 
\be{eq:singcondition}
\trch(0)<\sqrt{2}|\chih(0)|
\end{equation}
In that case, for $s=s_0=-\la(0)^{-1}$,  the numerator
of $\trch$ becomes 
\beaa
a(0)+2s_0 b(0)=\la(0)+\mu(0)-2\la(0)^{-1}\c\la(0)\c\mu(0)=\la(0)-\mu(0)
\eeaa
which can only vanish if $\la(0)=\mu(0)$ i.e. $\chih(0)=0$. Thus  if $\chih(0)\neq 0$
and condition \eqref{eq:singcondition} is satisfied then necessarily
$\trch$ must become infinite at $s_0=-\la(0)^{-1}$.
Consider now the volume element $dv_s=\sqrt{|\ga_s|}$. Recall that ,
see \eqref{eq:volumega},
$\frac{d}{ds}\sqrt{|\ga|}=\trch\sqrt{|\ga|}$. Thus,
\beaa
\frac{\sqrt{|\ga_s|}}{\sqrt{|\ga_0|}}&=&
\exp\int_0^s\big(\frac{\la(0)}{1+s'\la(0)}+\frac{\mu(0)}{1+s'\mu(0)}\big)ds'\\
&=&\big(1+s\la(0)\big)\c\big(1+s\mu(0)\big)=\De(s)
\eeaa
Now,
\beaa
\|\trch\|_{L^p(S_s)}^p=\int_{S_s} |\trch|^p dv_s=\int_{S_0}\frac{E^p(s)}{\De^p(s)}\De(s)dv_0
=\int_{S_0}\frac{E^p(s)}{\De^{p-1}(s)}dv_0
\eeaa
where
$E(s)=\trch(0)+s\big(\f12 \trch(0)^2-|\chih(0)|^2\big)$. 
Similarily,
$$\|\chih\|_{L^p(S_s)}^p=\int_{S_0} |\chih|^p\frac{1}{\De^{p-1}(s)}dv_0.$$
We infer that for  all  $p\ge 2$ the  $L^p$ norms of $\trch$ and $\chih$
are infinite if the $L^\infty$ norms are. Hence we can not hope
to control them unless we can control the $L^\infty$ norms.

\section{Motivation and second version of the Main Theorem}
\subsection{Further discussion of Main Theorem}
Using the null structure and Bianchi  equations
discussed in the previous section we can now give a more
detailed discussion
of  the main
ideas underlining the proof
of our main  theorem. As in the introduction
we write $\HH=\HH_1$ and $\RR_0=\RR[1]$.

{\bf Step 1.}\quad  To start with we recall that
$\trch$ verifies the transport equation
\be{eq:transporttrchi}
\frac{d}{ds}\trch=-|\chih|^2-\f12(\trch)^2.
\end{equation}
Integrating it we easily see that to bound
 $\trch$ it is necessary to have 
a  uniform  bound for $\int_0^s |\chih|^2 $.
 To obtain such an estimate we are allowed to  rely only on 
 the  total  curvature flux  $\RR_0$. The
obvious idea which comes in mind is a trace theorem. 
One  very much wishes that the following holds,
$$(\int_{0}^s|\chih|^2)^\f12\le C \|\nab\chih\|_{L^2(\HH)},$$
with a constant $C$ independent of $0\le s\le 1$  and 
the corresponding null geodesic on which $|\chih|^2$  is integrated .
  If such an estimate
were true we could combine it with   the Codazzi
equation
$$\div\chih =-\b+\f12 \nab\trch+\f12  \trch\c\ze-\ze\c\chih$$
 to estimate $\|\nab \chih\|_{L^2(\HH)}$ in terms of  $\RR_0$ using
only $L^2$ elliptic estimates on the leaves $S_s$. Unfortunately this
 fails to be true. As
explained in the introduction   we  could
overcome this difficulty if  
$\chih=\ddd_LQ+E$ with $Q$ a tensor verifying
 $\|\dddd^2 Q\|_{L^2}+\| Q\|_{L^2}\les
\RR_0$ and $E$ an appropriate error term. It  turns out,
miraculously,  that if we combine the Codazzi equation with the  null Bianchi identities
\eqref{eq:transportrho''},\,
\eqref{eq:transportsi''}:
\beaa
L(\rhoc)+\frac{3}{2}\trch\c\rhoc &=&\div\b-\ze\c\b
+\f12\chih\c(\nab\hot\ze+\f12\trchb\c\chih-\ze\hot\ze).\\
L(\sic)+\frac{3}{2}\trch\c\sic
&=&-\curl\b +\ze\wedge\b
+\f12\chih\wedge(\nab\hot\ze-\ze\hot\ze),
\eeaa
  we can    obtain such
 a decomposition. 
 Indeed, ignoring  for now  the quadratic and
 higher order  terms in the Bianchi identities,
 we  write, 
$$
\dcal\b=L(\rhoc,-\sic)+\ldots
$$
On the other hand 
the Codazzi equation can be written in the form,
$$\dcall\chih=-\b+\f12 \nab\trch+\ldots$$
Thus,
$$\chih=\dcall^{-1}\dcal^{-1}L(\rhoc,-\sic)+\ldots$$
or, ignoring the commutator term $[\dcall^{-1}\dcal^{-1},\ddd_L](\rhoc,-\sic)$,
 we can write,
\be{eq:simpletrace}\chih=L Q+\f12 \dcall^{-1}\c \nab
\trch+\ldots,\qquad Q=\dcall^{-1}\c \dcal^{-1}(\rhoc\,,\, -\sic) 
\end{equation}
Using this expression for $\chih$ we can return to \eqref{eq:transporttrchi}
and write\footnote{In reality we shall estimate
$\trch-\frac{2}{r}$. We choose to ignore this here
 to avoid additional technicalities. }, 
$$|\trch|\les \int_0^s| LQ|^2+\int_0^s|\dcall^{-1}\c \nab\trch|^2
+\int_0^s|\trch|^2+\int_0^s|\mbox{Error}|^2$$
In view of the trace theorem
mentioned above  we can hope to  estimate, for sufficiently small $\RR_0$, 
$$\big(\int_0^s |LQ|^2\big)^\f12\les \|\dddd^2\c\dcall^{-1}\c
\dcal^{-1}(\rhoc,\sic)\|_{L^2(\HH)}+\|(\rhoc,\sic)\|_{L^2(\HH)}\les\RR_0$$ All  other terms, with the exception
of $\int_\ga|\dcall^{-1}\c \nab\trch|^2$, seem  at first glance to be of lower order
 and thus
controllable\footnote{This is in fact  not  true, the commutator
term $[\dcall^{-1}\dcal^{-1},\ddd_L](\rhoc,-\sic)$ turns
out to  have the same level of regularity
as the principal terms.  This fact forces us to 
treat all  error terms very carefully in section 6. }. This
last however leads to serious complications. Indeed $\dcall^{-1}\c \nab$
is a nonlocal operator of order zero and therefore does not
map $L^\infty$ into $L^\infty$. To overcome this difficulty we are forced
to try to prove a stronger estimate for $\trch$. The idea,
as explained in the introduction,  is  to prove
the boundedness of $\trch$ not only in $L^\infty$ but rather in a Besov space
of type $B^1_{2,1}(S_s)$ which imbeds in $L^\infty(S_s)$ and
is stable under the action of  $\dcall^{-1}\c \nab$. 
Alternatively\footnote{ This   second approach is the one
we shall actually  take. } 
 we can estimate $\nab\trch$ in a  Besov norm  of type $B^0_{2,1}(S_s)$.  The
transport character of the equations satisfied  by $\trch$ or $\nab\trch$ leads us to
introduce 
 the following modified Besov
norms
 for  $S$-tangent  tensors $F$, $0\le \th\les 1$:
\beaa
\|F\|_{\BB^\th}&=&\sum_{k\ge 0} 2^{k\th} \sup_{0\le s\le 1}\| P_kF\|_{L^2(S_s)} +
\sup_{0\le s\le 1}\| P_{<0}F\|_{L^2(S_s)},\\
\|F\|_{\PP^\th}&=&\sum_{k\ge 0} 2^{k\th} \| P_kF\|_{L^2(\HH)}+\| P_{<0}F\|_{L^2(\HH)}
\eeaa
with $P_k$ a family of intrinsic  Littlewood-Paley projections
which  are properly defined,  see \cite{KR2}, in terms
of the heat flow on the surfaces $S_s$, see also \cite{S}.  Here\footnote{The low
frequency component of the norms $\BB^\th, \PP^\th$, with $\th\ge 0$,
is  much easier to treat and we shall often ignore it.} $P_{<0}=\sum_{k<0} P_k$.
To show the boundedness 
of $\|\trch\|_{B^1}$ we need:
\begin{enumerate} 
\item   Prove  a stronger trace type theorem of the form,
\be{eq:weaktrace}
\|\int_{0}^s |LQ |^2\|_{\BB^1}\les \|\dddd^2 Q\|_{L^2(H)}^2+ \|Q\|_{L^2(H)}^2\les \RR_0^2
\end{equation}
\item  Show that
$$\|\int_0^s|\dcall^{-1}\c \nab\trch|^2\|_{\BB^1}\les \|\trch\|_{\BB^1}^2 $$

\item Construct the LP-projections $P_k$, intrinsic
 to the leaves of the foliation, with
all desirable properties.
\item 
Conclude, based on a bootstrap argument, that 
$$\|\trch\|_{L^\infty(\HH)}\les \|\trch\|_{\BB^1},\quad
 \|\int_{0}^1|\chih|^2\|_{L^\infty(\HH)}+\|\nab\chih\|_{L^2(\HH)}^2\les\RR_0^2 $$
\item Make sure that all  the error terms are controllable.
It turns out in fact that the commutator $[\dcall^{-1}\dcal^{-1},\ddd_L](\rhoc,-\sic)$
is not in fact lower order, it has the same regularity 
properties as the  main term $\chih$ we have started with.
\end{enumerate}
But this is not all.  The error terms, denoted ``$\mbox{Error}$'' above,   depend not
only on $\trch$ and $\chih$ but also on $\ze$. Moreover
 the renormalized curvature  terms $\rhoc=\rho-\f12 \chih\c\chibh, \,\,
\sic=\si-\f12\chih\wedge\chibh$ depend also on $\chibh$ for which we have not yet any estimates.
As we need to show that $\rhoc$ and $\sic$ are bounded in $L^2(\HH)$ and,
 as we have explained below,  
we expect to be able to bound  the sup norm of $\int_0^1|\chih|^2$,
we infer that we  also need an estimate for $\|\sup_{0\le s\le 1}|\chib(s)|\|_{L^2(S_0)}$.
To obtain such an  estimate  we need  to consider the equations:
\beaa
\frac{d}{ds}\trchb&=&-\f12\trch\trchb-2\div\ze-\chih\c\chibh+2|\ze|^2+2\rho\\
\div\chibh&=&\f12\nab\trchb-\f12\trchb\ze+\ze\c\chibh+\underline{\b}
\eeaa
from which we  can easily convince ourselves 
 that, to control $\|\sup_{0\le s\le 1}|\chib(s)|\|_{L^2(S_0)}$,  we also 
 need an $L^2(\HH)$
estimate for
$\nab\ze$.

{\bf Step 2.}\quad To estimate $\nab \ze$ we use the transport equation
\eqref{eq:newmasstransportgeodesicfirstuse}
 for $\mu$, which we write  in
the form
\bea
\frac{d}{ds} \mu+\frac{3}{2}\trch \mu&=&\chih\c(\nab\hot\ze)+\f12
\trch\rhoc+2\ze\c\nab\trch+\ldots
\label{eq:newmasstransportgeodesicfirstuse}
\eea
together with the Hodge system \eqref{eq:hodgeze} which we write
 in the form,:
\be{eq:hodgezesimple}
\dcal\ze=-(\mu,0)-(\rhoc,-\sic)+\ldots.
\end{equation}
To estimate  $\|\nab\ze\|_{L^2(H)}$ we observe, by $L^2$ estimates,
$$\|\nab\ze\|_{L^2(H)}\les\|\dcal\ze\|_{L^2(H)} \les \|\mu\|_{L^2(\HH)}
+\|\rhoc\|_{L^2}+\|\sic\|_{L^2}
$$
To estimate\footnote{Note that there is no
additional gain of regularity for $\|\mu\|_{L^2(\HH)} $. } $\mu$ we write, integrating
\eqref{eq:newmasstransportgeodesicfirstuse},
\bea
\|\sup_s \mu(s)\|_{L^2(S_s)}&\les&\|\mu(0)\|_{L^2(S_0)}+\|\int_{0}^s \chih\c(\nab
\hot\ze)\|_{L^2(S_s)}\label{eq:transportmusimple}\\
&+ &  \|\int_0^s\ze\c\nab\trch\|_{L^2(S_s)}\nn
+\ldots
\eea
The first integral can be estimated  with the help of,
$$\|\int_{0}^s \chih\c(\nab
\hot\ze)\|_{L^2(S_s)}\les \|(\int_{0}^s |\chih|^2)^\f12\|_{L^\infty(S_s)}\|\nab
\hot\ze\|_{L^2(\HH)}
$$
The second integral can be treated, since we expect to control
 $\|\sup_t|\nab\trch|\|_{L^2(S_0)}$,  provided that we
can control $\|(\int_{0}^s |\ze|^2)^\f12\|_{L^\infty(S_s)}$. Thus to close 
we need to prove for $\ze$ the same type of trace theorem as we 
have discussed  above for $\chih$. Yet, while $\ze$  clearly plays the same role as 
 $\chih$ we do not have a quantity analogous to $\trch$. Indeed the scalar $\mu$ 
is at the same  level with $\nab\trch$ rather than $\trch$.  

The key observation for deriving trace type estimates for $\ze$
is that the  coupled transport-Hodge system  for $(\mu, \ze)$
has a similar structure to that of $(\nab\trch, \chih)$.  Indeed
taking a derivative of \eqref{eq:transporttrchi} we derive,
\be{eq:transportnabtrchi'}
\frac{d}{ds}\nab \trch=-\frac{3}{2} \trch\c\nab \trch-\chih\c \nab\trch-2\chih\c \nab \chih+\ldots
\end{equation}
The $B^1$ estimate for $\trch$ is, at first glance\footnote{For
technical reasons we cannot quite prove this equivalence. The main
difference is that $M$ is a scalar quantity while
$\nab\trch$ is a 1-form.}, 
equivalent to the $B^0$ estimate
for
$\nab\trch$. Using the transport equation for $\nab\trch$ we expect to show
that\footnote{We proceed quite formally here. In fact $\nab\trch$ is not a scalar
 and the integrals on the right hand side don't quite make sense. We only
use this formal  approach here to  make a parallel with the estimates for $(\mu,\ze)$, in
which case $\mu$ is a scalar and all corresponding integrals make sense.},
$$
\|\nab\trch\|_{\BB^0}\les \|\int_0^s \nab\chih\c \chih\|_{\BB^0} + \|\int_0^s \nab\trch\c
\chi\|_{\BB^0}
$$
To control the first integral we need to observe that
\be{eq:decompnabchih}
\chih\in \NN_1(\HH),\qquad \nab\chih = \ddd_L P + \f12 \nab\c \dcall^{-1}\c \nab
\trch + \mbox{Error},
\end{equation}
where $ P\in \NN_1(\HH)$.
This replaces the previous decomposition of $\chih$.
In turn, the trace inequality \eqref{eq:weaktrace} has to be replaced by a
bilinear version\footnote{This formulation is not quite appropriate since,
in this case, $F\c \ddd_LP$ is a tensor and it does not quite make
sense to take its integral, see previous footnote. 
We shall in fact only make  use  of this for scalars.}
\be{eq:strongtrace}
\|\int_0^s F\c \ddd_L P\|_{\BB^0}\les \big(\NN_1(F) +\|F\|_{L_x^\infty L_t^2}\big)\c\NN_1(P)
\end{equation}

The analogy between the systems for $(\mu, \ze)$ and $(\nab\trch, \chih)$ 
suggests that we should
look for 
 a $\BB^0$-Besov type estimate for $\mu$ and a bilinear trace inequality of the type
\eqref{eq:strongtrace} involving $\ze$. To estimate
$\|\int_{\ga_s}\chih\c \nab\hot \zeta\|_{\BB^0}$  we first
 need to express $\nab \ze$   as
we did in \eqref{eq:decompnabchih}
 for $\nab\chih$ . To achieve
this we combine the Hodge system \eqref{eq:hodgezesimple} with 
the null  Bianchi identity\footnote{We shall in fact need 
the renormalized Bianchi identity \eqref{eq:transportmodifiedbbgeod}.}
  connecting  $\bb$ with $(\rho, \sigma)$,
$$
\ddd_L\bb = \dcal^* (\rho, \sigma) + 2 \chibh \c \b + 3(\ze\c\rho -^*\ze\c \si)
$$ 
Thus setting  $P=
{\dcal^*}^{-1}\underline\beta$, with $ P\in \NN_1(\HH)$,
we can write, 
\be{eq:zedecomp}
\nab\ze = \ddd_LP + \nab \c\dcal^{-1} (-\mu, 0) + \ldots.
\end{equation}

Now that we have established the formal similarity between the couples 
$(\nab \trch\,,\,\chih)$ and $(\mu\,,\,\ze)$ we expect to be able
to prove the $B^0$ estimate for $\mu$.
 Indeed, using the transport equation for $\mu$  as in \eqref{eq:transportmusimple}
\footnote{Remark that all the integrands on the right hand side
are now scalars and thus the integrals along the null geodesics make sense.},
$$
\|\mu\|_{\BB^0}\les \|\mu(0)\|_{\BB^0}+\|\int_0^s \chih\c(\nab
\hot\ze)\|_{\BB^0}+   \|\int_0^s\ze\c\nab\trch\|_{\BB^0}
+\ldots
$$
The $B^0$ estimate for the integral 
$\int_0^s \chi\c(\nab
\hot\ze)$ uses the decomposition \eqref{eq:zedecomp}
$$
\|\int_0^s \chi\c(\nab\hot\ze)\|_{\BB^0}\les \|\int_0^s \chih\c \ddd_L P\|_{\BB^0} +
\|\int_0^s \chih\c\,\nab \dcal^{-1}(\mu,0)\|_{\BB^0} +\ldots
$$
The first term can be estimated with the help of the bilinear 
trace inequality \eqref{eq:strongtrace}.
The last  term, since it does not contain an $L$-derivative, can be  treated by
means of the following bilinear  inequality,
$$
\|\int_0^s F\c G\|_{\BB^0} \les \bigg(\NN_1(F) + 
\|\big (\int_0^s |F|^2\big )^{\f12}\|_{L^\infty}\bigg) \|G\|_{\PP^0},
$$
where $\PP^0$ is the integrated Besov norm introduced in \eqref{eq:Penrosenorm}.
This estimate forces us to control the $\PP^0$-norm (which itself is dominated by the 
$\BB^0$-norm) of $\nab\c\dcal^{-1}(\mu,0)$, which  can be estimated by the $\BB^0$-norm 
of $\mu$ and
therefore can be absorbed by a bootstrap argument.

In addition, to treat the error terms
 we are required to control the $\PP^0$-norm of the quantity $\chih\c \chibh$.
We thus need  yet another bilinear estimate of the form,
$$
\|F\c G\|_{\PP^0}\les \bigg(\NN_1(F)  + 
\|\int_0^s |F|^2\big )^{\f12}\|_{L^\infty}\bigg) \|G\|_{\BB^0}
$$
Thus to close our estimates we need a $B^0$-norm estimate for $\chibh$.

{\bf Step 3.}\quad 
The idea is to derive 
a   $B^0$ estimate for $\trchb$ 
from the transport equation,
 $$\frac{d}{ds}\trchb+\f12 \trch \trchb=-2\div\ze-\chih\c\chibh+2|\ze|^2+2\rho.$$
and combine it with
elliptic estimates in $\BB^0$ for
the Hodge system,
$$\div\chibh=\f12\nab\trchb-\f12\trchb\ze+\ze\c\chibh+\underline{\b}.$$
\subsection{ Sobolev and Besov type norms} 
Motivated by the discussion above we now give
precise definitions of the  Besov norms needed in
our work.
 This requires the introduction
of a family of  Littlewood-Paley type projections $(P_k)_{k\in \Bbb Z}\,\,$,
 $\sum_kP_k^2=I$, verifying suitable   properties. Such 
a family was defined in \cite{KR2} with the help
of the heat flow on  the surfaces $S_s$.  With the help
of this family we define, in the usual way, the Besov
norms  $B^\th_{2,1}$ on $S_s$. We also
recall below the definition of  the fractional Sobolev norms.

\begin{definition} Given an arbitrary tensor $F$ on a fixed  $S=S_s$ we  
 define for every $\ga\in{\Bbb R}$,
\be{eq:Sobnorms}
\|F\|_{H^a(S)}= \|\La^a F\|_{L^2(S)}\approx
\big(\sum_{k\ge 0}2^{2ka}\|P_k F\|_{L^2(S)}^2\big)^\f12 + \|P_{<0}F(S)\|_{L^2}
\end{equation}
where $\La =(I-\lap)^\f12$, see  definition in Appendix. 
We also define the Besov norms,
\be{eq:Besov}
\|F\|_{B^a_{2,1}(S)}=\big(\sum_{k\ge 0} 2^{a k}\|P_k F\|_{L^2(S)}\big)+
\|\sum_{k< 0} P_k F\|_{L^2(S)}.
\end{equation}
For $S$-tangent tensors $F$ on $\HH=\HH_1$, $0\le \th\le 1$,  we  also introduce the norms:
\bea
\|F\|_{\BB^\th}&=&\sum_{k\ge 0} 2^{k\th} \sup_{0\le s\le 1}\|
P_kF\|_{L^2(S_s)} +  \sup_{0\le s\le 1}\|
P_{<0}F\|_{L^2(S_s)}    ,\label{eq:Besovnorm}\\
\|F\|_{\PP^\th}&=&\sum_{k\ge 0} 2^{k\th} \| P_kF\|_{L^2(\HH)} +
\| P_{<0}F\|_{L^2(\HH)}\label{eq:Penrosenorm}
\eea
\end{definition}

\subsection{Precise version of the Main Theorem}
We start by making precise assumptions about
the original $2$-surface $S_0\subset \Si$. 
\begin{definition}\label{def:weakly-spherical}Let $S_0\subset\Si$ be a compact 2-surface 
diffeomorphic to ${\bf S^2}$ with $\ga$ a Riemannian metric on it. 
We say that $S_0$ is weakly regular if it can be covered by
 a finite number
of coordinate charts $U$ with
coordinates $\om^1,\om^2$ relative
to which, \be{eq:coordchart}
C^{-1}|\xi|^2\le \ga_{ab}(p)\xi^a\xi^b\le c|\xi|^2, \qquad \mbox{uniformly for  all }
\,\, p\in U
\end{equation}
\be{eq:gammaL2}
\sum_{a,b,c}\int_U|\pr_c\ga_{ab}|^2 dx^1dx^2\le C^{-1}
\end{equation}
 We say that the metric $\ga$ is weakly spherical  if 
in addition
\be{eq:gamma-sph}
\sum_{a,b,c}\int_U|\pr_c(\ga_{ab}-R^2\cga_{ab})|^2 dx^1dx^2\le I_0^2
\end{equation}
where $\cga$ and  is  the standard metric on ${\bf S^2}$ and $R$ 
 is a  constant, $\f12\le R
\le 2$. Here  $I_0$  is a sufficiently small constant.
We refer to surfaces, such as above, as {\bf WS} surfaces.
\end{definition}
 
Motivated by the discussion at the beginning
of this section we are now ready to define 
the initial data quantity 
$\II_0$, 
\bea  \label{eq:initialquantII0}
\II_0&=&I_0+ \|\trch-\frac{2}{r}\|_{L^\infty(S_0)} +
\|\nab \trch\|_{B^0_{2,1}(S_0)}+\|\mu\|_{B^0_{2,1}(S_0)}\\
&+&
\|\trchb-\frac{2}{r}\|_{B^0_{2,1}(S_0)}+\|\chibh\|_{B^0_{2,1}(S_0)}\nn
\eea
and give a precise version
of our main result:
\begin{theorem}[Main Theorem]Consider   an outgoing  null hypersurface $\HH=\HH_1$,
initiating on a compact $2$ surface $S_0\subset\Si$ diffeomorphic to ${\bf S^2}$
  and
foliated by the geodesic foliation associated to the affine parameter $s$
with $s|_{S_0}=0$. Assume  that $r_0=r(0)\ge 1$ and that  both the initial data quantity 
$\II_0$  and  the total curvature flux $\RR_0$  are 
sufficiently small. Then,
\be{eq:I12new}
\|\trch-\frac{2}{r}\|_{L^\infty}\les \II_0+\RR_0
\end{equation}
Moreover,
\beaa
\|\nab \trch\|_{\BB^0}+\|\mu\|_{\BB^0}&\les &\II_0+\RR_0\\
\|\int_0^1|\chih|^2\|_{L^\infty(S_0)}+ \|\int_0^1|\ze|^2\|_{L^\infty(S_0)} &\les &\II_0+\RR_0\\
\|\sup_{0\le s\le 1}|\nab\trch|\,\|_{L^2(S_0)}+\|\sup_{0\le s\le 1}|\mu|\,\|_{L^2(S_0)}&\les
&\II_0+\RR_0\\
\NN_1(\chih)+ \NN_1(\ze)+\NN_1(\trch-\frac{2}{r}) &\les &\II_0+\RR_0\\
\|\sup_{1\le s\le 1}|\trchb+\frac{2}{r}|\|_{L^2(S_0)}+\|\sup_{1\le s\le
1}|\chibh|\|_{L^2(S_0)}&\les &\II_0+\RR_0\\
\|\ddd_L\big(\trchb+\frac{2}{r}\big)\|_{L^2(\HH)}&\les& \II_0+\RR_0\\
\|\ddd_L\chibh\|_{L^2(\HH)}&\les &\II_0+\RR_0\\
\|\trchb+\frac{2}{r}\|_{\BB^0}+\|\chibh\|_{\BB^0}&\les &\II_0+\RR_0
\eeaa
\label{thm:Main}
\end{theorem}
The proof of the theorem relies on an elaborate bootstrap 
argument. In the next subsection we start by making a very mild 
boot-strap assumption and draw some  simple important conclusions
from it. The full set of bootstarp assumptions will be given in the next section.
The idea is to prove that, by choosing the free constants  $\II_0$
and $\RR_0$,  the bootstrap
assumptions can be improved. The proof
of the theorem then follows by a classical continuity
argument.
\subsection{Preliminary bootstrap assumption}
 We assume 
that  there exists a sufficiently  small positive constant $0<\De_0< \f12$
such that 

{\bf BA1.}\qquad \qquad 
$
 \sup_{\HH} r\,|\trchav-\frac{2}{r}|\le \De_0,\qquad   \sup_{\HH} r\,|\trch-\trchav|\le
\De_0.\qquad         
$

We  show  later, at the completion 
of the proof  that  the same inequalities hold true with $\De_0$
replaced by $\frac{\De_0}{2}$.

In view of \eqref{eq:drds} we write,
$$
\frac{dr}{ds}=1+\frac{r}{2}(\overline{\trch}-\frac{r}{2})
$$
According to {\bf BA1.} we infer by integration\footnote{Recall that
$r_0=r(0)\ge 1$.},
\beaa
|r(s)-r_0-s|& \le &\big|\int_0^s \frac{r}{2}(\trchav-\frac{r}{2})\big|\\
&\les&\f12 s\De_0
\eeaa
Thus, for $0\le s\le t\le 1$,
\be{eq:comparisonr}
r_0+(1-\f12 \De_0)s\le r(s)\le r_0+(1+\f12 \De_0)s
\end{equation}
or, simply since $\De_0\le \f12$, $$r_0+\f12 s\le r\le r_0+\frac{3}{2}s.$$
On the other hand,
$$|\trch-\frac{2}{r}|\le \frac{2\De_0}{r}.$$
Hence, for all $0\le s\le t$,
$$\frac{2(1-\De_0)}{r_0+(1+\f12 \De_0 )s}\le \trch\le \frac{2(1+\De_0)}{r_0+(1-\f12 \De_0 )s}
$$
Or, since $\De_0<\f12$,  
\be{eq:trchsupnorm}
 \frac{1}{r_0+2 s} \le \trch\le \frac{3}{r_0+\f12 s}
\end{equation}
Consequently,
$$\log\frac{r_0+2t}{r_0}\le \int_0^s\trch\le 6\log\frac{2r_0+t}{2r_0}$$

As a consequence of the formula
$\frac{d}{ds}\big(\log\sqrt{|\ga|}\big)= 2\trch$ we infer that,
$$v_s=\frac{\sqrt|\ga_s|}{\sqrt|\ga_0|}=\exp(2\int_0^s\trch)$$
Therefore,
$$\big(\frac{r_0+t}{r_0}\big)^2\le \frac{\sqrt|\ga_s|}{\sqrt|\ga_0|}\le 
\big(\frac{2r_0+t}{2r_0}\big)^{12}$$ In particular since $r_0\ge 1$, $0\le t\le 1\le r_0$,
\be{eq:controlvol}
1  \le v_s= \frac{\sqrt|\ga_s|}{\sqrt|\ga_0|}\le 2\big(\frac{3}{2}\big)^6,\qquad \mbox{
for all}\quad 0\le s\le t
\end{equation}
  We infer that the volumes of $S_s$ and $S_0$ remain
comparable for $0\le s\le t$.

As a consequence of  \eqref{eq:controlvol}   we also  infer  that the $L^2(\HH)$ norm,
 defined by the formula \eqref{eq:L2normH} is equivalent to  the product 
 norm on $[0,t]\times S_0$,
\be{eq:L2normHequiv}
\|F\|_{L^2}=\|F\|_{L_t^2L_x^2}=\big(\int_0^1\int_{S_0}|F|^2 \, dA_0 \, ds
\,\big)^\f12=
\bigg(\int_0^1 ds \int_{S_0}|F(s,\om)|^2 \,  \,\sqrt{|\ga_0|} d\om \bigg)^\f12
\end{equation}
We shall also make use of the following norms,
\be{eq:LxLt2norm}
\|F\|_{L_x^\infty L_t^2}=\sup_{\om\in S_0}\big(\int_0^1 ds \,|F(s,\om)|^2\big)^\f12
\end{equation}

\be{eq:Lx2Ltinftynorm}
\|F\|_{L_x^2 L_t^\infty}=\| \,\sup_{0\le s\le 1} |\,F(s,\om)|\,\,\|_{L^2(S_0)}
\end{equation}
as well as, for $1\le p\le \infty$,
\bea
\|F\|_{L_t^\infty L_x^p }&=& \,\sup_{0\le s\le 1}\|\,F(s)\,\|_{L^p(S_0)}
\label{eq:LtinftynormLxp}\\
\|F\|_{L_t^2 L_x^p }&=&\big( \,\int_0^1 \|\,F(s)\,\|_{L^p(S_0)}^2 ds\big)^{\f12}
\eea
Observe that $\|F\|_{L_t^\infty L_x^2 }\le \|F\|_{L_x^2 L_t^\infty}$.

\begin{lemma}
Consider the equation
$\ddd_L F+k\trch F=G$
for S-tangent tensors  $F, G$ on $\HH$. 
Then, for any $1\les p\les \infty$,
\be{eq:LtinftyLx2estimfortransport}
\|F\|_{L_x^pL_t^\infty}\les \|F(0)\|_{L^p(S_0)}+\|F\|_{L_x^pL_t^1}
\end{equation}
\label{le:transportformulaL2}
\end{lemma}
\begin{proof}:\quad Multiplying the transport  equation by $F$
we derive,
\beaa
\frac{d}{ds}|F|^2+2k\trch |F|^2=2G\c F
\eeaa
Applying to this     formula  \eqref{eq:niceintegralform} of
proposition 
\ref{prop:formulatranspgeneral},
$$|F(s,\om)|^2=v(s,\om)^{-2k}\bigg( |F(0,\om)|^2+ \int_0^s 
v(s',\om)^{2k} G\c F(s',\om) ds'  \bigg) $$
 where $v(s,\om)=v_s(\om)=\frac{\sqrt{|\ga_s|}}{\sqrt{|\ga_0|}}$. 
  Using  \eqref{eq:controlvol} we deduce,
\beaa
|F(s,\om)|^2\les |F(0,\om)|^2+\int_0^s|G(s',\om)|\, |F(s',\om)| ds'
\eeaa
or, by a standard argument,
$$|F(s,\om)|\les |F(0,\om)|+\int_0^s|G(s',\om)| ds'$$
We also have, 
\beaa
\sup_{0\le s\le 1}|F(s,\om)|\les |F(0,\om)|+\int_0^1 |G(s',\om)| ds'
\eeaa
Thus the desired  estimate \eqref{eq:LtinftyLx2estimfortransport} 
 follows now by taking the $L^p$
norm with respect to $S_0$.
\end{proof}

\section{Boot-strap assumptions. Preliminary estimates}
Based on the discussion of the previous section 
we make, in addition to {\bf BA1} the  following 
stronger  bootstrap assumptions 
for $\trch, \chih, \ze$

{\bf BA2.}\quad 
\beaa
\|\chih\|_{L_x^\infty L_t^2}\,\,,\,\,\,\, \|\ze\|_{L_x^\infty L_t^2}&\le& \De_0,\\
\|\nab\trch\|_{ L_x^2L_t^\infty}\,\,,\,\|\mu\|_{L_x^2 L_t^\infty}&\le&\De_0\\
\NN_1(\chih)\,\,,
\,\,\, \NN_1(\ze) &\le&\De_0
\eeaa
{\bf BA3.}
\beaa
\|(\trchb+\frac{2}{r})\|_{L_x^2L_t^\infty}, \quad\|\chibh\|_{L_x^2L_t^\infty}&\le& \De_0\\
\eeaa
{\bf BA4.}\qquad \qquad $$\|\trchb+\frac{2}{r}\|_{\BB^0},\quad \|\chibh\|_{\BB^0}\le \De.$$
In the proof of theorem\ref{thm:Main} we shall show 
that, by chosing $\II_0$ and $\RR_0$ small, we can
recover the estimates above with
$\De_0$ replaced by $\De_0/2$.
\begin{remark} Throughout this and the following section
we shall only rely on {\bf BA1}--{\bf BA3}. The assumption {\bf BA4}
will be useful in the last section of the paper.
\end{remark}

To simplify the discussion below it makes sense  
 to  introduce new notations which allow us to
write in a schematic form
the main structure and Bianchi equations, in view of the 
behavior of various components.

\begin{definition} We introduce the following schematic notations for 
the curvature and connection coefficients.
\begin{enumerate}
\item We denote by  $R$ the full collection  of null curvature components $\a, \b,\rho,\si,\bb$.
\item We denote by  $R_0$ the  collection  of null curvature components $ \b,\rho,\si,\bb$.
\item We denote by $\Roc$  the  collection  of the renormalized
 null curvature components
$(\rhoc,-\sic),\bboc$,
 see \eqref{eq:rhocsic}.
\item   We denote
by $\Gd$ the collection $\trch-\frac{2}{r},\chih,\ze$ .
\item We denote by  $\Bd$ the
full collection  $(\trch-\frac{2}{r}),\chih, \ze$  as well as the ``bad components ''
 $\trchb+\frac{2}{r},\,\,  \chib$
\item We denote by $M$ the collection $\nab\trch,\, \mu$.
\item We denote by $\nab \Gd$ the collection
 formed by the  first derivatives $\nab\trch, \nab\chih, \nab\ze$. 
\end{enumerate}
\end{definition}

In view of these notations the bootstrap assumptions {\bf BA1} - {\bf BA4} take the form:
\bea
\|\Gd\|_{L_x^\infty L_t^2}&\le &\De_0,\qquad
\NN_1(\Gd)\le \De_0\\
\|M\|_{L_x^2 L_t^\infty}&\le &\De_0,\qquad
\|\Bd\|_{L_x^2 L_t^\infty}\le \De_0\\
\|\Bd\|_{\BB^0}&\le& \De_0
\eea
We also have,
\be{eq:curvassumption}
\|R\|_{L_t^2L_x^2}\le\RR_0
\end{equation}
and assume that $\RR_0$ is sufficiently small, 
in particular
$\RR_0\le \De_0.$
We record below some simple consequences of the 
bootstrap assumptions {\bf BA1}--{\bf BA3}. 
\subsection{Weakly spherical surfaces} We recall
that the initial surface $S_0$ was weakly spherical {\bf WS}
in the sense of definition \ref{def:weakly-spherical}.
We shall show, using {\bf BA1}--{\bf BA2}, that this property is shared 
by all surfaces $S_s$ of the geodesic foliation of $\HH$. In other words
all surfaces $S_s$ verify {\bf WS}.
\begin{proposition} Consider the  transported local
coordinates $(s,\om)$, defined by the equation $x^\mu=x^\mu(s,\om)$
where $  x^\mu(s,\om)$ are  null geodesics parametrized by
the affine parameter $s$ and $\om$ are local coordinates on
$S_0$. Assume also that on $(S_0,\ga)$ is weakly spherical
in the sense of definition \ref{def:weakly-spherical}.
Then,  on all $\HH$, the metric $\ga$ remains weakly spherical,
i.e., 
\be{eq:weak-regular-cond} 
\|\ga\|_{L_t^\infty L_x^\infty},\,\, \|\ga^{-1}\|_{L_t^\infty L_x^\infty}\,\les 1,\qquad
\|\pr(\ga-(1+s)^2\cga)\|_{L_x^2L_t^\infty}\les \II_0+\De_0
\end{equation}
\label{prop:weak-regular-cond}
\end{proposition}
\begin{proof}:\quad
We have already proved that the metric coefficients $\ga$
remain bounded. To prove the second assertion we recall
the equation \eqref{eq:transportga},
$$\frac{d}{ds}\ga_{ab}=2\chi_{ab}.$$
Since the derivatices $\frac{\partial}{\pr\om^a}$
commute with $\frac{d}{ds}$ and $\chi_{ab}=\f12 \trch\c\ga_{ab}+\chih_{ab}$
relative to the coordinate system  $\om^1,\om^2$,
\beaa
\frac{d}{ds}\pr_c\ga_{ab}& =&2\pr_c\chi_{ab}=\pr_c\big(\trch\c\ga_{ab}+2\chih_{ab}\big)
\\
&=&\nab_c\trch\c\ga_{ab}+\trch \c\pr_c\ga_{ab} +  2\pr_c\chih_{ab}
\\
&=&+\trch \c\pr_c\ga_{ab}+\nab_c\trch\ga_{ab}+2\nab_c\chih_{ab}+(\pr\ga)\c\chih
\eeaa
with $(\pr\ga)\c\chih$ denoting sum of  terms involving
only products between derivatives of the metric coefficients
and components of $\chih$.
Therefore,
\beaa
\frac{d}{ds}\pr_c\ga_{ab}-\trch \c\pr_c\ga_{ab}=\nab_c\chi_{ab}+(\pr\ga)\c\chih
\eeaa
According to  proposition \ref{prop:formulatranspgeneral} 
\be{eq:niceintegralform'}
\pr\ga(s,\om)= v_s\bigg( \pr\ga (0,\om)+ \int_0^s 
v_{s'}^{-1} \,\,g(s',\om) ds  \bigg)
\end{equation}
where $g=\nab\chi+(\pr\ga)\c\chih$, relative to the transported coordinates $\om^a$.
Therefore, denoting  $\cga_s=v_s \cga=v_s\ga(0)$ and observing that the components
of $\pr\cga_s$  are uniformly bounded,
\beaa
\|\pr\ga- \pr\cga_s \|_{ L_x^2 L_t^\infty} &\les &\|\nab\chi\|_{L_x^2L_t^1}
+\|(\pr\ga)\c\chih\|_{L_x^2L_t^1}\\
 &\les&\De_0+\|\pr\ga- \pr\cga_s \|_{ L_x^2 L_t^\infty}
\c  \|\chih\|_{L_x^\infty L_t^1}\\
&\les&\De_0+\De_0\c\|\pr\ga- \pr\cga_s \|_{ L_x^2 L_t^\infty},
\eeaa
whence, in view of  our  weakly spherical assumption on $S_0$  and
  the smallness of $\De_0$,
$$\|\pr(\ga-\cga_s)\|_{L_x^2L_t^\infty}\les \De_0.$$
To end the proof of the theorem it only  remains to remark that,
for $0\le s\le 1$, $|v_s-(1+s)^2|\les \De_0.$
\end{proof}
In particular our surfaces $S_s$ remain weakly regular. This fact allows
us to apply the results of \cite{KR2}. In particular 
we have the following calculus inequalities  for scalar functions $f$
 and tensorfields $F$ on 
each $S=S_s\subset \HH$.
\bea
\|f\|_{L^2(S)}&\les&\|\nab f\|_{L^1(S)}+\|f\|_{L^1(S)}\label{eq:isoperimetric}\\
\|f\|_{L^\infty(S)}&\les&  \|\nab^2 f\|_{L^1(S)}+\|f(S)\|_{L^1}\label{eq:sharp-Sob-LinftyL1}
\eea
\bea
\|F\|_{L^p}&\les &\|\nab F\|_{L^2}^{1-\frac{2}{p}}\|F\|_{L^2}^{\frac 2 p}
+\|F\|_{L^2},\qquad 2\le p<\infty\label{eq:GNirenberg}\\
\|F\|_{L^\infty}&\les& \|\nab ^2 F\|_{L^2}^\f12 \|F\|_{L^2}^\f12+\|F\|_{L^2}
\label{eq:LinftyL2tensor}
\eea
\subsection{  Basic Estimates on $\HH$} The following simple calculus estimates
hold true on $\HH_t$, with any $0<t\le1 $. In applications we 
shall consider only the case
$t=1$.
\begin{lemma} Let $F$ be an arbitrary  $S$ tangent tensorfield on $\HH=\HH_t$. Then,
\be{eq:Lx4ineq}
\|F\|_{L_t^\infty L_x^4}\les  t^{-\f12}\|F\|_{L_t^2L_x^4}+
\|F\|_{L_x^\infty L_t^2}^\f12\c \|\ddd_LF\|_{L_t^2L_x^2}^\f12
\end{equation}
Also,
\be{eq:LtinftyLx4}
\|F\|_{L_t^\infty L_x^4}\les t^{-1/2}\big(\|F\|_{L_t^2L_x^2}+\|\nab F\|_{L_t^2L_x^2}\big)
+t^{1/2}\|\ddd_LF\|_{L_t^2L_x^2}
\end{equation}
\be{eq:SobL6}\|F\|_{L_t^6L_t^6}\les t^{-1/3}\big(\|F\|_{L_t^2L_x^2}+\|\nab
F\|_{L_t^2L_x^2}\big)+t^{2/3}\|\ddd_LF\|_{L_t^2L_x^2}.
\end{equation}
\bea
\|F\|_{L_t^\infty L_x^\infty
}&\les & t^{-\f12}\big(\|F\|_{L_t^2L_x^2}+\|\nab^2F\|_{L_t^2L_x^2}\big)
\label{eq:LxLtinftyineq}\\
&+&t^{\f12}\big(\|\ddd_LF\|_{L_x^\infty L_t^2}+ \|\nab\ddd_LF\|_{L^2}\big)\nn
\eea

\label{le:calculusestimates}
\end{lemma}
\begin{corollary} In the case of interest $t=1$ we have,
see definition \ref{def:NNnorms},
\beaa
 \|F\|_{L_t^\infty L_x^4},\quad \|F\|_{L_t^6 L_x^6} &\les&\NN_1(F)\\
\|F\|_{L_t^\infty L_x^\infty}&\les&\NN_2(F)
\eeaa
\label{corr:calcineq}
\end{corollary}
\begin{remark}
The proof below provides in fact the stronger estimate, $$\|F\|_{
L_x^4L_t^\infty}\les \NN_1(F).$$
\end{remark}
\begin{proof}:\quad Consider $|F|^2$. Integrating its derivative we write
\beaa
|F(s,\om)|^2-|F(0, \om)|^2&=&2\int_0^sF\c\ddd_LF(s',\om) ds'\\
&\les&\big(\int_0^s|F(s',\om)|^2\big)^\f12 \c\big(\int_0^s|\nab_LF(s',\om)|^2\big)^\f12
\eeaa
Therefore, 
\beaa
\|\,\,|F(s)|^2-|F(0)|^2\,\,\|_{L^2(S_0)}^2 &\les&\int_{S_0}
\big(\int_0^s|F(s',\om)|^2\big) \c\big(\int_0^s|\nab_LF(s',\om)|^2\big)\\
&\les&\|F\|_{L_x^\infty L_t^2}^2\c \|\ddd_LF\|_{L^2}^2
\eeaa
It remains to estimate $\|F(0)\|_{L^4(S_0)}$. Let $0\le \theta(s)\le 1 $ be a test  function 
on $[0,t]$ with $\th(0)=1$,  vanishing identically near $t$ and
$\sup_{0\le s\le t}|\th'(s)|\les t^{-1}$. Then,
\beaa
|F(0,\om)|^2&=&-\int_0^t\frac{d}{ds}\big( \th(s)|F(s, \om)|^2\big) ds\\
 &\les&\ t^{-1} \int_0^t |F(s,\om)|^2+\big(\int_0^s|F(s',\om)|^2\big)^\f12
\c\big(\int_0^s|\nab_LF(s',\om)|^2\big)^\f12\\
\eeaa
Thus,  taking 
 $L^2$ norm in $S_0$ and applying Minkowski's inequality:
\beaa
\|F(0)\|_{L^4(S_0)}\les t^{-\f12}\|F\|_{L_t^2L_x^4}+
\|F\|_{L_x^\infty L_t^2}^\f12\c \|\ddd_LF\|_{L^2}^\f12
\eeaa

Similarly, if we denote $|F|^q=(|F|^2)^{q/2}$ for $q\ge 2$,
\beaa
|F(0,\om)|^q&=&-\int_0^t\frac{d}{ds}\big( \th(s)|F(s, \om)|^q\big) ds\\
&\les&t^{-1}\|F\|_{L_t^qL_x^q}^q+
\|\ddd_LF\|_{L_t^2L_x^2}\c\|F\|_{L_t^{2(q-1)}L_x^{2(q-1)}}^{q-1}
\eeaa
Hence, reasoning as above, and by simple interpolation,
\bea
\|F\|_{L_t^\infty L_x^q}&\les &t^{-1/q}\|F\|_{L_t^q L_x^q}+
\|\ddd_LF\|_{L_t^2L_x^2}^{1/q}\c\|F\|_{L_t^{2(q-1)}L_x^{2(q-1)}}^{(q-1)/q}\nn\\
&\les&
\bigg((t^{-1}\|F\|_{L_t^2L_x^2})^{1/q}+\|\ddd_LF\|_{L_t^2L_x^2}^{1/q}\bigg)
\c\|F\|_{L_t^{2(q-1)}L_x^{2(q-1)}}^{(q-1)/q}\label{eq:simpleineqLtinfityLtq}
\eea
which, for $q=4$, becomes
\beaa
\|F\|_{L_t^\infty L_x^4}&\les &\bigg( t^{-1/4}\|F\|_{L_t^2 L_x^2}^{1/4}+
\|\ddd_LF\|_{L_t^2L_x^2}^{1/4}\bigg)\|F\|_{L_t^6L_x^6}^{3/4}
\eeaa
On the other hand we have the classical Sobolev inequality\footnote{See  section
 4.1  of \cite{Kl-Nic}, also chapter 3 of \cite{Chr-Kl}   for a direct  proof of
 this Sobolev inequality
relying on the weak  isoperimetric inequality \eqref{eq:isoperimetric}. } on $\HH$,
$$\|F\|_{L_t^6L_t^6}\les t^{-1/3}\big(\|F\|_{L_t^2L_x^2}+\|\nab
F\|_{L_t^2L_x^2}\big)+t^{2/3}\|\ddd_LF\|_{L_t^2L_x^2}.
$$
Therefore, combining with the previous inequality, 
\beaa
\|F\|_{L_t^\infty L_x^4}\les t^{-1/2}\big(\|F\|_{L_t^2L_x^2}+
\|\nab F\|_{L_t^2L_x^2}\big)
+t^{1/2}\|\ddd_LF\|_{L_t^2L_x^2}
\eeaa
which is precisely \eqref{eq:LtinftyLx4}.
The Sobolev inequality \eqref{eq:LxLtinftyineq} follows 
from \eqref{eq:LtinftyLx4}  combined  with,  
$$\|F\|_{L^4(S)}\les \|F\|_{L^2(S)}^{1/2}\c \big(\|\nab F\|_{L^2(S)}+ \|F\|_{L^2(S)}
\big)^{1/2}.$$ which is a particular case of \eqref{eq:GNirenberg}.
\end{proof}
Combining the results above with the bootstrap assumptions {\bf BA1}--{\bf BA2}
we deduce the following estimates for $\Gd$:
\be{eq:LtpLxqGd}
 \|\Gd\|_{L_t^\infty L_x^4},\quad
\|\Gd\|_{L_t^6 L_x^6}      \les \De_0
\end{equation}
Also for $2\le p<\infty$,
\be{eq:Lt2LxpG}
\|\Gd\|_{L_t^2L_x^p}\les\De_0
\end{equation}

\subsection{Estimates for $\|\Roc\|_{L_t^2L_x^2}$}
Using the   bootstrap assumptions {\bf BA1}--{\bf BA3} we can 
easily deduce the  following:
\begin{proposition} 
\label{prop:easyconsebootstrap}
The renormalized curvature components $\Roc$ verify:
\be{eq:normalizedcurvassumption}
\|\Roc\|_{L_t^2L_x^2}\le\RR_0+\De_0^2
\end{equation}
\end{proposition}
\begin{proof}:\quad To prove \eqref{eq:normalizedcurvassumption} it  suffices to check the
assertion for  $\rhoc=\rho-\f12 \chih\c\chibh$, the other components are similar.
Clearly,
\beaa
\|\rhoc\|_{L_t^2L_x^2}\les\|\rho\|_{L_t^2L_x^2}+\|\chih\|_{L_x^\infty L_t^2}\c
\|\chibh\|_{L_x^2 L_t^\infty}\les \RR_0+\De_0^2
\eeaa
\end{proof}

\subsection{Gauss Curvature}
Recall that the Gauss curvature $K$ of $S_s$
is given by the formula, see \eqref{eq:Gausseq},
$$K=-\frac{1}{4}\trch\trchb+\f12 \chih\c \chibh -\rho.$$
or,  writing  schematically with the help
of the definitions above, $\trch=\frac{2}{r}+\Gd$, $\trchb=-\frac{2}{r}+\Bd$    and,
$$K-\frac{1}{r^2}=\frac{1}{2r} \Gd-\frac{1}{2r} \Bd-\frac{1}{4} \Gd\c \Bd +\f12 \Gd\c \Bd+R.$$
Since the term $\frac{1}{2r} \Gd$  is better than $\frac{1}{2r} \Bd$
we may neglect it and, neglecting also signs and constants on the right hand side, we write
\be{eq:schematicgauss}
K-\frac{1}{r^2}=\frac{1}{2r} \Bd +\Gd\c \Bd+R
\end{equation}
Thus,
 according to our bootstrap assumptions and $r\ge 1$,
\beaa
\|K-\frac{1}{r^2}\|_{L_t^2L_x^2}&\les& \|\Bd\|_{L_t^2L_x^2} +\|\Gd\|_{L_x^\infty L_t^2}
\c\|\Bd\|_{L_x^2 L_t^\infty} +\|R\|_{L_t^2L_x^2}\\
&\les& \De_0+\De_0^2+\RR_0
\eeaa
Thus,
$$
\|K-\frac{1}{r^2}\|_{L_t^2L_x^2}\les \De_0
$$
This proves the  the following
\begin{proposition}
The Gauss curvature $K$ of $S=S_s\subset\HH$ verifies
the following estimates,
\be{eq:L2estimK}
\|K-\frac{1}{r^2}\|_{L_t^2L_x^2}\les \De_0
\end{equation}
\label{prop:Gauss-K1}
\end{proposition}
We shall also need the following important:

\maketitle
\begin{proposition} The following estimate
holds for all $\, \f12<a<1$,
\be{eq:L2estmateLa-gaK}
\|\La^{-a}(K-\frac{1}{r^2})\|_{L_t^\infty L_x^2}\les \De_0
\end{equation}
where $\La =(I-\lap)^\f12.$ 
\label{prop:Gauss-K2}
\end{proposition}
\begin{proof}:\quad Fractional
powers of the Laplacean $\La^a=(1-\De)^{a/2}$ on surfaces $S$ can be defined with 
the help of the heat flow on $S$.  This approach has been developed in \cite{KR2} based
on an assumption of {\sl weak regularity}. We note
that this assumption is implied by the  {\bf WS} property
derived in proposition \ref{prop:weak-regular-cond}.  We refer the
reader to the appendix  for 
a quick review of the main   properties of the heat flow   needed in our  proof of the 
proposition.  
 We shall make   repetead  use of  the composition  properties of   the operators $\La^a$
 $$\La^a\c\La^b=\La^{a+b}$$
Note also that the operators $\La^a$ are bounded in $L^2$ for $a\le 0$. Thus,
$$\|\La^a f\|_{L^2(S)}\les \|\La^b f\|_{L^2(S)}$$
for all $a\le b$.
We shall also make use of the following nonsharp Sobolev inequality:
\begin{lemma}
The following inequality holds for an arbitrary tensorfield $F$
on  a surface $S=S_s$, $0\le s\le 1$, with $2<p<\infty$, $\, s>1-\frac{2}{p}$ :
\be{eq:nonsharp-Sob}
\|F\|_{L^p(S)}\les \|\La^s F\|_{L^2(S)}
\end{equation}
\label{prop:nonsharp-Sob}
\end{lemma}
\begin{proof}:\quad 
The proof is given in \cite{KR2}. It only depends
on the {\sl weak regularity} assumption mentioned above and not on any bounds for $K$.
\end{proof}
We denote by $K_a$,  $\f12 <a <1$,  the constant we intend to bound,
\be{eq:define-Ka} K_a=\|\La^{-a}(K-\frac{1}{r^2})\|_{L_t^\infty L_x^2}.
\end{equation}

 Recall
that, see the derivation of  \eqref{eq:schematicgauss} and the definition
of $\rhoc$,
$$K-\frac{1}{r^2}=\frac{1}{2r}\Bd +\Gd\c\Bd +\rhoc$$
Clearly, on any $S=S_s$,  $$\|\La^{-a}(\frac{1}{2r}\Bd)\|_{L^2(S)}=\frac{1}{2r}
\|\La^{-a}\Bd\|_{L^2(S)}\les \|\Bd\|_{L^2(S)}\les \|\Bd\|_{L_t^\infty L_x^2}\les \De_0$$
On the other hand, to estimate the norm
$J_a(s)=\|\La^{-a}(\Gd\c\Bd)\|_{L^2(S_s)}$, $0\le s\le 1$ we proceed as follows,
\beaa
J_{a}^2(s)&=&\int_s\La^{-2a}(\Gd\c\Bd)\c(\Gd\c\Bd)
\les
\|\La^{-2a}(\Gd\c\Bd)\c \Gd\|_{L^2(S_s)}\c\|\Bd\|_{L^2(S_s)}\\
&\les&\De_0\c\|\La^{-2a}(\Gd\c\Bd)|_{L^4(S_s)}\c\| \Gd\|_{L^4(S_s)}
\eeaa
or, in view of \eqref{eq:LtpLxqGd},
\beaa
J_{a}^2(s)&\les&\De_0^2\c\|\La^{-2a}(\Gd\c\Bd)\|_{L^4(S_s)}
\eeaa
Applying the non-sharp Sobolev inequality \eqref{eq:nonsharp-Sob} with $p=4$
and $\f12<\de\le a$,
$$J_{a}^2(s)\les \De_0^2\c\|\La^{-2a+\de}(\Gd\c\Bd)\|_{L^2(S_s)}\les 
\De_0^2\c\|\La^{-a}(\Gd\c\Bd)\|_{L^2(S_s)}=\De_0^2\c J_a$$
and therefore, since $\De_0$ is small, $J_{a}\les \De_0^2$. 
Therefore,
\be{eq:K4}
\|\La^{-a}(\Gd\c\Bd)\|_{L_t^\infty L_x^2}
\les \De_0^2
\end{equation}
So far we have established the following estimate,
\be{eq:K5}
\|\La^{-a}\big(K-\frac{1}{r^2}-\rhoc\big)\|_{L_t^\infty L_x^2}\les \De_0
\end{equation}
It remains to estimate the more diffficult curvature term
 $\|\La^{-a}\rhoc\|_{L_t^\infty L_x^2}$. To tackle it
we need to integrate in $t$ and bring in a   $\ddd_L$ derivative. 
 Set $$W(s)=\|\La^{-a} \rhoc\|_{ L^2(S_s)}^2-\|\La^{-a} \rhoc\|_{ L^2(S_0)}^2.$$
and proceed as follows:
\beaa
W(s) &=&
\int_S\La^{-a} \rhoc\c \La^{-a}\rhoc-\int_{S_0}\La^{-a} \rhoc\c \La^{-a} \rhoc\\
& \approx&\int_{S_0}\La^{-a} \rhoc\c
\La^{-a}\rhoc(s,\om)-\La^{-a} \rhoc\c
\La^{-a}\rhoc(s,\om) \\
 &=&\int_{S_0}\int_0^s\ddd_L \big(\La^{-a} \rhoc\c \La^{-a}\rhoc(s,\om)\big) ds\\
&=&2 \int_{S_0}\int_0^s\big(\ddd_L (\La^{-a} \rhoc)\c \La^{-a}\rhoc(s,\om)\\
&=&\int_{S_0}\int_0^s  \La^{-2a}(\ddd_L\rhoc(s,\om))\c \rhoc+\int_{S_0}\int_0^s
 [\ddd_L,\La^{-a}]\rhoc\c \La^{-a}\rhoc\\
&\les&\| \La^{-2a}(\ddd_L\rhoc)\|_{L_t^2L_x^2}\c\|\rhoc\|_{L_t^2L_x^2}+
\| [\ddd_L,\La^{-a}]\rhoc\|_{L_t^1L_x^2}\c\|\La^{-a}\rhoc\|_{L_t^\infty L_x^2}
\eeaa
We can also estimate $\|\La^{-a} \rhoc\|_{ L^2(S_0)}^2$ in the same manner 
 using a cut-off
function as in the proof of \eqref{eq:Lx4ineq}.
Thus,
\beaa
\|\La^{-a} \rhoc\|_{L_t^\infty L_x^2}^2\les \|
\La^{-2a}(\ddd_L\rhoc)\|_{L_t^2L_x^2}\c\|\rhoc\|_{L_t^2L_x^2}+
\| [\ddd_L,\La^{-a}]\rhoc\|_{L_t^1L_x^2}\c\|\La^{-a}\rhoc\|_{L_t^\infty L_x^2}
\eeaa
 and therefore,
\be{eq:K6}
\|\La^{-a} \rhoc\|_{L_t^\infty L_x^2}\les \| \La^{-2a}(\ddd_L\rhoc)\|_{L_t^2L_x^2}+\|\rhoc\|_{L_t^2L_x^2}
+\| [\ddd_L,\La^{-a}]\rhoc\|_{L_t^1L_x^2}
\end{equation}

To estimate the commutator
term in \eqref{eq:K6} we need the following:
\begin{lemma} The following estimate holds
true for scalar functions $f$ on $\HH$ with some  snall $\ep>0$:
\be{eq:A2}
\|\,[\ddd_L,
\La^{-a}]\,f\|_{L_t^1L_x^2}\les I_{a}^{\ep}\c\|f\|_{L_t^2L_x^2}
\end{equation}
with  $I_{a}=1+K_a^{\frac 1{1-a}} + K_{a}^\f12$
\label{le:pcomm-postponed}
\end{lemma}
\begin{proof}:\quad
We postpone the proof, which requires
some properties of the heat flow defining 
the fractional operators $\La^{a}$, to  the appendix.  
\end{proof}
In view of the lemma we infer that,
\be{eq:K7}
\|\La^{-a} \rhoc\|_{L_t^\infty L_x^2}\les \| \La^{-2a}(\ddd_L\rhoc)\|_{L_t^2L_x^2}
+(1+I_a^\ep)\c\|\rhoc\|_{L_t^2L_x^2}
\end{equation}

To estimate $\|\La^{-2a}\ddd_L\rhoc\|_{L_t^2 L_x^2}$
we need to use the renormalized\footnote{Note that
the standard Bianchi identity \ref{eq:bianchro} for $\rho$
 would lead to dificculties because of the 
presence of terms like $\Bd\c R$.}  null  Bianchi \eqref{eq:transportrho''}, 
which we write in the following symbolic  form
\beaa
L(\rhoc)+\frac{3}{2}\trch\c\rhoc &=&\div\b-\ze\c\b
+\f12\chih\c(\nab\hot\ze+\f12\trchb\c\chih-\ze\hot\ze)\\
&=&\div\b +\Gd\c( R+\nab\Gd +\Gd\c\Gd)+\f12(\trchb+\frac{2}{r})\c|\chih|^2+\frac{1}{r}\Gd\c\Gd\\
&=&\div\b+\Gd\c( R+\frac{1}{r}\Gd+\nab\Gd +\Gd\c\Bd)
\eeaa
Denote $\bar{R}=( R+\frac{1}{r}\Gd+\nab\Gd +\Gd\c\Bd)$ 
and observe that, in view of {\bf BA1}--{\bf BA3}, it satisfies the estimate,
$$
\|\bar{R}\|_{L_t^2L_x^2}\les \RR_0+\De_0+\|\Gd\c\Bd\|_{L_t^2L_x^2}\les \De_0
$$
Therefore,
\beaa
\|\La^{-2a}\ddd_L\rho\|_{L_t^2 L_x^2}&\les&\|\La^{-2a}\div\b\|_{L_t^2 L_x^2}+
\|\La^{-2a}(\Gd\c \bar{R})\|_{L_t^2 L_x^2}
\eeaa
To show that  $R_{a}=\|\La^{-2a}(\Gd\c \bar{R})\|_{L_t^2 L_x^2}\les\De_0^2$
we use the same trick as we have used above for $J_{a}$. 
\beaa
R_{a}^2&=&\int_{\HH}\La^{-4a}(\Gd\c \bar{R})\c(\Gd\c \bar{R})\les
 \|\La^{-4a}(\Gd\c \bar{R})\|_{L_t^2L_x^4}\c\|\Gd\|_{L_t^\infty L_x^4}\c \|\bar{R}\|_{L_t^2L_x^2}\\
&\les&\De_0^2\c \|\La^{-4a+\de}(\Gd\c \bar{R})\|_{L_t^2L_x^2}
\les\De_0^2\c\|\La^{-2a}(\Gd\c
\bar{R})\|_{L_t^2L_x^2}=
 \De_0^2\c R_{a}
\eeaa
Thus, recalling   \eqref{eq:K7}, 
\beaa 
\|\La^{-a} \rhoc\|_{L_t^\infty L_x^2}&\les& \|\La^{-2a}\ddd_L\rhoc\|_{L_t^2
L_x^2}+(1+I_a^\ep)\c\|\rhoc\|_{L_t^2L_x^2}\nn\\
&\les&\|\La^{-2a}\div\b\|_{L_t^2 L_x^2}+(1+I_a^\ep)\c\De_0
\eeaa
In the next lemma we shall show that
$\|\La^{-2a}\div\b\|_{L_t^2 L_x^2}\les \De_0$. 
Assuming it for a moment we derive,
\be{eq:A5'}\|\La^{-a} \rhoc\|_{L_t^\infty L_x^2}\les (1+I_a^\ep)\c\De_0.
\end{equation}
Thus, together with \eqref{eq:K5},
\beaa
 K_a &\les &\|\La^{-a}\big(K-\frac{1}{r^2}-\rhoc\big)\|_{L_t^\infty L_x^2}
+ \|\La^{-a} \rhoc\|_{L_t^\infty L_x^2} \les  (1+I_a^\ep)\c\De_0\\
&\les &\De_0  + \De_0\c  (1+K_a^{\frac 1{1-a}} + K_{a}^\f12)^{\ep}
\eeaa
Thus, for every $a$, $1/2<a<1$, we can pick $\ep$ sufficiently small 
so that $\ep\c\frac{1}{1-a}<1$ and therefore $K_a\les \De_0$ as desired.
It only remains to prove the following:
\begin{lemma}
\label{le:estimA4} For a 1-form $F$ and $a\ge \f12$ we have
\be{eq:A4}
\|\La^{-2a}\div F\|_{L_t^2 L_x^2}\les\|F\|_{L_t^2L_x^2}
\end{equation}
\end{lemma}
\begin{proof}:\quad
By duality it suffices to prove, for scalar functions $f$ on $\HH$,
\be{eq:A5}
\|\nab\c\La^{-2a}f\|_{L_t^2 L_x^2}\les \|f\|_{L_t^2L_x^2}
\end{equation}
Now, on each $S=S_s$, since $\lap$ commutes with $\La^{-2a}$,
\beaa
\|\nab\La^{-2a}f\|_{L^2(S)}^2&=&\int_S\lap \La^{-2a}f\c\La^{-2a}f
=\int_S\lap \La^{-4a}f\c f\\
&\les&\|\lap \La^{-4a}f\|_{L^2(S)}\c\|f\|_{L^2(S)}\les\|f\|^2_{L^2(S)}
\eeaa
\end{proof}
This ends the proof of proposition \ref{prop:Gauss-K2}.
\end{proof}
The bound on $\|\La^{-a}(K-\frac{1}{r^2})\|_{L_t^2L_x^2} $
allows us to prove the following elliptic estimates:
 
\begin{proposition}
Under the bootstrap assumptions {\bf BA1}-{\bf BA3} we infer that,
\be{eq:K22}
\|\nab^2 f\|_{L^2(S)}+\|\nab f\|_{L^2(S)}\les \|\lap f\|_{L^2(S)}
\end{equation}
Thus, in particular 
\be{eq:K44}
\|\nab^2(\lap^{-1} f)\|_{L^2(S)}+\|\nab^2(\lap^{-1}f)\|_{L^2(S)}\les \|f\|_{L^2(S)}
\end{equation}
\label{prop:Hodge-estim-lap}
\end{proposition}
\begin{proof}
The proof is based on  the following,
\begin{proposition}
The following Bochner identity holds
true on any  2-surface $S=S_s$ and scalar function $f$,
\begin{equation}
\label{eq:scalBoch}
\int_{S} |\nab^{2} f|^{2} = \int_{S} |\lap f|^{2} - 
\int_{S} K |\nab f|^{2}
\end{equation}
\end{proposition}
\begin{proof}:\quad
Recall that on a 2-surface the Riemann tensor
$
R_{abcd} = (\ga_{ac}\ga_{bd}- \ga_{ad}\ga_{bc})K
$
In particular,
$
R_{ab}=\ga_{ab} K.
$
 We start with the identity 
$$
|\nab^{2} f|^{2} = f_{;ab} f^{;ab}= 
 f_{;ab} f^{;ba}= \nab^{a} (f_{;ab} \nab^{b}f) -
\nab^{a} f_{;ab} \nab^{b}f
$$
Furthermore,
$$
\nab^{a} f_{;ab} \nab^{b}f = 
\nab_{b} \lap f \nab^{b} f + R_{ab} \nab^{a} f\nab^{b}f=
\nab_{b} (\lap f \nab^{b} f) -|\lap f|^{2} + K |\nab f|^{2}
$$
Integrating over $S$ and using the Stokes' formula we obtain
\eqref{eq:scalBoch}.
\end{proof}
We now write $K=(K-\frac{1}{r^2})+\frac{1}{r^2}$ in \eqref{eq:scalBoch},
and estimate with the help of proposition \ref{prop:Gauss-K2},
\bea
\int_{S} |\nab^{2} f|^{2}+\frac{1}{r^2}\int_S |\nab f|^2& = &\int_{S} |\lap f|^{2} - 
\int_{S} (K-\frac{1}{r^2}) |\nab f|^{2}\nn\\
&=&\int_{S} |\lap f|^{2}-\int_S\La^{-a}(K-\frac{1}{r^2}) \c \La^{a} |\nab f|^{2}\nn\\
&\les&\int_{S} |\lap f|^{2}+\De_0\c\|\La^{a} |\nab f|^{2}\|_{L^2(S)}\label{eq:K1}
\eea
To estimate the term $\|\La^{a} |\nab f|^{2}\|_{L^2(S)}$ we need 
the following product estimates,
\begin{proposition}The following product estimate holds true
for arbitrary tensors $F, G$ on $S$,
\be{eq:prodga} 
\|\La^{a} (F\c G)\|_{L^2(S)}\les 
\|\La F\|_{L^{2}(S)} \c\|\La^{a} G\|_{L^2(S)}+\|\La^{a} F\|_{L^2(S)} 
\c \|\La G\|_{L^{2}(S)}
\end{equation}
\label{prop:product-estim-S}
\end{proposition}
\begin{proof}:\quad The proof, based on geometric LP theory,
is given in \cite{KR2}. It is interesting  to
note that these product  estimates do not depend on the 
quantity
$\quad K_a=\|\La^{-a}(K-\frac{1}{r^2})\|_{L_t^\infty L_x^2}$. 
\end{proof}
In view of proposition \ref{prop:product-estim-S} we derive:
\beaa
\|\La^{a} |\nab f|^{2}\|_{L^2(S)}&\les&\|\La\nab  f\|_{L^2(S)}\|\La^{a}\nab f\|_{L^2(S)}
\les \|\La\nab 
f\|_{L^2(S)}^2\\&\les&(\,\|\nab^2 f\|_{L^2(S)}+\|\nab f\|_{L^2(S)}\,)^2
\eeaa
Thus, back to \eqref{eq:K1}, we infer  in view of 
 $r\ge 1$ (see\eqref{eq:comparisonr}),
\beaa
\|\nab^2 f\|_{L^2(S)}+\|\nab f\|_{L^2(S)}\les \|\lap f\|_{L^2(S)}
\eeaa
This ends the proof of proposition \ref{prop:Hodge-estim-lap}.
\end{proof}

\subsection{ $L^2$ Estimates for Hodge systems}
We recall the following identities for Hodge systems:
\begin{proposition} Let $(S,\ga)$ be
a compact manifold with Gauss curvature $K$.

{\bf i.)}\quad The following identity holds for vectorfields  $  F  $ 
 on $S$:
\be{eq:hodgeident1}
\int_S\big(|\nab   F  |^2+K|  F  |^2\big)=\int_S\big( |\div   F  |^2+|\curl  F  |^2\big)=\int_S|\dcal  F  |^2
\end{equation} 

{\bf ii.)}\quad The following identity holds for symmetric, traceless,
 2-tensorfields   $  F  $ 
 on $S$:
\be{eq:hodgeident2}
\int_S\big(|\nab   F  |^2+2K|  F  |^2\big)=2\int_S |\div   F  |^2=2\int_S |\dcall   F  |^2
\end{equation} 

{\bf iii.)}\quad The following identity holds for pairs of functions $(\rho,\si)$
 on $S$:
\be{eq:hodgeident3}
\int_S\big(|\nab \rho|^2+|\nab\si|^2\big)=\int_S |-\nab\rho+(\nab\si)^\star |^2=\int_S|\dcalll
(\rho,\si)|^2
\end{equation} 

{\bf iv.)}\quad  The following identity holds for vectors $  F  $ on $S$,
\be{eq:hodgeident3*}
\int_S \big(|\nab  F  |^2-K|  F  |^2\big)=2\int_S|\dcallll   F  |^2
\end{equation}

\label{prop:hodgeident}
\end{proposition}
\begin{proof}:\quad All the above identities follow
easily from the formulas \eqref{eq:dcalident} and \eqref{eq:dcallident}.
Indeed let $\Dcal$ be any of the operators $\dcal$, $\,\dcall$, $\,\dcalll$, $\dcallll$.
Then, for a corresponding tensor $U$,
$$\|\Dcal U\|_{L^2(S)}^2=<^*\Dcal\c \Dcal\, U\,, U>.$$
We give below a direct derivation of the identities 
\eqref{eq:hodgeident3} and  \eqref{eq:hodgeident3*}, the only ones which were not 
    not dicussed in \cite{Chr-Kl}.
 Observe
that $|-\nab\rho+(\nab \si)^\star|^2=|\nab\rho|^2+|\nab\si|^2-2\nab\rho\c(\nab\si)^\star$
 and that, integrating by parts,
 $$\int_S\nab\rho\c(\nab\si)^\star=\int_S\in^{ab}\nab_a\rho\,\nab_b\si
=\int_S\in^{ab}(\nab_b\nab_a\rho)\,\si=0.$$
Recall that $(\dcallll   F  )_{ab}=-\f12 \big(\nab_b  F  _a+\nab_a  F  _b-\div  F  \ga_{ab}\big)$.
Therefore, relative to an  arbitrary ortonormal frame on $S$,
\beaa
4|\dcallll   F  |^2=2|\nab  F  |^2+2|\div  F  |^2+2\sum_{a,b=1,2}\nab_b  F  _a\nab_a  F  _b-4|\div  F  |^2
\eeaa
and therefore,
\beaa
4\int_S|\dcallll   F  |^2&=&2\int_S|\nab  F  |^2-
2\int_S|\div  F  |^2+2\sum_{a,b=1,2}\int_S\nab_b  F  _a\c\nab_a  F  _b\\
&=&2\int_S|\nab  F  |^2 -2\sum_{a,b=1,2}\int_S(\nab_a\nab_b-\nab_b\nab_a)  F  _a \c  F  _b\\
&=&2\int_S\big(|\nab  F  |^2 -
 K|  F  |^2\big)
\eeaa
as desired.
\end{proof}
In view of \eqref{eq:L2estimK} we rewrite \eqref{eq:hodgeident1}, on each $S=S_s\subset\HH$,
in the form,
\beaa
\int_{S}\big(|\nab   F  |^2+\frac{1}{r^2}|  F  |^2\big) &=&2\int_{S}\big( |\div   F  |^2+|\curl  F  |^2\big)-
\int_{S}(K-\frac{1}{r^2})|  F  |^2\\
&=&\int_{S}\big( |\div   F  |^2+|\curl  F  |^2\big)-
\int_{S}\La^{-a}(K-\frac{1}{r^2})\La^{a} |  F  |^2\\
&\les&\|\dcal   F  \|^2_{L^2(S)}+\|\La^{-a}(K-\frac{1}{r^2})\|_{L^2(S)}\c\|\La^{a}|  F 
|^2\|_{L^2(S)}
\eeaa
for some $\f12 <a<1$ sufficiently close to $\f12$. Thus,
in view of the estimate\eqref{eq:L2estmateLa-gaK}
 of proposition \ref{prop:Gauss-K2} and the product estimates of
proposition \ref{prop:product-estim-S},
\beaa
\int_{S}\big(|\nab   F  |^2+\frac{1}{r^2}|  F  |^2\big)&\les&
\|\dcal   F  \|^2_{L^2(S)}+\De_0\c \|\La  F  \|_{L^2(S)}\|\La^a  F  \|_{L^2}\\
&\les&\|\dcal   F  \|^2_{L^2(S)}+\De_0\c \|\La  F  \|_{L^2(S)}^2\\
&\les&\|\dcal   F  \|^2_{L^2(S)}+\De_0\big(\|  F  \|_{L^2(S)}^2+\|\nab  F  \|_{L^2}^2\big)
\eeaa
Since $r\le r_0+\frac{3s}{2}\le 3r_0\les 1 $ we deduce that,
 \beaa
\|  F  \|_{L^2(S)}^2+\|\nab  F  \|_{L^2(S)}^2\les\|\dcal
  F  \|^2_{L^2(S)}+\De_0\big(\|  F  \|_{L^2(S)}^2+\|\nab  F  \|_{L^2}^2\big)
\eeaa
and therefore, for sufficiently small  $\De_0$, 
$$
\|  F  \|_{L^2(S)}^2+\|\nab  F  \|_{L^2}^2\les\|\dcal
  F  \|^2_{L^2(S)}
$$

Proceeding in the same manner, from the
 Hodge identity \eqref{eq:hodgeident2},
and recalling the definition of $\dcall$
 we deduce, for traceless , symmetric 2-tensors $  F  $,
$$
\|  F  \|_{L^2(S)}^2+\|\nab  F  \|_{L^2}^2\les\|\dcall
  F  \|^2_{L^2(S)}
$$

On the other hand, recalling the definition of $\dcalll$ we deduce from
 \eqref{eq:hodgeident3},
$$ \|\nab \rho\|_{L^2(S)}^2+\|\nab\si\|_{L^2(S)}^2=
\int_{\HH}|\dcalll (\rho, \si)|^2$$
We summarize these results in the following 
\begin{proposition}
The following estimates hold on  an arbitrary  2-surface $S=S_s\subset \HH$:

{\bf i.)}\quad The operator $\dcal$  is invertible  on its range and its inverse $\dcal^{-1}$
takes pair of functions $f=(\rho,\si)$  (in the range of $\,\dcal)  $ into  $S$-tangent 1-forms $  F  $
with
 $\div  F  =\rho,\,\, \curl  F  =\si$ with estimate,
\be{eq:estimdcal-1}
\|\nab\c\dcal^{-1}F\|_{L^2(S)}+\|\dcal^{-1}F\|_{L^2(S)}\les\|F\|_{L^2(S)}
\end{equation}

{\bf ii.)}\quad The operator $\dcall$  is invertible on its range  and its inverse $\dcall^{-1}$
takes $S$ tangent 1-forms $F$ (in the range of $\,\dcall$)  into  $S$-tangent symmetric, traceless,
2-tensorfields 
$  Z  $ with
 $\div  Z =F$ with estimate,
\be{eq:estimdcall-1}
\|\nab\c\dcall^{-1}F\|_{L^2(S)}+\|\dcall^{-1}F\|_{L^2(S)}\les\|F\|_{L^2(S)}
\end{equation}

{\bf iii.)}\quad The operator $(-\lap)$ is invertible on its range\footnote{Which consists 
of functions of mean zero.} and its inverse $(-\lap)^{-1}$ verifies the
estimate 
\be{eq:-lap-1}
\|\nab^2(-\lap)^{-1}f\|_{L^2(S)}+\|\nab(-\lap)^{-1}f\|_{L^2(S)}\les \|f\|_{L^2(S)}
\end{equation}

{\bf iv.)}\quad The operator $\dcalll$  is invertible as an operator
defined from pairs of $H^1$ functions with mean zero(i.e. the quotient of $H^1$
by the kernel of $\dcalll$)  and its inverse
$\dcalll^{-1}$
 takes $S$-tangent $L^2$  1-forms $ F$ (i.e. the full range  of $\,\dcalll$) into pair of
functions
$(\rho,\si)$ with mean zero, such that
$-\nab\rho+(\nab\si)^\star=F$, with estimate,
\be{eq:estimdcalll-1}
\|\nab\c \dcalll^{-1}F \|_{L^2(S)}\les\|F\|_{L^2(S)}
\end{equation}

{\bf v.)}\quad The operator $\dcallll$  is invertible as an operator
defined on the quotient of $H^1$-vectorfields by the kernel
of $\dcallll$ . Its inverse $\dcallll^{-1}$
 takes $S$-tangent,  $L^2$,   2-forms $ Z$  (i.e. the full  range of $\,\dcallll$) into $S$
tangent 1-forms
$  F  $ ( orthogonal to the
 kernel of $\dcall$) , such that
$\dcallll   F  =Z$, with estimate,
\be{eq:estimdcallll-1}
\|\nab\c \dcallll^{-1}Z \|_{L^2(S)}\les\|Z\|_{L^2(S)}
\end{equation}
\label{prop:Hodgeestimates}
\end{proposition}
\begin{proof}:\quad 
It  remains to prove\footnote{We have already
proved  {\bf iii} in proposition \ref{prop:Hodge-estim-lap}. }{\bf v)}. 
 The complication in this case is that the curvature term
on the left hand side of \eqref{eq:hodgeident3*} has a negative sign.
This correponds to the nontriviality of the kernel of
$ \dcallll$.
 The proof of \eqref{eq:estimdcallll-1} requires a   contradiction argument.
 As we  will not need the estimate in this paper we omit the argument here.  
\end{proof}

We consider now commutators between $L$ and $\dcal^{-1}$, $\dcall^{-1}$, $\dcalll^{-1}$.
Consider first
 $$[L,\dcal^{-1}](\rho,\si):=\ddd_L\dcal^{-1}(\rho,\si)-\dcal^{-1}( L\rho, L\si).$$
Let $\hat{  F  }=\dcal^{-1}(\rho,\si)$ and $  F  =\dcal^{-1}( L\rho, L\si)$, i.e.
\beaa
L(\div\hat{  F  })&=&L\rho,\qquad\quad L(\curl\hat{  F  })=L\si\\
\div  F  &=&L\rho,\qquad\qquad\,\,\curl  F  =L\si
\eeaa
Therefore,
\beaa
\div(\ddd_L\hat{  F  }-  F  )&=&[\div, L]\hat{  F  }\\
\curl(\ddd_L\hat{  F  }-  F  )&=&[\curl, L]\hat{  F  }
\eeaa
whence, $$\ddd_L\hat{  F  }-  F  =\dcal^{-1}\big([\div, L]\hat{  F  }\,\,,\,\, [\curl, L]\hat{  F  }\big)
=\dcal^{-1}\big([\dcal , L]\hat{  F  }\big)=\dcal^{-1} [\dcal , L]\dcal^{-1}(\rho,\si) $$
or,
$$[L,\dcal^{-1}]=\dcal^{-1} [\dcal , L]\dcal^{-1}$$
Proceeding in precisely same  manner for the other Hodge operators we 
infer that,
\begin{lemma} Let ${\cal D}^{-1}$ be any of the operators $\dcal^{-1}$, $\dcall^{-1}$, $\dcalll^{-1}$.
Then,
\be{eq:commLdcal-1}
[L,{\cal D}^{-1}]={\cal D}^{-1} [{\cal D} , L]{\cal D}^{-1}
\end{equation}
\label{le:commutationL}
\end{lemma}
We shall also need   second derivative estimates. We start with 
$\dcal  F  =(\rho,\si)$. 
In view of  \eqref{eq:dcalident},
$$(-\lap+K)  F  = \dcalll(\dcal   F  )=\dcalll(\rho,\si).$$
We now recall the Bochner identity for 1-forms $  F  $ on $S=S_s$,
$$\int_{S} |\nab^{2}   F  |^{2} = \int_{S} |\lap   F  |^{2} -
\int_{S} K (2|\nab   F  |^{2}-|\div   F  |^{2}) + \int_{S} K^{2} |  F  |^{2}\\
$$
Therefore, proceeding as for the first derivative estimates, 
since $\|\La^{-a}( K-\frac{1}{r^2}) \|_{L^2(S)}\les \De_0$, and
with the help
of the product estimates of proposition \ref{prop:product-estim-S}, 
\beaa
\big|\int_{S}( K-\frac{1}{r^2}) (2|\nab   F  |^{2}-|\div   F  |^{2})\big|
&\les& \De_0\c
\|\La^a\big(2|\nab   F  |^{2}-|\div   F  |^{2}\big)\|_{L^2(S)}\\
&\les&\De_0 \big(\|\nab^2  F  \|_{L^2(S)}^2+\|\nab  F  \|_{L^2(S)}^2\big)
\eeaa
we derive,
\beaa
\|\nab^2  F  \|_{L^2(S)}^2+\|\nab  F  \|_{L^2(S)}&\les&\|\dcalll (\rho,\si)\|_{L^2(S)}^2
+\De_0 \big(\|\nab^2  F  \|_{L^2(S)}^2+\|\nab  F  \|_{L^2(S)}^2\big)\\& +&
\|  F  \|_{L^\infty(S)}^2\| K\|_{L^2(S)}^2
\eeaa
and therefore, for sufficiently small $\De_0$,
$$\|\nab^2  F  \|_{L^2(S)}^2+\|\nab  F  \|_{L^2(S)}\les 
\|\dcalll (\rho,\si)\|_{L^2(S)}^2+\|  F  \|_{L^\infty(S)}^2\| K\|_{L^2(S)}^2
$$
Thus integrating over  $s$,
\beaa
\|\nab^2  F  \|_{L_t^2L_x^2}^2+\|\nab  F  \|_{L_t^2L_x^2}^2
&\les&\|\dcalll (\rho,\si)\|_{L_t^2L_x^2}^2
+\|  F  \|_{L^\infty}^2\| K\|_{L_t^2L_x^2}^2\\
&\les&\|\dcalll (\rho,\si)\|_{L_t^2L_x^2}^2 +\De_0^2\|  F  \|_{L^\infty}^2
\eeaa
Taking
into account the estimate, see \eqref{eq:LxLtinftyineq},
$$\|  F  \|_{L^\infty}^2\les \NN_2(F)^2\les \|  F  \|_{L_t^2L_x^2}^2+\|\nab^2  F 
\|_{L_t^2L_x^2}^2 +\|\nab_L  F  \|_{L_t^2L_x^2}^2+\|\nab\nab_L  F 
\|_{L_t^2L_x^2}^2$$ we infer that,
\beaa \|\nab^2  F  \|_{L_t^2L_x^2}^2+\|\nab  F  \|_{L_t^2L_x^2}^2
\les\|\dcalll (\rho,\si)\|_{L_t^2L_x^2}^2
+\De_0^2\big\|\nab_L  F  \|_{L_t^2L_x^2}^2+\|\nab\nab_L  F  \|_{L_t^2L_x^2}^2\big)
\eeaa
or, using  the definition of $\dcal^{-1}$,
$$\|\nab^2\dcal^{-1}(\rho,\si)\les\|\dcalll (\rho,\si)\|_{L_t^2L_x^2}
+\De_0\big(
\|\nab_L  F  \|_{L_t^2L_x^2}+\|\nab\nab_L  F  \|_{L_t^2L_x^2}\big)
$$
Similar estimates can be derived for $\dcall^{-1}$ and $\dcalll^{-1}$.
\begin{proposition} Let $\Dcal$ denote either $\,\dcal$  or $\dcall$
with $\,^\star\Dcal$ and $\,\Dcal^{-1}$ the corresponding adjoints and
inverses. Then, for $F=(\rho, \si)$ pairs of scalar functions in the first case 
 and  $S$ -tangent one forms for the second case,
the  following second order
estimates hold true:
\bea
\|\nab^2\Dcal^{-1} F\|_{L_t^2L_x^2}&\les&\|^\star\Dcal F\|_{L_t^2L_x^2}+
\De_0\big(
\|\nab_L\Dcal^{-1} F\|_{L_t^2L_x^2}+\|\nab\nab_L\Dcal^{-1} F\|_{L_t^2L_x^2}\big)\nn
\eea
\label{prop:secondorderHodge}
\end{proposition}

We shall also need the following :
\begin{proposition} Let $\Dcal^{-1}$ denote either $\dcal^{-1}$ or $\dcall^{-1}$
 Then, for $F=(\rho, \si)$ pairs of scalar functions in the first case 
 and  $S$ -tangent one forms for the second case. For any $1<p\le 2$,
\be{eq:Lt2LxpDcal-1}
\|\Dcal^{-1} F\|_{L^2(S)}\les \|F\|_{L^p(S)}
\end{equation}
\label{prop:Lt2LxpDcal-1}
\end{proposition}
\begin{proof}:\quad By duality it suffices to
prove the estimate
$$\|^\star\Dcal^{-1} F\|_{L^q(S)}\les \|F\|_{L^2(S)}$$
for $q^{-1}+p^{-1}=1$. In view of \eqref{eq:nonsharp-Sob}
of proposition \ref{prop:nonsharp-Sob},
\beaa
\|^\star\Dcal^{-1}F\|_{L^q(S)}&\les& 
\|\La ^\star\Dcal^{-1}F\|_{L^2(S)}\\
&\les& \| ^\star \Dcal^{-1}F\|_{L^2(S)}+\| \nab^\star
\Dcal^{-1}F\|_{L^2(S)}\les \|F\|_{L^2(S)}
\eeaa
in view of proposition \ref{prop:hodgeident}
\end{proof}

\section{ Main Estimates in Besov Norms} We present in this
section
our main lemmas
concerning the Besov norms $\BB^0$ and $\PP^0$. These
estimates, based on the properties
of the LP projections $P_k$ and the bootstrap assumptions  ${\bf BA1}-{\bf
BA2}$,  are quite involved and require
 a separate paper for their proof,
 see \cite{KR2}, \cite{KR3}.

We start with  two propositions concerning the  Besov
spaces $B^\ga_{2,1}(S)$  which were proved in \cite{KR2}.
We remark that the estimates of \cite{KR2} were proved based only on the assumption
of weak regularity. This assumption holds for our surfaces $S_s$ in view of 
proposition \ref{prop:weak-regular-cond}.

\begin{proposition}\label{prop:estimB1-nabB0}
For scalar functions  $f$ on $S$
the following sharp Sobolev inequality holds,
\be{eq:Besov-Sobolev}
\|f\|_{L^\infty}\les \|f\|_{B^1_{2,1}(S)}
\end{equation}
Morever, 
\be{eq:estimB1-nabB0}
\|f\|_{B^1_{2,1}(S)}\les \|f\|_{L^2(S)}+\|\nab f\|_{B^0_{2,1}(S)}
\end{equation} 
As a corollary we have for scalar functions on  $\HH$,
\be{eq:Besov-Sobolev-BB0}
\|f\|_{L_t^\infty L_x^\infty}\les \|f\|_{\BB^1}\les 
\|f\|_{L_t^\infty L_x^2}+\|\nab f\|_{\BB^0}
\end{equation}
\end{proposition}
We shall need the following  sharp
version of product estimates. 
\begin{proposition} Let $f$ be a scalar and  $U$ an  arbitrary tensor 
   on a fixed  $S=S_s$.
 The following bilinear estimate holds true:
\be{eq:firstbilBesov}
\| f\c U\|_{B^0_{2,1}(S)}\les\big( \|\nab
f\|_{L^2(S)}+\|f\|_{L^\infty(S)}\big)\|U\|_{B^0_{2,1}(S)}
\end{equation}
Also,  on  all $\HH$,
\be{eq:secondbilBesov}
\| f\c U\|_{\BB^0}\les \big(\|\nab f\|_{L_t^\infty L_x^2} +
\|f\|_{L^\infty}\big)\|U\|_{\BB^0}
\end{equation}
\label{le:firstbilBesov}
\end{proposition}
We shall need the following non sharp
embedding estimate,
\begin{proposition}
Let $F$ be as $S$-tangent tensorfield on $\HH$.
The following non sharp embedding result holds
 for any $0\le \th<\f12$:
\be{eq:GoodestimBthetanotGd}
\|F\|_{\BB^\th}\les \NN_1(F)
\end{equation}
As a consequence 
 we infer the following estimate for $\Gd$:
\be{eq:GoodestimBtheta}
\|\Gd\|_{\BB^\th}\les \De_0
\end{equation}
\label{prop:easyconsebootstrap2}
\end{proposition}
\begin{proof}:\quad
See \cite{KR3}.
\end{proof}

\subsection{Main Lemmas}

\begin{lemma}[Main Lemma] Consider the transport equation along $\HH$ for
 an arbitrary $S$-tangent  tensor  $U$:
\be{eq:transportmainle}
\ddd_L U=F
\end{equation}
with $F$ given tensor of same type. Assume $U(0)=0$ and that $F$
is one of the following types of bilinear expressions:

{\bf i.} \quad $F=G\c\ddd_L P$. Then,
\be{eq:mainlemma1}
\|U\|_{\BB^0}\les \big(\NN_1(G)+\|
G\|_{L_x^\infty L_t^2}\big)\c\NN_1(P)
\end{equation}

{\bf ii.} \quad $F=G\c W$. Then,
\be{eq:mainlemma2}
\qquad\qquad\,\,\|U\|_{\BB^0}\les \big(\NN_1(G)+\|
G\|_{L_x^\infty L_t^2}\big)\c\|W\|_{\PP^0}\qquad\qquad
\end{equation}

{\bf iii.}\quad If $\, U(0)\neq 0$ and $F=0$,
\be{eq:homogtrbesov}
\|U\|_{\BB^0}\les \|U(0)\|_{B^0_{2,1}(S_0)}
\end{equation}

{\bf iv.}\quad 
If the tensorfield  $\,\, W$ verifies the transport equation,
\be{eq:transport2mainle}
\ddd_LW=E
\end{equation}
 the following product estimate holds true,
\be{eq:mainlemma4}
\|G\c W\|_{\PP^0}\les \big(\NN_1(G)+\|
G\|_{L_x^\infty L_t^2}\big)\c\big(\|W(0)\|_{B^0_{2,1}(S_0)}+\|E\|_{\PP^0}\big)
\end{equation}
\label{le:mainlemma}
\end{lemma}

We shall also make use of the following nonsharp  product estimates:
\begin{proposition} For any $S$-tangent tensors $F,G$ we have,
\label{prop:5.7}
 for $0\le\ep<1/2$,
\bea
\|F\c G\|_{\PP^\ep} &\les& \big(\,\|F\|_{L_t^2 L_x^2}+
\|\nab F\|_{L_t^2 L_x^2}\big)\c\|G\|_{\BB^\ep}\\
\| F\c G\|_{\PP^\ep} &\les& \NN_2(F)\c\|G\|_{\PP^\ep}\label{eq:bilin2newnew}\\
\|F\c G\|_{\PP^\ep}&\les& \NN_1(F)\c
\big(\|G\|_{L_t^2 L_x^2}+
\|\nab G\|_{L_t^2 L_x^2}\big)\label{eq:anotherbilPP0}
\eea
\end{proposition}

As a consequence of the main  lemma we 
prove the following:

\begin{proposition}
Assume that the scalar function $U$  satisfies the transport
 equation along  $\HH$, for some positive number $k$.
$$\frac{d}{ds}U+k\trch U =F_1\c\ddd_L P+F_2\c W.$$
Then,
\bea
\|U\|_{\BB^0}&\les&\|U(0)\|_{B^0_{2,1}(S_0)} +
\big(\NN_1(F_1)+\|F_1\|_{L_x^\infty L_t^2}\big)\c\NN_1(P)
\label{eq:mainleprop1}\\
&+&\big(\NN_1(F_2)+\|F_2\|_{L_x^\infty L_t^2}\big)\c\|W\|_{\PP^0}\nn
\eea
\label{prop:mainlemmaapplication}
\end{proposition}
\begin{remark} We shall see from the proof below that the precise
form
of the left hand side of the  transport equations above does not matter;
 we can replace $k\trch$ by any scalar  function $\kappa$ such that $\int_0^s\kappa$
has the same properties of the integral factor $v_s$ below.
\label{rem:mainlemmaapplication}
\end{remark}
\begin{proof}:\quad Let
$v_s(\om)=v(s,\om)=\frac{\sqrt{|\ga_s(\om)|}}{\sqrt{|\ga_0(\om)|}}$. According
to  formula  \eqref{eq:integratingfactor},
$$\ddd_L(v^k U)=v^k\big(F_1\c\ddd_L P+F_2\c H\big)
$$
or, setting $U'=v^k U$, $F'_1=v^k F_1$, $F_2'=v^k F_2$,
$$\ddd_L U'=F'\ddd_L P+G'\c H.$$
In view of main lemma \ref{le:mainlemma},
\beaa
\|U'\|_{\BB^0}&\les& \|U'(0)\|_{B^0_{2,1}(S_0)}+
\big(   \NN_1(F_1')+\|F'_1\|_{L_x^\infty L_t^2}\big)\c\NN_1(P)\\
 &+&\big(\NN_1(F_2')+\|F_2'\|_{L_x^\infty L_t^2}\big)\|H\|_{\PP^0}
\eeaa
Now, recall that 
$$\frac{d}{ds} v=\trch\,  v, \qquad v(0)=1$$
We have already established a pointwise  bound for $v$, see \eqref{eq:controlvol}.
We can also provide a bound for $\|\nab v\|_{L_x^2L_t^\infty}$.
Indeed,
 with the help of the commutation 
formula \ref{prop:commutationL},
$\ddd_L\nab v-\nab \frac{d}{ds} v=-\chi\c \nab v$.
Henceforth,
$$\ddd_L\nab v=\nab(\trch\, v)-\chi\c 
\nab v=\f12 \trch \nab v -\chih\c \nab v+v\nab\trch.$$ Therefore,
$\frac{d}{ds}|\nab v|^2=\trch|\nab v|^2-\chih\c\nab v\c\nab v.$
We can now apply  formula \eqref{eq:niceintegralform}
to derive,
$$|\nab v(s,\om)|^2=v(s,\om)\bigg( |\nab v(0)|^2+\int_0^s v(s', \om)^{-1}h(s',\om) \bigg)$$
where $h(s',\om)=-\chih\c \nab v+v\nab\trch$. In view of  $\nab v(0)=0$ and 
our bootstrap
assumptions {\bf BA1} and {\bf BA2} we deduce,
\beaa \|\sup_{0\le s\le t}|\nab v|\|_{L^2(S_0)}&\les& \|h\|_{L_x^2L_t^1}\les
\|\chih\|_{L_x^\infty L_t^1}\c\|\nab v\|_{L_t^\infty L_x^2}+\|\nab\trch\|_{L_x^2L_t^1}\\
&\les&\De_0\|\nab v\|_{ L_x^2 L_t^\infty} +\De_0
\eeaa
and consequently, since $\De_0$ is sufficiently small,
\be{eq:L2estimnabv}\|\nab v\|_{ L_x^2L_t^\infty}\le \De_0.
\end{equation}
Using these bounds for $v$ and $\nab v$ we infer the following,
\beaa
\|F'_1\|_{L_x^\infty L_t^2}&\les&     \|F_1\|_{L_x^\infty L_t^2}   \\
\|\nab F'_1\|_{ L_t^2L_x^2}&\les&\|\nab F_1\|_{ L_t^2L_x^2}+
\|\nab v\|_{ L_x^2L_t^\infty}\|F_1\|_{L_x^\infty L_t^2}\\
&\les&\|\nab F_1\|_{ L_t^2L_x^2}+\De_0\|F_1\|_{L_x^\infty L_t^2}\\
\|\ddd_L F'_1\|_{ L_t^2L_x^2}&\les&\|\ddd_L F_1\|_{ L_t^2L_x^2}+\De_0\|F\|_{L_x^\infty L_t^2}
\eeaa
and similar for $F_2'$.
Thus,
\beaa
\|U'\|_{\BB^0}&\les&\|U'(0)\|_{B^0_{2,1}(S_0)} +
\big(\NN_1(F_1)+\|F_1\|_{L_x^\infty L_t^2}\big)\c\NN_1(P)\\
&+&\big(\NN_1(F_2)+\|F_2\|_{L_x^\infty L_t^2}\big)\|H\|_{\PP^0}
\eeaa
To finish it remains to check
that $\|U\|_{\BB^0}\les \|U'\|_{\BB^0}$ and
 $\|U'(0)\|_{B^0_{2,1}(S_0)}\les \|U(0)\|_{B^0_{2,1}(S_0)}$.
These follow easily from from the pointwise bounds on $v$ and $v^{-1}$,
the $L^2$ bound on $\nab v$ and  a straightforward application
of lemma \ref{le:firstbilBesov}.
\end{proof}
The following is a also an easy consequence of  part iv) of the main lemma
\begin{proposition}
 Assume that the $S$ tangent  tensorfield $W$ verifies the transport equation,
$$\ddd_L W+k\trch W= E$$
for some positive $k$.
Then for any other $S$ tangent tensorfield $F$,
\be{eq:mainlemma-application2}
\|G\c W\|_{\PP^0}\les \big(\NN_1(G)+\|
G\|_{L_x^\infty L_t^2}\big)\c\big(\|W(0)\|_{B^0_{2,1}(S_0)}+\|E\|_{\PP^0}\big)
\end{equation}
\label{prop:mainlemmaapplication2}
\end{proposition}
\begin{proof}:\quad
Can be proved from part iv) of the main lemma
by following the same steps as in the proof of proposition
\ref{prop:mainlemmaapplication}.
\end{proof}
As a corollary 
to proposition \ref{prop:mainlemmaapplication}
we can prove the following version
of the sharp
classical trace theorem.
\begin{corollary} Assume $F$ is an $S$-tangent tensor which
admits a decomposition
of the form, $\nab F=\ddd_LP+E$, with $P,E$ tensors of the same type.
Then,
\be{eq:funny-classical-trace}
\|F\|_{L_x^\infty L_t^2}\les \NN_1(F)+\NN_1(P)+\|E\|_{\PP^0}
\end{equation}
\label{corr:funny-classical-trace}
\end{corollary}
\begin{proof}:\quad
The scalar function  $f(t)=\int_0^t|F|^2$ verifies the transport
equation,
$$\ddd_L f=|F|^2,\qquad U(0)=0$$
Differentiating and applying the commutator formula
\eqref{eq:commscalarf} we derive,
\beaa
\ddd_L (\nab f)+\f12 \trch (\nab f)&=&2 F \c\nab F-\chih \c(\nab f)\\
&\approx&F\c \ddd_LP +F\c E+\Gd \c(\nab f)
\eeaa

Applying  \eqref{eq:mainleprop1} we deduce,
 \beaa
\|\nab f\|_{\BB^0}&\les&
\big(\NN_1(F)+\|F\|_{L_x^\infty L_t^2}\big)\c\NN_1(P)
\\
&+&\big(\NN_1(F)+\|F\|_{L_x^\infty L_t^2}\big)\c\|E\|_{\PP^0}\nn\\
&+&\big(\NN_1(\Gd)+\|\Gd\|_{L_x^\infty L_t^2}\big)\c\|\nab f\|_{\PP^0}
\eeaa
Therefore in view of our bootstrap assumptions for $\Gd$, 
\beaa
\|\nab f\|_{\BB^0}&\les&
\big(\NN_1(F)+\|F\|_{L_x^\infty L_t^2}\big)\c\big(\NN_1(P)+\|E\|_{\PP^0}\big)
+\De_0\|\nab f\|_{\BB^0}
\eeaa
which implies, for $\De_0$ sufficiently small.
\be{eq:interm-classical-trace}
\|\nab f\|_{\BB^0}\les
\big(\NN_1(F)+\|F\|_{L_x^\infty L_t^2}\big)\c\big(\NN_1(P)+\|E\|_{\PP^0}\big).
\end{equation}
Now, in view of  estimate \eqref{eq:Besov-Sobolev-BB0},
 we infer that,
\beaa
\|f\|_{L_t^\infty L_x^\infty}&\les& \|f\|_{\BB^1}\les 
\|f\|_{L_t^\infty L_x^2}+ \|f\|_{\BB^1}\\
&\les&\|f\|_{L_t^\infty L_x^2}+\big(\NN_1(F)+\|F\|_{L_x^\infty
L_t^2}\big)\c\big(\NN_1(P)+\|E\|_{\PP^0}\big)\\
&\les&\|F\|_{L_t^2L_x^4}^2+\big(\NN_1(F)+\|F\|_{L_x^\infty
L_t^2}\big)\c\big(\NN_1(P)+\|E\|_{\PP^0}\big)
\eeaa
Thus, recalling the definition of $f=\int_0^t |F|^2$,
and the estimate  $\|F\|_{L_t^2L_x^4}\les \NN_1(F)$  of
 corollary \ref{corr:calcineq},
 we infer that, 
\beaa
\|F\|_{L_x^\infty L_t^2}^2\les \big(\NN_1(F)+\|F\|_{L_x^\infty
L_t^2}\big)\c\big(\NN_1(P)+\|E\|_{\PP^0}\big)+\NN_1(F)^2
\eeaa
whence,
\beaa
\|F\|_{L_x^\infty L_t^2}\les \NN_1(F)+\NN_1(P)+\|E\|_{\PP^0}
\eeaa
as desired.

\end{proof}

We shall also need elliptic estimates in Besov norms applied to the Hodge operators
$\dcal,\dcall, \dcalll, \dcallll$. More precisely we want to extend the estimates
of proposition \ref{prop:Hodgeestimates} from $L^2(S)$  to $\PP^0$. The proofs, 
 which are quite subtle because of the commutators between the LP projections $P_k$
and  the Hodge operators,
are to be found in \cite{KR2}, \cite{KR3}.

\begin{proposition} Let $\Dcal^{-1}$ denote  one of the operators
$ \dcal^{-1}$, $\dcall^{-1}$ or  $\dcalll^{-1},  \dcallll^{-1}$.
  Then, for a
corresponding  $S$-tangent tensor
$F$ on
$\HH$, and any $0\le\th<1$,
\be{eq:funnyhodgeestimPP}
\|\nab\c\Dcal^{-1} F\|_{\PP^\th}\les \|F\|_{\PP^\th}
\end{equation}
\label{prop:funnyhodgeestimPP}
\end{proposition}
We shall also need the following
version of proposition \ref{prop:Lt2LxpDcal-1}.
\begin{proposition} Let $\Dcal^{-1}$ denote either $\dcal^{-1}$,
 $\dcalll^{-1} $,  $\dcall^{-1}$. 
 Then, for appropriate $S$-tangent tensors on $\HH$  and   any $1<p\le 2$,
\be{eq:Lt2LxpDcal-1besov}
\|\Dcal^{-1} F\|_{B^0_{2,1}(S)}\les \|F\|_{L^p(S)}
\end{equation}
Also, along $\HH$,
\be{eq:Lt2LxpDcal-1BB}
\|\Dcal^{-1} F\|_{\BB^0}\les \|F\|_{L_t^\infty L_x^p}.
\end{equation}
Moreover, for $0\le\th<\f12$ and $\frac{2}{2-\th}<p\le 2 $,
\be{eq:Lt2LxpDcal-1PP}
\|\Dcal^{-1} F\|_{\PP^\th}\les \|F\|_{L_t^2 L_x^p}
\end{equation}
The same estimates hold true
if we replace $\Dcal^{-1}$ by $\nab\c\Dcal^{-2}$
where $\Dcal^{-2}$ denotes $\dcall^{-1}\c\dcal^{-1}$ or $\dcal^{-1}\c\dcalll^{-1}$.
\label{prop:SobolevBesovLp}
\end{proposition}
We shall also make
use of the following:
\begin{proposition} Let $F$ be an $S$-tangent  tensor on $\HH$. We have the estimate,
\be{eq:estimBBL2H}
\|F\|_{\BB^0}\les\NN_1(F) 
\end{equation}
Thus, in particular,
 \be{eq:estimDcal-1BBL2H}
\|\Dcal^{-1} F\|_{\BB^0}\les \|F\|_{L_t^2L_x^2}+\|\nab_L\c\Dcal^{-1} F\|_{L_t^2L_x^2}
\end{equation}
where $\Dcal^{-1}$ is one of the inverse Hodge operators $\dcal^{-1}, \dcall^{-1}, 
\dcalll^{-1}$, $\dcallll^{-1}$ encountered above.  
 The same  estimates hold true
if we replace $\Dcal^{-1}$ by $\nab\c\Dcal^{-2}$
where $\Dcal^{-2}$ denotes $\dcall^{-1}\c\dcal^{-1}$ or $\dcal^{-1}\c\dcalll^{-1}$.
\label{prop:estimDcal-1BBL2H}
\end{proposition}

\section{ Structure of error terms}
In this section we shall study various type of error terms
which we will have to deal with in the next section. Some
of these error terms are easy to treat while others, such
as various commutators between  $\ddd_L$ and  operators of order $-1$,  are of the same
 order of regularity  as the principal terms and therefore require 
a lot of care.

The
main results of the section are  proposition \ref{prop:main-comm-result}   
 and the decomposition lemma \ref{le:decompbandmore}, both 
will play an important  role in the next section.

We recall our main   bootstrap assumptions: 
\bea
\|\Gd\|_{L_x^\infty L_t^2},\quad \|M\|_{L_x^2 L_t^\infty},\quad \NN_1(\Bd)&\le &\De_0,\qquad
\\
\|\Bd\|_{L_x^2 L_t^\infty}, \quad \|\Bd\|_{\BB^0}&\le& \De_0
\eea
In view of equations \eqref{eq:LtpLxqGd}, \eqref{eq:Lt2LxpG} as well as 
propositions \ref{prop:easyconsebootstrap}, \ref{prop:easyconsebootstrap2}
we also have,
\beaa
 \|\Gd\|_{L_t^\infty L_x^4},\quad 
 \,\,\,\|\Gd\|_{L_t^6 L_x^6} 
&\les& \De_0\\
\|\Gd\|_{L_t^2L_x^p}, \qquad\  \|\Gd\|_{\PP^{\th}}   &\les&\De_0 ,\quad \qquad2\le
p<\infty,\quad
\th< 1\\
\|\Gd\|_{\BB^\th}&\les& \De_0,\qquad\qquad\quad\,\, 0\le \th<1/2.
\eeaa
Also, with a constant $\RR_0\le \De_0$,
$$\|R\|_{L_t^2L_x^2}\le\RR_0,\qquad \|\Roc\|_{L_t^2L_x^2}\le\RR_0+\De_0^2.
$$

\subsection{Symbolic notations}
Throghout this section
we shall denote by $\Dcal^{-1}$ any of the Hodge operators
$\dcal^{-1}, \dcall^{-1}$, $\dcalll^{-1}, \dcallll^{-1}$ encountered 
before. We also denote by $\Dcal^{-2}$ the operators  $\dcall^{-1}\c\dcal^{-1}$
or  $\dcal^{-1}\c\dcalll^{-1}$. To avoid confusion in the application
of the Bianchi
identities  we need to be  more precise
in the meaning of expressions such as  $\Dcal^{-1}\Roc$ or $\Dcal^{-2}\Roc$.
\begin{definition} Throughout this section we use the following conventions:
\begin{itemize}
\item $\Roc$ denotes either  the pair $(\rhoc,-\sic)$ or $\bboc$.
\item $\Dcal^{-1}\Roc$ denotes either $\dcal^{-1}(\rho,-\si)$ or $\dcalll^{-1}\bboc$.
\item $\Dcal^{-2}\Roc$ denotes either $\dcall^{-1}\c\dcal^{-1}(\rho,-\si)$
or\footnote{or $\dcal^{-1}\c J\c \dcalll^{-1}\bboc$ where
$J$ is the involution which takes the pair  $(\rho, \si)$
 into $(-\rho,\si)$. This involution will appear later
on but it clearly does not affect the estimates. } $\dcal^{-1}\c\dcalll^{-1}\bboc$.
\item $\Dcal^{-1}\c\ddd_L(\Roc)$ denotes $\dcalll^{-1}\c\ddd_L\bboc$ or $\dcal^{-1}L(\rhoc,-\sic)$
\item $\Dcal^{-2}\c\ddd_L(\Roc)$ denotes $\dcal^{-1}\c\dcalll^{-1}\c\ddd_L\bboc,$ 
or $,\dcall^{-1}\c\dcal^{-1}L(\rhoc,-\sic)$
\end{itemize}  
\label{def:symbolicnot}  
\end{definition}
We recall the  renormalized Bianchi  equations  \eqref{eq:transportrho''}--
 \eqref{eq:transportmodifiedbbgeod},   
\beaa
L(\rhoc)+\frac{3}{2}\trch\c\rhoc &=&\div\b-\ze\c\b
+\f12\chih\c(\nab\hot\ze+\f12\big((\trchb+\frac{2}{r})-\frac{2}{r}\big)\c\chih-\ze\hot\ze).\\
L(\sic)+\frac{3}{2}\trch\c\sic
&=&-\curl\b +\ze\wedge\b
+\f12\chih\wedge(\nab\hot\ze-\ze\hot\ze)\\
\ddd_L(\bboc) &=&-\nab\rho+(\nab \si)^\star -2 (\nab\hot\ze)\c\ze   
+3(\ze\c\rho-\ze^\star
\si)-\trch \bb\\
&+&2\ze\c(-\f12\trch\c\chibh -\f12 \trchb\c\chih
+\ze\hot\ze)-4\chi\c\chibh\c\ze\nn
\eeaa
which can be written schematically in the form,
\beaa
\ddd_L\bboc&=&\dcalll (\rho,\si)+\frac{1}{r} R_0+\Gd\c(R_0+\nab \Gd+\Gd\c \Bd)
\label{eq:schematicmodB1}\\
L(\rhoc,
-\sic) & =&\dcal\b+\Gd\c(R_0+\nab \Gd+\frac{1}{r}\Gd +\Gd\c \Bd)
\eeaa
or, according to definition \ref{def:symbolicnot},
\be{eq:schematicmodB2}
\Dcal^{-1}\c\ddd_L\Roc=R_0+\Dcal^{-1}\bigg(\frac{1}{r} R_0+\Gd\c(R_0+\nab
\Gd+\frac{1}{r}\Gd +\Gd\c
\Bd)\bigg)
\end{equation}

We shall need  the following commutators, 
\begin{definition}
 \beaa 
C_1(\Roc)&=&\nab\c\dcall^{-1}\c [\ddd_L\,,\,\dcal^{-1}](\rhoc,-\sic)\quad \mbox{or} \quad
 \nab\c[\ddd_L\,,\,\dcal^{-1}]\c\dcalll^{-1}(\bboc)\\
C_2(\Roc)&=&\nab\c [\ddd_L\,,\,\dcall^{-1}]\c\dcal^{-1}(\rhoc,-\sic)\quad \mbox{or} \quad
 \nab\c[\ddd_L\,,\,\dcal^{-1}]\c\dcalll^{-1}(\bboc)\\
C_3(\Roc)&=&[\ddd_L\,,\,\nab]\c\dcall^{-1}\c\dcal^{-1}(\rhoc,-\sic) \quad \mbox{or} \quad
[\ddd_L\,,\,\nab]\c\dcall^{-1}\c\dcal^{-1}(\bboc)
\eeaa
We shall treat  the two different cases together by simply
writing,
\beaa
C_1(\Roc)&=&\nab\c\Dcal^{-1}\c [\ddd_L\,,\,\Dcal^{-1}], \quad C_2(\Roc)=\nab\c
[\ddd_L\,,\,\Dcal^{-1}]\c\Dcal^{-1}\\\quad C_3(\Roc)&=&[\ddd_L\,,\,\nab]\c\Dcal^{-2}\Roc
\eeaa
We also write,
$$\,\,C(\Roc)=\big(C_1(\Roc),C_2(\Roc),C_3(\Roc)\big). $$
\label{def:maincommform}
\end{definition}
\subsection{Error type terms}
To analyze various error terms which appear below 
we make the following definition:
\begin{definition}
We say that an expression $E$  is a strong error type term, denoted  $\err$ in what follows, if 
\be{eq:deferror}
\|E\|_{\PP^\ep}\les \De_0^2
\end{equation}
for some fixed\footnote{In fact  all our strong
error terms  are of type $\err$ with $\ep=\frac{1}{10}$. } $\ep>0$ .

We say that an expression $E$  is a weak 
 error type term, denoted  $\Err$ in what follows, if 
\be{eq:deferror-weak}
\|E\|_{\PP^0}\les \De_0^2
\end{equation}
\label{def:errortype}
\end{definition}
\begin{proposition}
 The following expressions are   error type in the sense of  \eqref{eq:deferror}
\beaa
E_1&=&\Dcal^{-1}(\Gd\c\nab \Gd), \quad \nab\c\Dcal^{-2}(\Gd\c \nab \Gd)\\
E_2&=&\Dcal^{-1}( \Gd\c R),\quad \quad \nab\c\Dcal^{-2}(\Gd\c R),\quad \Dcal^{-1}( \Gd\c\nab\c\Dcal^{-1}
R)\\ E_3&=&\Dcal^{-1}( \Gd\c \Roc),\quad \quad \nab\c\Dcal^{-2}(\Gd\c \Roc),\quad \Dcal^{-1}(
\Gd\c\nab\c\Dcal^{-1}
\Roc)\\
E_4&=& \Dcal^{-1}(\Gd\c \Gd\c \Gd),\,\,\, \Dcal^{-1}(\Gd\c \Gd\c \Bd),  \,\,\,
 \nab\c\Dcal^{-2}(\Gd\c \Gd\c \Gd),\,\,\,\\
&&\nab\c\Dcal^{-2}(\Gd\c \Gd\c \Bd)\\
E_5&=&\Gd\c\Gd,\quad  \nab\c\Dcal^{-1} (\Gd\c\Gd)
\eeaa
The following expressions are of weak error type in the sense of 
\eqref{eq:deferror-weak} :
\beaa
E_6&=& \Gd\c \Bd,\quad  \nab\c\Dcal^{-1}(\Gd\c \Bd), \\
\eeaa
\label{prop:errortype}
\end{proposition}
\begin{proof}:\quad Recall  that $\nab \c\dcal^{-1}$, $\nab\c\dcall^{-1}$  and 
$\nab\c\dcalll^{-1}$ are  bounded
operators on $\PP^\th$, see lemma \ref{prop:funnyhodgeestimPP}.
 We shall also make use
of the estimate \eqref{eq:Lt2LxpDcal-1PP} of lemma \ref{prop:SobolevBesovLp}, with
any $\frac{2}{2-\th}<p\le 2$ and $\Dcal^{-1}=\dcal^{-1}, \dcall^{-1}$ or $\dcalll^{-1}$, 
\be{eq:Lt2LxpDcal-1PPagain}
\|\Dcal^{-1} f\|_{\PP^\th}\les \|f\|_{L_t^2 L_x^p}
\end{equation}
Also,
\be{eq:Lt2LxpDcal-1PPagainn}
\|\nab\c\Dcal^{-2} f\|_{\PP^\th}\les \|f\|_{L_t^2 L_x^p}
\end{equation}
Thus, for  $\th=\ep<\f12$
\beaa
\|E_{1}\|_{\PP^\th}&\les& \|\Gd\c\nab \Gd\|_{L_t^2L_x^{\frac{4}{3}}}
\les\|\Gd\|_{L_t^\infty L_x^4}\c\|\nab \Gd\|_{L_t^2L_x^2}
\les \De_0^2\\
\|E_{2}\|_{\PP^\th}&\les& \|\Gd\c R\|_{L_t^2L_x^{\frac{4}{3}}}
\les\|\Gd\|_{L_t^\infty L_x^4}\c  \|R\|_{L_t^2L_x^2}
\les \De_0^2\\
\|E_{3}\|_{\PP^\th}&\les& \|\Gd\c \Roc\|_{L_t^2L_x^{\frac{4}{3}}}
\les\|\Gd\|_{L_t^\infty L_x^4}\c  \|\Roc\|_{L_t^2L_x^2}
\les \De_0^2
\eeaa
Similarily, for $\th=\ep<\frac{1}{4}$,
\beaa
\|E_{4}\|_{\PP^\th}&\les& \|\Gd\c \Gd\c \Bd\|_{L_t^2L_x^{\frac{8}{7}}}
\les\|\Gd\|_{L_t^2L_x^8}\c \|\Gd\|_{L_t^\infty L_x^4}\c \| \Bd\|_{L_t^\infty L_x^2}\\
&\les&\|\Gd\|_{L_t^2L_x^8}\c\|\Gd\|_{L_t^\infty L_x^4}\c 
\| \Bd\|_{
L_x^2L_t^\infty}\les
\De_0^3
\eeaa
To estimate $E_5$ we need to apply the nonsharp  product
 estimate \eqref{eq:bilin2newnew}, or \eqref{eq:anotherbilPP0},
of proposition \ref{prop:5.7}
 as well as 
the boundedness of $\nab\c\Dcal^{-1}$ on $\PP^\ep$,
\beaa
\|E_5\|_{\PP^\ep}&\les&\|\Gd\c \Gd\|_{\PP^\ep}\les \NN_1(\Gd)
  \c \|\Gd\|_{\BB^\ep}\les \NN_1(\Gd)^2\les \De_0^2
\eeaa
as desired.

The proof that    $E_6$ is a weak  error type term
is much more subtle. This will be proved only at the end of the next section
by an elaborate bootstrap argument. In fact we shall assume, by the  auxilliary 
boostrap argument {\bf BA5}, that $\|E_6\|_{\PP^0} \le \De_0^2$
and show that in fact $E_6\le \f12 \De_0^2$. We note however that error terms
of type $E_6$ are not needed in the discussions of this section.
\end{proof}
\begin{proposition}The commutator   $[\Dcal^{-1},\ddd_L](\Roc)$
is of strong error type,
$$\|[\Dcal^{-1},\ddd_L](\Roc)\|_{\PP^\ep}\les \De_0^2$$
\label{prop:comm[DL]first}
\end{proposition}
\begin{proof}:\quad
To calculate the precise form of the commutator $[\Dcal^{-1},\ddd_L]$
we need to make use of 
 lemma \ref{le:commutationL},
\be{eq:commDcal-1}
[\ddd_L,\Dcal^{-1}]=\Dcal^{-1}[\ddd_L,\Dcal]\Dcal^{-1}.
\end{equation}
To calculate $[\ddd_L,\Dcal]$
commutation formulas of proposition
\ref{prop:commutationL} which we 
write schematically\footnote{The dot products appearing in
\eqref{eq:commutationlnab} refer to apropriate tensor products and traces consistent with 
 the precise commutation formulas of \eqref{prop:commutationL}.}, for arbitrary $S$-tangent
tensor $F$ on $\HH$ as follows,
\be{eq:commutationlnab}
[\ddd_L,\nab ]F=\Gd\c\nab F+\Gd\c \Gd\c F+\b\c F
\end{equation}
According to   formulas \eqref{eq:commDcal-1}, \eqref{eq:commutationlnab}, 
\bea
 [\ddd_L,\Dcal^{-1}]\Roc
&=&\Dcal^{-1}\bigg(\Gd\c(\nab\c \Dcal^{-1}\Roc)+\Gd\c \Gd\c(\Dcal^{-1}\Roc)+
\b\c(\Dcal^{-1}\Roc)\bigg)\nn\\
&=&I_1+I_2+I_3\label{eqI1I2I3commL}
\eea 
Therefore,
\be{eq:comm100}
\|[\ddd_L,\Dcal^{-1}]\Roc\|_{\PP^\ep}\les \|I_1\|_{\PP^\ep}
+ \|I_2\|_{\PP^\ep}+\|I_3\|_{\PP^\ep}
\end{equation} and 
\beaa
\|I_1\|_{\PP^\ep}&\les&
\|\Gd\c(\nab\c\Dcal^{-1}) \Roc\|_{L_t^2L_x^{4/3}}
\les\|\Gd\|_{L_t^\infty L_x^4}\c  \|\Roc\|_{L_t^2L_x^2}\les \De_0^2\label{eq:commI1}\\
\|I_2\|_{\PP^\ep}&\les&\|\Gd\c \Gd\c \Dcal^{-1}\Roc\|_{L_t^2 L_x^{4/3}}
\les\|\Gd\|_{L_t^\infty L_x^4}\c \|\Gd\|_{L_t^\infty L_x^4}\c \| \Dcal^{-1}\Roc\|_{L_t^2 L_x^4}\nn\\
&\les&
\De_0^2\c \|\Dcal^{-1}\Roc \|_{L_t^2
L_x^4}\les \De_0^3 \label{eq:commI2}\\
\|I_3\|_{\PP^\ep}
&\les&\|\b\c \Dcal^{-1}\Roc\|_{L_t^2L_x^{4/3}}
\les\|R\|_{L_t^2L_x^2}\c\|\Dcal^{-1}\Roc\|_{L_t^\infty
L_x^4}\\
&\les&\De_0\|\Dcal^{-1}\Roc \|_{L_t^\infty L_x^4}
\eeaa
In view of  \eqref{eq:LtinftyLx4}\label{eq:commI3}
$
\|\Dcal^{-1}\Roc\|_{L_t^\infty L_x^4}\les\NN_1(\Dcal^{-1}\Roc)$,
whence,
$
 \|I_3\|_{L_t^2L_x^2}
\les\De_0\c\NN_1(\Dcal^{-1}\Roc).
$
Thus, back to \eqref{eq:comm100},
\be{eq:estimateI3comm}
\|[\ddd_L,\Dcal^{-1}]\Roc\|_{\PP^\ep}\les \De_0^2+\De_0\c\NN_1(\Dcal^{-1}\Roc)
\end{equation}
To finish the proof we need the estimate $\NN_1(\Dcal^{-1}\Roc)\les\De_0$,
to be proved in the next proposition.
\end{proof}
\begin{proposition} 
 We have the following estimates: 
\bea
\|\Dcal^{-1}\ddd_L\Roc\|_{L_t^2L_x^2} &\les&\RR_0+\De_0^2\label{eq:RocD-1L2estim}\\
\|[\ddd_L,\Dcal^{-1}]\Roc\|_{L_t^2L_x^2}&\les&\De_0^2\label{eq:RocD-1L2commestim}\\
\NN_1(\Dcal^{-1}\Roc)&\les&\RR_0+\De_0^2\label{eq:RocD-1NN1estim}
\eea
\label{prop:RocD-1L2commestim}
\end{proposition}
\begin{proof}:\quad 
To derive   \eqref{eq:RocD-1L2estim}, 
 we use the symbolic expression of the  renormalized Bianchi  equations 
\eqref{eq:schematicmodB2},
$$\Dcal^{-1}\c\ddd_L\Roc=\Roc+\Dcal^{-1}\bigg(\frac{1}{r} R_0+\Gd\c(R_0+\nab
\Gd+\frac{1}{r}\Gd +\Gd\c
\Bd)\bigg)$$
Observe that, in view of  the error estimates
of proposition \ref{prop:errortype} the term 
$$\Dcal^{-1}\bigg(\Gd\c(R_0+\nab
\Gd+\frac{1}{r}\Gd +\Gd\c
\Bd)\bigg)$$ is a strong error term.
Therefore,
\beaa
\|\Dcal^{-1}\c\ddd_L\Roc\|_{L_t^2L_x^2}&\les& \|R\|_{L_t^2L_x^2}+\|\Dcal^{-1}(\frac{1}{r}
R_0)\|_{L_t^2L_x^2}+\De_0^2\\
&\les&\RR_0+\De_0^2
\eeaa

To prove \eqref{eq:RocD-1L2commestim} we go back to
\eqref{eq:estimateI3comm} from where, in particular,
\be{eq:estimateI3comm'}
\|[\ddd_L,\Dcal^{-1}]\Roc\|_{L_t^2L_x^2}\les \De_0^2+\De_0\c\NN_1(\Dcal^{-1}\Roc)
\end{equation}
On the other hand, according to the elliptic estimates
of proposition \ref{prop:Hodgeestimates}  as well as \eqref{eq:RocD-1L2estim}
\beaa
\NN_1(\Dcal^{-1}\Roc)&=&\|\Dcal^{-1}\Roc\|_{L_t^2L_x^2}+\|\nab\c\Dcal^{-1}\Roc\|_{L_t^2L_x^2}
+\|\ddd_L\c\Dcal^{-1}\Roc\|_{L_t^2L_x^2}\\
&\les&\|\Roc\|_{L_t^2L_x^2}+   \|\Dcal^{-1}\c\ddd_L\Roc\|_{L_t^2L_x^2}  +
\|[\ddd_L,\Dcal^{-1}]\Roc\|_{L_t^2L_x^2}\\
&\les&\RR_0+\De_0^2+\|[\ddd_L,\Dcal^{-1}]\Roc\|_{L_t^2L_x^2}
\eeaa
Thus, back to \eqref{eq:estimateI3comm'}
\beaa
\|[\ddd_L,\Dcal^{-1}]\Roc\|_{L_t^2L_x^2}&\les&
\De_0^2+\De_0\c\big(\RR_0+\De_0^2+\|[\ddd_L,\Dcal^{-1}]\Roc\|_{L_t^2L_x^2}\big)\\
&\les&\De_0^2+\De_0\c \|[\ddd_L,\Dcal^{-1}]\Roc\|_{L_t^2L_x^2}
\eeaa
Thus, for small $\De_0$,
 \beaa
\|[\ddd_L,\Dcal^{-1}]\Roc\|_{L_t^2L_x^2}&\les& \De_0^2\\
\NN_1(\Dcal^{-1}\Roc)&\les&\RR_0+\De_0^2
\eeaa
as desired.
\end{proof}
\subsection{ Preliminary estimates for the commutators $C(\Roc)$.}\quad\label{sec:comm}

In the next section we shall need
an extension of proposition \ref{prop:comm[DL]first}
to commutators of the form
$C(\Roc)=[\ddd_L\,,\,\nab\c\Dcal^{-2}](\Roc)$. Unfortunately this is not possible;
we need to perform in fact a nontrivial renormalization. Clearly
$C(\Roc)$ can be expressed as a sum of the  commutators,
 see  definition \ref{def:maincommform}:
$
C_1(\Roc)=\nab\c\Dcal^{-1}\c [\ddd_L\,,\,\Dcal^{-1}]$,\,\,  $C_2(\Roc)=\nab\c
[\ddd_L\,,\,\Dcal^{-1}]\c\Dcal^{-1}$, \,\, $C_3(\Roc)=[\ddd_L\,,\,\nab]\c\Dcal^{-2}\Roc$.
 Recall that
  $$\,\,C(\Roc)=\big(C_1(\Roc),C_2(\Roc),C_3(\Roc)\big) $$ and
schematically, 
\bea
\ddd_L\c(\nab\c\Dcal^{-2}\Roc)&=&\nab\c\Dcal^{-2}\c\ddd_L(\Roc)+C(\Roc)\label{eq:maincommform}\\
\nab\c\Dcal^{-1}\c \ddd_L\c \Dcal^{-1}(\Roc)&=&\nab\c\Dcal^{-2}\c\ddd_L(\Roc)+C(\Roc), \nn\\ 
\nab \c\ddd_L\c\Dcal^{-2}(\Roc)&=&\nab\c\Dcal^{-2}\c\ddd_L(\Roc)+C(\Roc).\nn
\eea

In view of lemma \ref{le:commutationL}
and commutation  formula \eqref{eq:commutationlnab},
\bea
C_1(\Roc)&=&\nab\c\Dcal^{-1}\c [\ddd_L\,,\,\Dcal^{-1}](\Roc)
=\nab\c\Dcal^{-1}\c\Dcal^{-1} [\ddd_L\,,\,\Dcal](\Dcal^{-1}\Roc)\nn\\
&=&\nab\c\Dcal^{-2}\bigg(\Gd\c \nab\c(\Dcal^{-1}\Roc)+\Gd\c \Gd\c (\Dcal^{-1}\Roc)+\b\c
(\Dcal^{-1}\Roc)\bigg)\label{eq:C1}
\eea
\bea
C_2(\Roc)&=&\nab\c [\ddd_L\,,\,\Dcal^{-1}]\c\Dcal^{-1}(\Roc)
=\nab\c\Dcal^{-1}\c[\ddd_L\,,\,\Dcal]\c\Dcal^{-2}(\Roc)\nn\\
&=&\nab\c\Dcal^{-1}\bigg(\Gd\c \nab \c\Dcal^{-2}(\Roc)+\Gd\c \Gd
\c\Dcal^{-2}(\Roc)\bigg)\nn\\
&+&\nab\c\Dcal^{-1}\c\big(\b \c\Dcal^{-2}(\Roc)\big)\label{eq:C2}
\eea
\bea
C_3(\Roc)&=&[\ddd_L,\nab]\c \Dcal^{-2}(\Roc)\nn\\
&=&\Gd\c\nab \c\Dcal^{-2}(\Roc)+\Gd\c \Gd\c \Dcal^{-2}(\Roc)
+\b\c \Dcal^{-2}(\Roc)\label{eq:C3}
\eea
\begin{proposition} The commutator $C_1(\Roc)$ is of strong error type. The other
commutators $C_2(\Roc), C_3(\Roc)$ verify:
\bea
C_2(\Roc)&=&\nab\c\Dcal^{-1}\c\big(\b \c\Dcal^{-2}(\Roc)\big)+\err\label{eq:C2(Roc)form}\\
C_3(\Roc)&=&\b\c \Dcal^{-2}(\Roc)+\err\label{eq:C3(Roc)form}
\eea
\label{prop:C1C2C3errorestim}
\end{proposition}
\begin{proof}:\quad The three  terms  in the expression of 
$C_1(\Roc)$ look precisely like the terms  $I_1,I_2,I_3$ in the formula 
\eqref{eqI1I2I3commL}, with the  $\Dcal^{-1}$ in front  replaced by
$\nab \c\Dcal^{-2}$,  
and therefore can be estimated as in 
\eqref{eq:comm100}.

 To estimate
the terms $J_1=\nab\c\Dcal^{-1}\big(\Gd\c \nab
\c\Dcal^{-2}(\Roc)\big)$ and $J_2=\nab\c\Dcal^{-1}\big(\Gd\c \Gd
\c\Dcal^{-2}(\Roc)\big)$ on the right hand side of \eqref{eq:C2}
we  make use of  the nonsharp product   estimate \eqref{eq:bilin2newnew}, or
\eqref{eq:anotherbilPP0} of  proposition \ref{prop:5.7} 
 as well as 
the boundedness of $\nab\c\Dcal^{-1}$ on $\PP^\ep$. The estimate
for $J_1$ is similar to that of 
 $E_5$ in proposition \ref{prop:errortype},
\beaa
\|J_1\|_{\PP^\ep}&\les&
\|\Gd\c(\nab\c\Dcal^{-2}) \Roc\|_{\PP^\ep}
\les\NN_1(\Gd)\c 
\|\nab\c\Dcal^{-2}\Roc\|_{\BB^\ep}\\
&\les&\De_0\c\NN_1(\nab\c\Dcal^{-2}\Roc)
 \eeaa
To estimate $J_2$ we also use\footnote{alternatively one can use
 \eqref{eq:anotherbilPP0}
of the main lemma.}  the product estimate \eqref{eq:bilin2newnew} 
 as follows,
\beaa
\|J_2\|_{\PP^\ep}&\les& \|\Gd\c \Gd
\c\Dcal^{-2}(\Roc)\|_{\PP^\ep}
\les\NN_2(\Dcal^{-2} \Roc)\c
\|\Gd\c \Gd\|_{\PP^\ep}
\eeaa
or, since   according to proposition \ref{prop:errortype}, $\|\Gd\c \Gd\|_{\PP^\ep}\les\De_0^2$,
\be{eq:C22Roc}\|J_2\|_{\PP^\ep}\les
\De_0^2\c\NN_2(\Dcal^{-2} \Roc)
\end{equation}
Consequently,
\be{eq:C2(Roc)form'}
\|C_2(\Roc)-\nab\c\Dcal^{-1}\c\big(\b \c\Dcal^{-2}(\Roc)\big)\|_{\PP^\ep}\les
\De_0 \c\NN_1(\nab\c\Dcal^{-2}\Roc)+   \De_0^2 \c\NN_2(\Dcal^{-2} \Roc)
\end{equation}
The estimates for $C_3(\Roc)$ are exactly the same as for $C_2(\Roc)$.
Therefore,
\be{eq:C3(Roc)form'}
\|C_3(\Roc)-(\b \c\Dcal^{-2}(\Roc)\|_{\PP^\ep}\les
\De_0 \c\NN_1(\nab\c\Dcal^{-2}\Roc)+   \De_0^2 \c\NN_2(\Dcal^{-2} \Roc)
\end{equation}
Thus the proof of the proposition
reduces to the following estimates:
\beaa
\NN_1(\nab\c\Dcal^{-2}\Roc)&\les& \RR_0+\De_0^2\\
\NN_2(\Dcal^{-2}\Roc)&\les& \RR_0+\De_0^2
\eeaa
We shall prove in fact the following proposition.
\end{proof}
\begin{proposition}
 We have the following estimates: 
\bea
\|\nab\c\Dcal^{-2}\c\ddd_L\Roc\|_{L_t^2L_x^2} &\les&\RR_0+\De_0^2\label{eq:RocnabD-2L2estim}\\
\|C(\Roc)\|_{L_t^2L_x^2}&\les&\De_0^2\label{eq:RocnabD-2L2commestim}\\
\NN_1(\nab\c\Dcal^{-2}\Roc)&\les&\RR_0+\De_0^2\label{eq:RocnabD-2NN1estim}\\
\NN_2(\Dcal^{-2}\Roc)&\les& \RR_0+\De_0^2\label{eq:RocD-2NN2estim}
\eea
As a corollary of the above estimates, combined with corollary \ref{corr:calcineq},  we also have,
\bea
\|\Dcal^{-1}\Roc\|_{L_t^\infty L_x^4}&\les& \RR_0+\De_0^2\label{eq:D-1RocLtinftyLx4}\\
\|\nab\c\Dcal^{-2}\Roc\|_{L_t^\infty L_x^4}&\les& \RR_0+\De_0^2\label{eq:D-1RocLtinftyLx4'}\\
\|\Dcal^{-2}\Roc\|_{L_t^\infty L_x^\infty}&\les& \RR_0+\De_0^2\label{eq:RocD-2Linftyestim}
\eea
\label{prop:L2C(Roc)}
\end{proposition}
\begin{proof}:\quad The proof of \eqref{eq:RocnabD-2L2estim} is 
similar
to that of \eqref{eq:RocD-1L2estim}. 
Indeed, in view of \eqref{eq:schematicmodB2},
$$
\nab\c\Dcal^{-2}\c\ddd_L\Roc=\nab\c\Dcal^{-1}\Roc+\nab\c\Dcal^{-2}\bigg(\frac{1}{r} R_0+\Gd\c(R_0+\nab
\Gd+\frac{1}{r}\Gd +\Gd\c
\Bd)\bigg)
$$
Therefore, in view of the elliptic-Hodge  estimates of
proposition \ref{prop:Hodgeestimates} and the error estimates
of proposition \ref{prop:errortype}
$$\|\nab\c\Dcal^{-2}\c\ddd_L\Roc\|_{L_t^2L_x^2}\les\RR_0+\De_0^2.
$$
as desired.

To estimate the  $L^2$ norms of the comutators $C_i(\Roc)$, $i=1,2,3$
 we make use of the estimates \eqref{eq:C2(Roc)form'} and
 \eqref{eq:C3(Roc)form'} derived in the proof of
proposition \ref{prop:C1C2C3errorestim},
\beaa
\|C_1(\Roc)\|_{L_t^2L_x^2}&\les&\De_0^2\\
\|C_2(\Roc)\|_{L_t^2L_x^2}&\les&
\|\nab\c\Dcal^{-1}\c\big(\b \c\Dcal^{-2}(\Roc)\big)\|_{L_t^2L_x^2}\\
&+&
\De_0\c \NN_1(\nab\c\Dcal^{-2}\Roc)+   \De_0^2 \c\NN_2(\Dcal^{-2} \Roc)\\
\|C_3(\Roc)\|_{L_t^2L_x^2}&\les&
\|\b \c\Dcal^{-2}(\Roc)\|_{L_t^2L_x^2}\\
&+&
\De_0\c \NN_1(\nab\c\Dcal^{-2}\Roc)+   \De_0^2 \c\NN_2(\Dcal^{-2} \Roc)
\eeaa
On the other hand, making use of the estimate
$\|F\|_{L_t^\infty L_x^\infty}\les\NN_2(F)$
of corollary   \ref{corr:calcineq}
\beaa
\|\nab\c\Dcal^{-1}\c\big(\b \c\Dcal^{-2}(\Roc)\big)\|_{L_t^2L_x^2}&\les& \|\b
\c\Dcal^{-2}(\Roc)\|_{L_t^2L_x^2}\les \|\b\|_{L_t^2L_x^2}\c\|\Dcal^{-2}(\Roc)\|_{L_t^\infty L_x^\infty}\\
&\les&\De_0\c \NN_2(\Dcal^{-2}(\Roc))
\eeaa
Consequently,
\be{eq:L2estimC(roc)1}
\|C(\Roc)\|_{L_t^2L_x^2}\les\De_0^2+\De_0\c\big(\NN_1(\nab\c\Dcal^{-2}\Roc)+ \NN_2(\Dcal^{-2}\Roc)\big)
\end{equation}
It remains to estimate
$\NN_1(\nab\c\Dcal^{-2}\Roc), \,\, \NN_2(\Dcal^{-2}\Roc)$.
From the definition of $C(\Roc)$ and the estimate
\eqref{eq:RocnabD-2L2estim} we infer that,
\beaa
\|\nab\c\ddd_L\c\Dcal^{-2}\Roc\|_{L_t^2L_x^2}&\les&
\RR_0+\|C(\Roc)\|_{L_t^2L_x^2}+\De_0^2\label{eq:RocD-2L2estim''}\\
\|\ddd_L\c\nab\c\Dcal^{-2}\Roc\|_{L_t^2L_x^2} 
&\les&\RR_0+\|C(\Roc)\|_{L_t^2L_x^2}+\De_0^2\label{eq:RocD-2L2estim'''}
\eeaa
On the other hand, using also  the elliptic
estimates
of proposition \ref{prop:secondorderHodge},
$$\|\nab^2\c\Dcal^{-2}\Roc\|_{L_t^2L_x^2}\les\RR_0+
\De_0\|C(\Roc)\|_{L_t^2L_x^2}$$
Therefore,
\beaa
\NN_1(\nab\c\Dcal^{-2}\Roc)&\les&\RR_0+\|C(\Roc)\|_{L_t^2L_x^2}+\De_0^2\\
\NN_2(\Dcal^{-2}\Roc)&\les&\RR_0+\|C(\Roc)\|_{L_t^2L_x^2}+\De_0^2
\eeaa
Combining these with \eqref{eq:L2estimC(roc)1}
and using the smallness of $\De_0$, we conclude:
\beaa
\|C(\Roc)\|_{L_t^2L_x^2}&\les&\De_0^2\\
\NN_1(\nab\c\Dcal^{-2}\Roc)&\les&\RR_0+\De_0^2\\
\NN_2(\Dcal^{-2}\Roc)&\les&\RR_0+\De_0^2
\eeaa
as desired.
\end{proof}
\subsection{Decomposition and correction estimates for $C(\Roc)$.}\quad

The terms $\nab\c\Dcal^{-1}\c\big(\b \c\Dcal^{-2}(\Roc)\big)$
and $\b\c \Dcal^{-2}(\Roc)$ appearing on the right hand side 
of the formulas \eqref{eq:C2(Roc)form} and \eqref{eq:C3(Roc)form}
in proposition \ref{prop:C1C2C3errorestim} are not  weak  error type.
This is a serious technical  difficulty requiring another renormalization.
The idea is to express both of them, modulo strong error terms,  as $\ddd_L$ derivatives of
other tensors controllable in the $\NN_1$ norm. This is done
in the following:

\begin{lemma}[Decomposition Lemma]
Consider the expression $\b\c F$ with $F$ a tensor of arbitrary order on $\HH$ verifying
\be{eq: niceestimL2forH} \NN_2(F) <\infty.
\end{equation} 
{\bf i.)}:\quad We claim that,
\be{eq:b=DL+E0}
\b=\ddd_L\big(\dcal^{-1}(\Roc)\big)+E_0\qquad \mbox{with } \quad\|E_0\|_{\PP^\ep}\les \De_0^2.
\end{equation}
{\bf ii.) }\quad There exists a decomposition,
\bea
\b\c F&=&\ddd_L P+E\label{eq:decompb}
\eea
where $P, E$ are tensors( of the same type as $\b\c F$) verifying,
\bea
\NN_1(P),\, \|E\|_{\PP^\ep} &\les& \De_0\c \NN_2(F)\label{eq:L2estimP}
\eea
{\bf iii.)}\quad There exist tensors $\bar{P}, \bar{E}$
verifying the same conditions \eqref{eq:L2estimP} as above
such that,
\bea
\nab\c\Dcal^{-1}(\b\c F)&=&\ddd_L \bar{P}+\bar{E}\label{eq:decompbwithnabDcal-1}
\eea 

{\bf iv.)}\quad  There exists a decomposition,
\bea
\Dcal^{-1}(\b\c F)&=&\ddd_L \tilde{P}+\tilde{E}\label{eq:decompb'}
\eea
where $\tilde{P}, \tilde{E}$ are tensors verifying,
\bea
\NN_2(\tilde{P}),\, \|\nab \tilde{E}\|_{\PP^\ep} &\les& \De_0\c \NN_2(F)\label{eq:L2estimP'}
\eea
\label{le:decompbandmore}
\end{lemma}
\begin{proof}:\quad
  The first statement {\bf i}) follows
immediately from  the renormalized Bianchi, see
\eqref{eq:schematicmodB2}, and propositions \ref{prop:errortype},
\ref{prop:comm[DL]first}, 
\beaa
\b&=&\dcal^{-1}\c\ddd_L\Roc+\err=\ddd_L\c\dcal^{-1}(\Roc)+\err
\eeaa
 Therefore
$
\b=\ddd_L\big(\Dcal^{-1}(\Roc)\big)+E_0$, with $ \|E_0\|_{\PP^\ep}\les \De_0^2
$\,\,
as desired.

 To prove {\bf ii}) we use  \eqref{eq:b=DL+E0} and  write,
\beaa
\b\c F&=&\big(\ddd_L(\Dcal^{-1} \Roc) +E_0\big)\c F=\ddd_L(\Dcal^{-1}\Roc\c F)-(\Dcal^{-1}\Roc) \c \ddd_L
F +E_0\c F
\eeaa
We now  observe   that $E_0\c F$ is an error type term ;
 to show this we need to make use of the product estimate
\eqref{eq:bilin2newnew} of proposition \ref{prop:5.7}, i.e.
$$
\| F\c E_0\|_{\PP^\ep} \les 
\NN_2(F)\c\|E_0\|_{\PP^\ep}\les \De_0 \c \NN_2(F).
$$
 We next show that
$(\Dcal^{-1}\Roc) \c \ddd_L F$ is also error type. To do this we need to appeal  to
 the bilinear estimate  \eqref{eq:anotherbilPP0}
of the proposition \ref{prop:5.7},
 \beaa
\|(\Dcal^{-1}\Roc) \c \ddd_L F\|_{\PP^\ep}&\les& \NN_1(\Dcal^{-1}\Roc)\c \big(\|\ddd_L F\|_{L_t^2
L_x^2}+\|\nab \c\ddd_L F\|_{L_t^2 L_x^2}\big)\\
&\les&(\RR_0+\De_0^2)\c\NN_2(F)\les \De_0\c\NN_2(F)
\eeaa
as desired.
Consequently, setting $P=(\Dcal^{-1}\Roc)\c F$, we infer that,
$$\b\c F=\ddd_L P +E,\qquad \|E\|_{\PP^\ep}\les \De_0 \c \NN_2(F)$$
 It remains to check that $P$ verifies \eqref{eq:L2estimP}. Using 
\eqref{eq:D-1RocLtinftyLx4'} and 
\eqref{eq:RocD-2Linftyestim},
\beaa
\| \ddd_LP\|_{L_t^2L_x^2}&\les& \|\b\c F\|_{L_t^2L_x^2}+\|E\|_{L_t^2L_x^2}\\
&\les&\|\b\|_{L_t^2L_x^2} \c \|F\|_{L_t^\infty L_x^\infty}+\De_0\c \NN_2(F)\les \De_0\c \NN_2(F)\\
\|\nab P\|_{L_t^2L_x^2}&\les& 
\|(\nab \c\Dcal^{-1}\Roc)\c F\|_{L_t^2L_x^2}+\|(\Dcal^{-1}\Roc)\c \nab  F\|_{L_t^2L_x^2}\\
&\les&\|\nab \c\Dcal^{-1}\Roc\|_{L_t^2L_x^2}\c\|F\|_{L_t^\infty L_x^\infty}+
\|\nab  F\|_{L_t^4L_x^4}\|\Dcal^{-1}\Roc\|_{L_t^4L_x^4}
\les \De_0\c \NN_2(F)
\eeaa
 Thus, $\NN_1(P)\les   \De_0\c \NN_2(F)$  as desired.

Part {\bf iii.)} of the lemma is a little more involved
 as we will have to prove our decomposition \eqref{eq:decompbwithnabDcal-1} by
an iteration procedure.

{\bf Step 1.}\quad  In view of the decomposition \eqref{eq:decompb} above,
we have  
$$\b\c H=\ddd_L P_1+E_1, \qquad \NN_1(P_1),\,\, \|E_1\|_{\PP^\ep}\les \De_0^2$$

{\bf Step k.} \quad   We construct  by recurrence two 
sequences  $P_2, P_3\ldots P_k$ and  $E_2, E_3\ldots E_k$ 
such  that,
\beaa
\b \c(\Dcal^{-1} P_{i})&=&\ddd_L P_{i+1}+E_{i+1},\quad i=1,\ldots k-1\\
\NN_1(P_{i+1})&\les&
(C\De_0)^{i}\c\NN_2(H)
\\
\|E_{i+1}\|_{\PP^\ep}&\les&
(C\De_0)^{i}\c\NN_2(H)
\eeaa
Assume $P_i, E_i$ have been already  constructed for $i=2,3\ldots k$. 
We construct $P_{k+1}, E_{k+1}$ by applying \eqref{eq:decompb} to $\b\c\Dcal^{-1}P_k$.  
Thus,
$$\b \c(\Dcal^{-1} P_k)=\ddd_L P_{k+1}+E_{k+1}$$
with, 
\bea
\NN_1(P_{k+1})&\les&\De_0\c \NN_2(\Dcal^{-1} P_k)\les \De_0\c \NN_1(P_k)\le
(C\De_0)^{k}\c\NN_2(H)\label{eq:stepkforP}
\\
\|E_{k+1}\|_{\PP^\ep}&\les&\De_0\c \NN_2(\Dcal^{-1} P_k)
\les \De_0\c \NN_1(P_k)\le
(C\De_0)^{k}\c\NN_2(H)
\label{eq:stepkforE}
\eea
where   $C$ is  a universal constant. These estimates
are based on the following:
\begin{lemma} The following  estimate holds for any $S$-tangent tensor
$P$ on $\HH$,
\bea
\NN_2(\Dcal^{-1} P)&\les& \NN_1(P)\\
\NN_2(\nab\c\Dcal^{-1} P)&\les&  \NN_1(P)
\eea
\label{le:preservofN1N2}
\end{lemma}
\begin{proof}:\quad
We present  the proof at   the end of  this section.
\end{proof}
Now consider,
\beaa
I_{k}&=&\nab\c\Dcal^{-1}\big(\b \c \Dcal^{-1}(P_{k})\big)=\nab\c\Dcal^{-1}\big(
\ddd_L P_{k+1}+E_{k+1})\\
&=&
\ddd_L\big(\nab\c\Dcal^{-1}P_{k+1}\big)+[\nab\c\Dcal^{-1}, \ddd_L] P_{k+1} + 
\nab\c\Dcal^{-1}  E_{k+1}
\eeaa
On the other hand 
\beaa
[\ddd_L, \nab\c\Dcal^{-1}] P&=&[\ddd_L, \nab]\c\Dcal^{-1} P
+\nab\c[\ddd_L, \Dcal^{-1}] P\\
&=&\Gd\c\nab\c \Dcal^{-1} P+\Gd\c \Gd\c\Dcal^{-1} P+\b\c \Dcal^{-1} P\\
&+&\nab\c\Dcal^{-1}\big(\Gd\c\nab\c \Dcal^{-1} P+\Gd\c \Gd\c\Dcal^{-1} P+\b\c
\Dcal^{-1} P\big)\\ &=&[\ddd_L, \nab\c\Dcal^{-1}]_g P+ \b\c \Dcal^{-1} P+
\nab\c\Dcal^{-1}\big(\b\c \Dcal^{-1} P\big)
\eeaa
where we denote by $[\ddd_L, \nab\c\Dcal^{-1}]_g$ the good part of the commutator
$[\ddd_L, \nab\c\Dcal^{-1}]$ i.e. the one which preserves error type terms. More precisely,
\begin{lemma} The ``good commutator'' $[\ddd_L\,,\, \nab\c\Dcal^{-1}]_g$ verifies the following
estimate,
\be{eq:commLforPinP0}
\|\,[\ddd_L\,, \,\nab\c\Dcal^{-1}]_g P\,\|_{\PP^\ep}\les \De_0\c 
\NN_1(P)
\end{equation}

\label{le:goodcomm}
\end{lemma}
We postpone the proof until the end of this subsection.
We have shown that,
\beaa
I_{k}&=&\ddd_L(\nab\c\Dcal^{-1} P_{k+1})+\b\c \Dcal^{-1} P_{k+1}+
 \nab\c\Dcal^{-1}\big(\b\c \Dcal^{-1}
P_{k+1}\big)\\
&+&[\ddd_L, \nab\c\Dcal^{-1}]_g P_{k+1}+\nab\c\Dcal^{-1}  E_{k+1}
\eeaa
We continue as follows,
\beaa
I&=&\nab\c\Dcal^{-1}\big(\b \c H\big)=\nab\c\Dcal^{-1}\big(\ddd_L P_1+E_1)\\
&=&
\ddd_L\big(\nab\c\Dcal^{-1}P_1\big)+\b\c \Dcal^{-1} P_1+ \nab\c\Dcal^{-1}\big(\b\c \Dcal^{-1} P_1\big)\\
&+&[\nab\c\Dcal^{-1}, \ddd_L]_g P_1 + \nab\c\Dcal^{-1}  E_1\\
&=&\ddd_L\big(\nab\c\Dcal^{-1}P_1\big)+\ddd_L P_2+E_2+\nab\c\Dcal^{-1}\big(\ddd_L P_2+E_2\big)\\
&+&[\nab\c\Dcal^{-1}, \ddd_L]_g P_1 + \nab\c\Dcal^{-1}  E_1\\
&=&\ddd_L\big(\nab\c\Dcal^{-1}P_1 +P_2 \big)+\nab\c\Dcal^{-1}\big(\ddd_L P_2\big)+[\nab\c\Dcal^{-1}, \ddd_L]_g P_1 \\
&+&E_2+\nab\c\Dcal^{-1}  E_2+\nab\c\Dcal^{-1}  E_1\\
&=&\ddd_L\big(\bar{P_2})+\nab\c\Dcal^{-1}\big(\ddd_L P_2\big) +\bar{E}_2
\eeaa
i.e.,
$$I=\ddd_L\big(\bar{P_2})+\nab\c\Dcal^{-1}\big(\ddd_L P_2\big) +\bar{E}_2$$
where 
\beaa
\bar{P_2}&=&\nab\c\Dcal^{-1}P_1 +P_2 \\
\bar{E}_2&=&[\nab\c\Dcal^{-1}, \ddd_L]_g P_1+E_2+\nab\c\Dcal^{-1}  E_2+\nab\c\Dcal^{-1}  E_1
\eeaa
Continuing to decompose $\nab\c\Dcal^{-1}\big(\ddd_L P_2\big)$ as above
we arrive to the following formula,
\beaa
I&=&\ddd_L\big(\bar{P_k})+\nab\c\Dcal^{-1}\big(\ddd_L P_k\big) +\bar{E}_k\\
\bar{P_k}&=&\nab\c\Dcal^{-1}\big(P_1+\ldots P_{k-1}\big)+P_2+\ldots P_k\\
\bar{E}_k&=&[\nab\c\Dcal^{-1}, \ddd_L]_g( P_1+\ldots P_{k-1})+\nab\c\Dcal^{-1}\big(E_1+\ldots E_{k}\big) +E_2+\ldots E_{k}
\eeaa
Observe that, in view of lemmas \eqref{le:preservofN1N2}, \eqref{le:goodcomm},
the estimates \eqref{eq:stepkforP}, \eqref{eq:stepkforE}, and the
smallness of
$\De_0$,
\beaa
\NN_1(\bar{P_k}-\bar{P_j})&\les&\NN_1(P_{j+1})+\ldots \NN_1(P_k)\les 
\NN_2(H)\c\sum_{j\ge 1}(C\De)^{j}\les \De_0\c\NN_2(H)\\
\|\bar{E_k}-\bar{E_j}\|_{\PP^\ep}&\les&\|E_{j+1}\|_{\PP^\ep}+\ldots\|E_{k}\|_{\PP^\ep}
+\De_0\big(\NN_1(P_{j+1})+\ldots \NN_1(P_k)\big)\\
&\les&\NN_2(H)\c \sum_{j\ge 1}(C\De)^{j}\les \De_0\c\NN_2(H).
\eeaa
In other words $\bar{P_k}$ forms a Cauchy sequence relative to the norm $\NN_1$,
while  $\bar{E_k}$ forms a Cauchy sequence relative to $\PP^\ep$. Denote by $\bar{P}$ and $\bar{E}$
the corresponding limits. Clearly
 $$\NN_1(\bar{P})\les \De_0\c \NN_2(H),\qquad \|\bar{E}\|_{\PP^\ep}\les \De_0\c\NN_2(H)$$
We also  observe that, for sufficiently small $\De_0$,
\beaa
\|I-\ddd_L\big(\bar{P_k})-\bar{E}_k\|_{L_t^2L_x^2}=
\|\nab\c\Dcal^{-1}\big(\ddd_L P_k\big)\|_{L_t^2L_x^2}\les\NN_1(P_k)\les 
 (C\De_0)^{k}\c\NN_2(H)\rightarrow 0,
\eeaa
as   $k\rightarrow\infty$.
  Therefore  $\|I-\ddd_L\big(\bar{P})-\bar{E}\|_{L_t^2L_x^2}=0$ whence,
$$I=\ddd_L\bar{P}+\bar{E}$$
as desired.
\end{proof}
It remains to prove the last part iv) of  lemma \ref{le:decompbandmore}.
In view of part ii.) 
we can write,
$\b\c H=\ddd_LP+E$, with $\NN_1(P),\,\|E\|_{\PP^\ep}\les \De_0\c\NN_2(H)$.
As before  we can construct a sequence  
$P_1,\ldots P_k,\ldots $ and 
$E_1,\ldots E_k,\ldots $ such that 
$$\b \c(\Dcal^{-1} P_k)=\ddd_L P_{k+1}+E_{k+1}$$
with, 
\bea
\NN_1(P_{k+1}),\,\,\|E_{k+1}\|_{\PP^0} &\les&
(C\De_0)^{k}\c\NN_2(H)\label{eq:stepkforP'}
\eea
where   $C$ is  a universal constant.
Consider now, 
\beaa
I_{k+1}&=&\Dcal^{-1}(\b\c \Dcal^{-1} P_{k+1})=\Dcal^{-1}(\ddd_LP_k+E_k)\\
&=&\ddd_L(\Dcal^{-1}P_k)+[\Dcal^{-1},\ddd_L]P_k+\Dcal^{-1}E_k\\
&=&\ddd_L(\Dcal^{-1}P_k)+\Dcal^{-1}E_k+\Dcal^{-1}\big(G\c \nab \c \Dcal^{-1}P_k 
 +G\c G \c \Dcal^{-1}P_k+\b\c \Dcal^{-1}P_k  \big)\\
&=&\ddd_L(\Dcal^{-1}P_k)+\Dcal^{-1}(\b\c \Dcal^{-1}P_k) +\Dcal^{-1}E_k+
[\Dcal^{-1},\ddd_L]_gP_k
\eeaa
where
\beaa
[\Dcal^{-1},\ddd_L]_gP_k&=&\Dcal^{-1}\big(G\c \nab \c \Dcal^{-1}P_k 
 +G\c G \c \Dcal^{-1}P_k  \big)
\eeaa
verifies,
\beaa
\|\nab\c [\Dcal^{-1},\ddd_L]_g P\|_{\PP^\ep}\les \NN_1(P)
\eeaa
as it is easy to check. The proof now follows precisely as before.

We prove lemmas \ref{le:preservofN1N2}, \ref{le:goodcomm} in the next subsection.
\subsection{Proof of lemmas \ref{le:preservofN1N2}, \ref{le:goodcomm}}
To prove the first inequality of  lemma  \ref{le:preservofN1N2} we write,
\beaa
\NN_2(\Dcal^{-1} P)&=&\|\Dcal^{-1} P\|_{L_2}+\|\nab \Dcal^{-1} P\|_{L_t^2L_x^2}+
\|\nab^2 \Dcal^{-1} P\|_{L_t^2L_x^2}\\
&+&\|\ddd_L\Dcal^{-1} P\|_{L_t^2L_x^2}+\|\nab\ddd_L\Dcal^{-1} P\|_{L_t^2L_x^2}\\
&\les&\|P\|_{L_t^2L_x^2}+\|\nab P\|_{L_t^2L_x^2}+   \|\ddd_L\Dcal^{-1}
P\|_{L_t^2L_x^2}+\|\nab\ddd_L\Dcal^{-1} P\|_{L_t^2L_x^2}\\
&\les&\NN_1(P)+\|[\ddd_L,
\Dcal^{-1}]P\|_{L_t^2L_x^2}+\|\nab\c[\ddd_L,
\Dcal^{-1}]P\|_{L_t^2L_x^2}
\eeaa
 To get the desired estimate it suffices to estimate
$\|\nab\c[\ddd_L,
\Dcal^{-1}]P\|_{L_t^2L_x^2}$.
We calculate the commutator as before,
\beaa
\nab\c[\ddd_L,\Dcal^{-1}]P=\nab\c\Dcal^{-1}
\big(\Gd\c\nab\c\Dcal^{-1}P+\Gd\c \Gd\c\Dcal^{-1}P+\b\c\Dcal^{-1}P\big)
\eeaa
Therefore,
\beaa
\|\nab\c[\ddd_L,
\Dcal^{-1}]P\|_{L_t^2L_x^2}&\les&\|\Gd\c\nab\c\Dcal^{-1}P\|_{L_t^2L_x^2}
+\|\Gd\c \Gd \c\Dcal^{-1}P\|_{L_t^2L_x^2}+\|\b \c\Dcal^{-1}P\|_{L_t^2L_x^2}\\
&\les&\|\Gd\|_{L_t^\infty L_x^4}\c \|\nab \c\Dcal^{-1}P\|_{L_t^2 L_x^4}
\\&+&\big(\|\Gd\c \Gd\|_{L_t^2L_x^2}+
\|\b\|_{L_t^2L_x^2}\big)\c\|\Dcal^{-1}P\|_{L_t^\infty L_x^\infty}\\
&\les& \De_0\c\NN_2(\Dcal^{-1}P)
\eeaa
Hence,
\beaa
\NN_2(\Dcal^{-1}P)&\les&\NN_1(P)+\De_0\c\NN_2(\Dcal^{-1}P),
\eeaa
and therefore $\NN_2(\Dcal^{-1}P)\les \NN_1(P)$ for sufficiently small $\De_0$.
The proof of the second part of lemma \ref{le:preservofN1N2} is similar 
but simpler.

To prove lemma \ref{le:goodcomm} we write,
\beaa
[\ddd_L, \nab\c\Dcal^{-1}]_g P
&=&\Gd\c\nab \c\Dcal^{-1} P+\Gd\c \Gd\c\Dcal^{-1} P \\
&+&\nab\c\Dcal^{-1}\big(\Gd\c\nab \c\Dcal^{-1} P+\Gd\c \Gd\c\Dcal^{-1} P\big)\\
&=&J_1 +J_2+\nab\c\Dcal^{-1}(J_1+J_2)
\eeaa
In view of the $\PP^\ep$ boundedness of $\nab\c\Dcal^{-1}$
it sufffices to estimate $\|J_1\|_{\PP^\ep}$ and $\|J_2\|_{\PP^\ep}$.
These two can be estimated exactly
as in the proof of proposition \ref{prop:C1C2C3errorestim}.
Thus, using the non sharp  bilinear estimate  \eqref{eq:anotherbilPP0} of 
proposition \ref{prop:5.7},
\beaa
\|J_1\|_{\PP^\ep}&=&\|\Gd\c\nab \c\Dcal^{-1} P\|_{\PP^\ep}\les
N_1( \Gd)\c \big (\|\nab^2\c \Dcal^{-1}P\|_{L_t^2L_x^2}+\|\nab\c
\Dcal^{-1}P\|_{L_t^2L_x^2}\big)\\ &\les &
\De_0\c
\big(\|\nab^2\c\Dcal^{-1}P\|_{L_t^2L_x^2}+\|P\|_{L_t^2L_x^2}\big)\les
\De_0\c
\NN_1(P)
\eeaa
 in view of the second order elliptic
estimates of proposition  \ref{prop:secondorderHodge}.

Also, with the help
of \eqref{eq:bilin2newnew} of proposition  \ref{prop:5.7}
and proposition \ref{prop:errortype},
\beaa
\|J_2\|_{\PP^\ep}&=&\|\Gd\c \Gd
\c\Dcal^{-1}P\|_{\PP^\ep}\les\NN_2(\Dcal^{-1} P)\c
\|\Gd\c \Gd\|_{\PP^\ep}\\
&\les&\De_0^2\c \NN_2(\Dcal^{-1} P)\les \De_0^2\c \NN_1( P)
\eeaa
as desired.
\subsection{Main commutator result }
Using the results of the decomposition lemma \ref{le:decompbandmore}
together with  proposition \ref{prop:C1C2C3errorestim} we deduce
that all the commutators which we have denoted by $C(\Roc)$
 can be decomposed into the sum between an error type term
and another  which is  a perfect $\ddd_L$  derivative.
More precisely,
\begin{proposition} All the commutators $C(\Roc)$
 can be expressed as follows
\be{eq:main-comm-result}
C(\Roc)=\ddd_LP+E
\end{equation}
where $P$ and $E$ are tensors verifying:
$$\NN_1(P),\,\, \|E\|_{\PP^\ep}\les \De_o^2$$
for some $\ep>0$.
\label{prop:main-comm-result}
\end{proposition}
\begin{proof}:\quad
Indeed in view of proposition  \ref{prop:C1C2C3errorestim}
all commutators we have denoted by $C(\Roc)$ can
be expressed as sum of terms in $\PP^\ep$ with the 
exception of terms of the type 
$\nab\c\Dcal^{-1}\c\big(\b \c\Dcal^{-2}(\Roc)\big)$
and $\b \c\Dcal^{-2}(\Roc)$. The proposition
follows immediately by applying parts ii) and iii) of
 the decomposition lemma to these exceptional terms.
We also need to make use of the estimate
 $\NN_2(\Dcal^{-2}(\Roc))\les \RR_0+\De_0^2$
 established in proposition \ref{prop:L2C(Roc)}.
\end{proof}

\section{Main Estimates}In this section we 
provide the proof of theorem \ref{thm:Main}.
To do this we write down the equations
satisfied by our main  quantities in a  schematic form
to which we can then  apply the  main lemma \ref{le:mainlemma}.
We shall rely on all the bootstrap assumptions {\bf BA1}-{\bf BA3}
 as well as {\bf BA4}. In addition we shall make the following 
auxilliary bootstrap assumption,

{\bf BA5.} \qquad\qquad\qquad\qquad  $\|\Gd\c\Bd\|_{\PP^0}\le \De_0^2$

This additional  assumption will be used only to derive the desired estimates for 
$\ze, \mu, \trchb, \chibh$ and not for $\chih, \trch,\nab\trch$.
 We shall show, at the end of
this section,   that in fact 
$$\|\Gd\c\Bd\|_{\PP^0}\les\De_0\c\big(\II_0+\RR_0+\De_0^2\big)\le \f12 \De_0^2$$
and therefore {\bf BA5} was justified.

 We shall  make systematic use of 
 the results of the previous section concerning error type terms.

\subsection{Estimates for $\trch$, $\chih$.} \quad

{\bf Step 1}:\,\,\,{\sl Estimates for  $\|\trch-\frac{2}{r}\|_{L_t^\infty
L_x^\infty}$.}

Recall  that $\trch-\frac{2}{r}=(\trch-\overline{\trch})+(\overline{\trch}-\frac{2}{r})=V+W$.
According to proposition \ref{prop:travertrch} $V$ and $W$ verify 
the transport equation,
\beaa
\frac{d}{ds} W+\f12\trch W&=&\f12( V\c W + V^2)-\overline{\,|\chih|^2}
\\
\frac{d}{ds} V+\f12(\trch+\overline{\trch}) V&=&-\big(|\chih|^2-\overline{\,|\chih|^2}\big)
\eeaa
Applying  lemma \ref{le:transportformulaL2} as well as 
corollary \ref{corr:calcineq} and remembering also
the bootstarp
assumption  $\trch-\frac{2}{r}\les \De_0$ and  $\NN_1(\chih)+\|\chih\|_{L_x^\infty L_t^2}\les \De_0$,
 we derive,
\beaa
\|V\|_{L_x^\infty L_t^\infty }&\les& \|V(0)\|_{L^\infty(S_0)}+
\|  |\chih|^2-\overline{\,|\chih|^2}    \|_{L_x^\infty L_t^1}^2\les
\|V(0)\|_{L^\infty(S_0)}+\|\chih\|_{L_x^\infty L_t^2}^2+\|\chih\|_{L_t^2 L_x^2}\\
&\les&\|V(0)\|_{L^\infty(S_0)}+\De_0^2
\\
\|W\|_{ L_x^\infty  L_t^\infty }&\les&\|W(0)\|_{L^\infty(S_0)}+
\| V\c W   \|_{L_x^\infty L_t^1}^2+\| V\c V  \|_{L_x^\infty L_t^1}^2+\|
\overline{\,|\chih|^2}\|_{L_x^\infty L_t^1}^2\\
&\les&\|W(0)\|_{L^\infty(S_0)}+\De_0^2
\eeaa
Thus,
\be{eq:final-estimate-trch-Linfty}
\|\trch-\frac{2}{r}\|_{L_t^\infty L_x^\infty}\les\II_0+\De_0^2
\end{equation}
We can also easily check
the following,
\be{eq:final-estim-Ltrch}
\|\ddd_L(\trch-\frac{2}{r})\|_{ L_x^\infty L_t^1}\les\De_0^2
\end{equation}

{\bf Step 2}: {\sl Estimates for $\nab\trch$.}

Recall the transport equation \eqref{eq:transportnabtrch} for
$\nab\trch$,
$$\ddd_L(\nab\trch)= -\frac{3}{2} \trch\c\nab \trch-\chih\c \nab\trch-2\chih\c \nab \chih$$
and the  Codazzi equations \eqref{eq:dvichihgeod}, 
$$\div \chih=-\b+    \f12\,\nab\trch -\chih\c\ze.$$
We need also the renormalized  Bianchi equations,
\beaa
L(\rhoc)+\frac{3}{2}\trch\c\rhoc &=&\div\b-\ze\c\b
+\f12\chih\c(\nab\hot\ze+\f12\trchb\c\chih-\ze\hot\ze).\\
L(\sic)+\frac{3}{2}\trch\c\sic
&=&-\curl\b +\ze\wedge\b
+\f12\chih\wedge(\nab\hot\ze-\ze\hot\ze)
\eeaa
 Schematically,
$$\dcall\b=L(\rhoc,-\sic)+\frac{1}{r} R_0+\Gd\c\big( R_0+\frac{1}{r} \Gd+\nab \Gd +\Gd\c \Bd\big).$$
or,
$$\b=\dcall^{-1}L(\rhoc,-\sic)+\dcall^{-1}\bigg(\frac{1}{r} R_0+
\Gd\c\big( R_0+\frac{1}{r} \Gd+\nab \Gd +\Gd\c
\Bd\big)\bigg)$$ Therefore, setting $\Dcal^{-2}=
\dcall^{-1}\dcalll^{-1}$, $\Dcal^{-1}=\dcall^{-1}$ and $L(\rhoc,-\sic)=\ddd_L\Roc$
as in the previous section,
\beaa
\chih&=&\Dcal^{-2}\ddd_L\Roc+\Dcal^{-2}\bigg( \frac{1}{r} R_0+  
\Gd\c\big(\frac{1}{r} \Gd+\nab \Gd +\Gd\c \Bd\big)\bigg)
\\
&+&\Dcal^{-1}(\nab \trch+\Gd\c \Gd)
\eeaa
or, taking the covariant derivative,
\beaa
\nab\chih&=&\nab\c\Dcal^{-2}\ddd_L\Roc+\nab\c\Dcal^{-2}\bigg( \frac{1}{r} R_0+  
\Gd\c\big(\frac{1}{r} \Gd+\nab \Gd +\Gd\c \Bd\big)\bigg)
\\
&+&\nab\c\Dcal^{-1}(\nab \trch+\Gd\c \Gd)
\eeaa
Observe that 
$F_1=\nab\c\Dcal^{-2}\bigg(\Gd\c\big( \frac{1}{r} \Gd+\nab \Gd +\Gd\c \Bd\big)\bigg)$ and
$F_2=\Dcal^{-1}(\Gd\c
\Gd)$ are of strong error type,  
$\| F_1\|_{\PP^\ep},\,\| F_2\|_{\PP^\ep}\les\De_0^2$.

Thus, since clearly $\|\nab \c\Dcal^{-2}(\frac{1}{r} R_0)\|_{\PP^\ep}\les \RR_0+\De_0^2$,
\be{eq:estimchi1}
\nab\c \chih=-\nab\c\Dcal^{-2}\c\ddd_L(\Roc)+\nab\c\Dcal^{-1}(\nab\trch) +F,\quad \| F\|_{\PP^\ep}\les
\RR_0+\De_0^2
\end{equation}
We now write, using the commutator notation of definition \ref{def:maincommform},
$$\nab\c\Dcal^{-2}\c\ddd_L(\Roc)=\ddd_L\big(\nab\c\Dcal^{-2}(\Roc)\big)+C(\Roc)$$
According to proposition \ref{prop:C1C2C3errorestim},
$$C(\Roc)=\nab\c\Dcal^{-1}(\b\c\Dcal^{-2}\Roc)+\b\c\Dcal^{-2}\Roc+\err$$
On the other hand since,  in view of the decomposition lemma
 \ref{le:decompbandmore}, there
exists   tensors $P', E'$  such that,
\beaa
C(\Roc)&=&\ddd_L P'+E'\\
\NN_1(P'),\, \|E'\|_{\PP^\ep} &\les& \De_0\c \NN_2(\Dcal^{-2}\Roc)\les\De_0^2.
\eeaa
with the last inequality due to   proposition \ref{prop:L2C(Roc)}. Therefore,
\beaa
\nab\c\Dcal^{-2}\c\ddd_L(\Roc)&=&\ddd_L P+\err\\
P&=&\nab\c\Dcal^{-2}(\Roc)+P'\\
\NN_1(P)&\le&\RR_0+\De_0^2
\eeaa
Therefore, back to \eqref{eq:estimchi1} and setting $M=\nab\trch$,
\bea
\nab\c \chih&=&\ddd_L P+\nab\c\Dcal^{-1}M+E \label{eq:decom-nabchi}\\
\NN_1(P),\,\,\|E\|_{\PP^0}&\le&\RR_0+\De_0^2\nn
\eea
We also rewrite  the transport equation for $M=\nab\trch$
in the form,
\be{eq:transport-M}
\ddd_L M+\frac{3}{2}\trch M
=\Gd\c\big( \ddd_L P+M+ \nab\c\Dcal^{-1}M+E\big)
\end{equation}
We  are now in a position
to apply proposition \ref{rem:mainlemmaapplication}
to \eqref{eq:transport-M}. Using also the boundedness of
$\nab\c\Dcal^{-1}$ on $\PP^0$ and the bootstrap
bound  $\NN_1(\Gd)+\|\Gd\|_{L_x^\infty L_t^2}\les \De_0$ we  derive,
\beaa
\|M\|_{\BB^0}&\les&\|M(0)\|_{B^0_{2,1}(S_0)} +
\De_0\c\big(\NN_1(P)+\|M\|_{\PP^0}+
\|\nab\c\Dcal^{-1} M\|_{\PP^0}+\|E\|_{\PP^0}\big)\\
&\les&\II_0+\De_0\c \big(\NN_1(P)+\|E\|_{\PP^0}+\|M\|_{\PP^0}\big)
\les \II_0+\De_0\c \big(\RR_0+\De_0^2+\|M\|_{\BB^0}\big)
\eeaa
Thus, for small $ \De_0$, we deduce 
\be{eq:estim-for-Mchi}
\|M\|_{\BB^0}\les \II_0+\De_0^2
\end{equation}

We can also derive an estimate for $\|\nab\trch\|_{L_x^2 L_t^\infty}$
directly from  the transport equation \eqref{eq:transportnabtrch}
which we write in the form, 
$$\ddd_L M+\frac{3}{2}\trch M=\Gd\c M+\Gd\c\Gd$$
Applying to it lemma \ref{le:transportformulaL2} we derive,
\beaa
\|M\|_{L_x^2 L_t^\infty}&\les &\|M(0)\|_{L^2(S_0)}+\|\Gd\c M+\Gd\c\Gd\|_{L_x^2 L_t^1}\\
&\les&\II_0+\|M\|_{L_x^2 L_t^\infty}\c\|\Gd\|_{L_x^\infty L_t^1}+\|\Gd\|_{L_x^\infty L_t^2}^2\\
&\les&\II_0+\De_0\c \|M\|_{L_x^2 L_t^\infty}+\De_0^2
\eeaa
Thus, for sufficiently small $\De_0$,
\be{eq:final-estim-L-funnynorm}
\|M\|_{L_x^2 L_t^\infty}\les \II_0+\De_0^2
\end{equation}

{\bf Step 3}:\quad  {\sl Estimates for $\NN_1(\chih)$  and  $\|\chih\|_{ L_x^\infty L_t^2}$.}

We first estimate
$\NN_1(\chih)=\|\chih\|_{L_t^2L_x^2}+\|\nab\chih\|_{L_t^2L_x^2}+\|\ddd_L\chih\|_{L_t^2L_x^2}$. 
To estimate $\|\chih\|_{L_t^2L_x^2}+\|\nab\chih\|_{L_t^2L_x^2}$ we make
use of the elliptic estimates
of proposition \ref{prop:Hodgeestimates} to
the 
 Codazzi equations \eqref{eq:dvichihgeod}, 
$$\div \chih=-\b+    \f12\,\nab\trch -\chih\c\ze.$$
Thus, since  in view of \eqref{eq:estim-for-Mchi}  we have  $\|\nab\trch\|_{L_t^2L_x^2}\les
\II_0+\De_0^2$, we deduce
\beaa
\|\chih\|_{L_t^2L_x^2}+\|\nab\chih\|_{L_t^2L_x^2}&\les& \|\b\|_{L_t^2L_x^2}+
\|\nab\trch\|_{L_t^2L_x^2}+ \|\Gd\c\Gd\|_{L_t^2L_x^2}\\
&\les&\RR_0+\II_0+\De_0^2+\|\Gd\|_{L_t^4L_x^4}^2\les \RR_0+\II_0+\De_0^2
\eeaa
To estimate $\|\ddd_L\chih\|_{L_t^2L_x^2}$ we make use
of the transport equation \eqref{eq:chihgeod}
$$ \nab_L\chih =-\trch\c \chih+\a $$
Thus,
\beaa
\|\ddd_L\chih\|_{L_t^2L_x^2}\les\|\trch\|_{L^\infty}\|\chih\|_{L_t^2L_x^2}+\RR_0\les
\II_0+\RR_0
\eeaa
Therefore,
\be{eq:final-estimate-NN1(chih)}
\NN_1(\chih)\les  \RR_0+\II_0+\De_0^2
\end{equation}

To estimate the sharp trace norm $\|\chih\|_{ L_x^\infty L_t^2}$
we make use of  corollary \ref{corr:funny-classical-trace} 
and the decomposition formula \eqref{eq:decom-nabchi} to deduce,
\beaa
\|\chih\|_{L_x^\infty L_t^2}&\les& \NN_1(\chih)+\NN_1(P)+
\|\nab\c\Dcal^{-1} M\|_{\PP^0}+\|E\|_{\PP^0}\\
&\les& \NN_1(\chih)+\NN_1(P)+\|E\|_{\PP^0}+
\| M\|_{\BB^0}.
\eeaa

Therefore,
in view of the estimates \eqref{eq:estim-for-Mchi} for $M=\nab\trch$,
\eqref{eq:decom-nabchi} for $P, E$ and \eqref{eq:final-estimate-NN1(chih)}
for $\NN_1(\chih)$ already derived,
\be{eq:final-estim-trace-chih}
\|\chih\|_{L_x^\infty L_t^2}\les\RR_0+\II_0+\De_0^2 
\end{equation}

We record the results obtained so far in the following,
\begin{proposition} Assuming the boot-strap assumptions {\bf BA1} --{\bf BA3}
 we can derive 
 the following estimates for $\trch$ and $\chih$,
\bea
\|\trch-\frac{2}{r}\|_{L_t^\infty L_x^\infty}&\les&\II_0+\De_0^2
\\
\|\ddd_L(\trch-\frac{2}{r})\|_{ L_x^\infty L_t^2}&\les&\De_0^2
\\
\|\nab \trch\,\|_{L_x^2 L_t^\infty}&\les& \II_0+\De_0^2\\
\|\chih\|_{L_x^\infty L_t^2}+\NN_1(\chih)&\les& 
\II_0+\RR_0+\De_0^2 \\
\|\nab\trch\|_{\BB^0}&\les& \II_0+\RR_0+\De_0^2
\eea
In particular we can choose $\II_0, \RR_0$ sufficiently small 
and check that the bootstrap assumptions {\bf BA1}, {\bf BA2}
 for $\trch$ and $\chih$ 
are verified with  $\De_0$ replaced by $\De_0/2$.
We also deduce  all the estimates
for $\trch$, $\chih$ in theorem \ref{thm:Main}.
\label{prop:mainstructurenabch}
\end{proposition}

\subsection{Estimates  for $\mu$, $\nab \ze$}\quad

As in the previous section the most difficult estimate is
that of $\|\mu\|_{\BB^0}$. we shall
start with it.

{\bf Step 1}:\quad {\sl Estimates for $\mu$.}

 Recall the transport
equation \eqref{eq:newmasstransportgeodesic} for $\mu=-\div\ze+
\f12\chih\c\chibh -\rho+|\ze|^2$
\beaa
\frac{d}{ds} \mu+\frac{3}{2}\trch \mu&=&\chih\c(\nab\hot\ze)+\f12 \trch\rhoc+
2\ze\c\nab\trch
\\ &-& 4\chih\c\ze\c\ze + \trch|\ze|^2 -\frac 1 4
\trchb|\chih|^2\\&=&\chih\c(\nab\hot\ze)+\f12 \trch\rhoc+ 2\ze \c \nab\trch +
 \frac{1}{r}(\Gd\c \Gd +\Gd\c
\Bd)\\&+& \Gd\c \Gd\c \Gd+\Gd\c \Gd\c \Bd
\eeaa
\begin{remark}
Clearly the terms  $\frac{1}{r} \Gd \c \Gd$ and  $\Gd\c \Gd\c \Gd$  are better than 
 $\frac{1}{r} \Gd \c \Bd$, respectively $\Gd\c \Gd\c \Bd$ and therefore can be omitted.
To simplify notations we shall also  drop the term of the form $\frac{1}{r} \Gd \c \Bd$. 
Indeed it is  easier to handle  than $\Gd\c \Gd\c \Bd$ and though it is only quadratic
 relative 
to the bootstrap parameter $\De_0$, unlike $\Gd\c \Gd\c \Bd$ which is cubic, this will 
in no way affect our proof.
\label{rem:drop1r}
\end{remark}
Thus, schematically,
\be{eq:schematictransportM}
\frac{d}{ds} \mu+\frac{3}{2}\trch \mu=\chih\c(\nab\hot\ze)+\f12 \trch\rhoc+\Gd\c\nab\trch
+\Gd\c \Gd\c \Bd
\end{equation}
We now want to express  the terms $\chih\c(\nab\hot\ze)$, and $ \trch\c\rhoc$
in a suitable form to which we can apply our sharp bilinear trace lemma,
as we have done in the transport equation for $\nab\trch$ in the previous  subsection.
We do this with the help of the Hodge system \eqref{eq:hodgeze}:

\beaa
\div\ze&=&-\mu-\rhoc+|\ze|^2,\qquad
\curl\ze= \sic\nn
\eeaa  which we write in the form,
$$
\dcal\ze=-( \rhoc, -\sic)-(\mu,0)=-(\rho, -\si)-(\mu,0)+\Gd\c  \Bd
$$
or, with $M=(\mu,0)$,
\be{eq:schematicze}\ze=-\dcal^{-1}(\rho, -\si)-\dcal^{-1}M+\dcal^{-1}(\Gd\c  \Bd)
\end{equation}
On the other hand we write the Bianchi
equation \eqref{eq:transportmodifiedbbgeod} for $\bboc$,
\beaa
\ddd_L(\bboc) &=&-\nab\rho+(\nab \si)^\star -2 (\nab\hot\ze)\c\ze   
+
3(\ze\c\rho-\ze^\star
\si)-\trch \bb\\
&+&2\ze\c(-\f12\trch\c\chibh -\f12 \trchb\c\chih
+\ze\hot\ze)-4\chi\c\chibh\c\ze\nn
\eeaa
 in the form, see remark \ref{rem:drop1r},
\be{eq:schematictranspbboc}
\ddd_L\bboc=\dcalll (\rho,\si)+\frac{1}{r} R_0+\Gd\c(R_0+\nab \Gd+\Gd\c \Bd)
\end{equation}
or
$$
(-\rho,\si)=J\c\dcalll^{-1}\c\ddd_L\bboc+J\c\dcalll^{-1}(\frac{1}{r}
R_0)+J\c\dcalll^{-1}\big(\Gd\c(R_0+\nab \Gd+\Gd\c \Bd)\big)
$$
where $J$ is the involution $(\rho,\si)\longrightarrow (-\rho,\si)$.
Therefore, combining it with \eqref{eq:schematicze}
\beaa
\ze&=&-\dcal^{-1}\c J\c\dcalll^{-1}(\ddd_L\bboc+\frac{1}{r} R_0)+ \dcal^{-1} M +\dcal^{-1}(\Gd\c \Bd)\\
&+&\dcal^{-1}\c J\c \dcalll^{-1}\bigg(\Gd\c(R_0+\nab \Gd+\Gd\c \Bd)\bigg).
\eeaa
Now  introduce  the operators,
\be{eq:caldd}
\Dcal^{-2}=\dcal^{-1}\c J\c\dcalll^{-1},\qquad \Dcal^{-1}=\dcal^{-1}
\end{equation}
and observe that, see proposition \ref{prop:errortype}, the
the term  $  \nab\c\Dcal^{-2}\big(\Gd\c(R+\nab \Gd+\Gd\c \Bd)\big)$   
is a strong error term. To deal with the term $\nab\c  \Dcal^{-1}(\Gd\c \Bd)$
we need to make use of the auxilliary assumption {\bf BA5}
made at the beginning of this section,
$$\|\Gd\c \Bd\|_{\PP^0}\les \De_0^2$$ as well as the boundedness
of $\nab\c  \Dcal^{-1}$ on $\PP^0$.
Hence also,
$$\|\nab\c  \Dcal^{-1}(\Gd\c \Bd)\|_{\PP^0}\les \De_0^2$$
Remark also that,
  $$\|\nab\Dcal^{-2}(\frac{1}{r} R_0)\|_{\PP^\ep}\les \RR_0+\De_0^2.$$ 
 Thus  we can  write, 
 \be{eq:derivzeschematic}
\nab\ze=- \nab\c\Dcal^{-2}(\ddd_L\bboc )+\nab\c \Dcal^{-1} M +E,
 \quad \|E\|_{\PP^0}\les\RR_0+ \De_0^2
\end{equation}
Next we   commute $\ddd_L$ with $\nab\c\Dcal^{-2}$  in \eqref{eq:derivzeschematic} exactly as in the previous subsection.
 Recall that, see  definition
 \ref{def:maincommform} and \eqref{eq:maincommform},  
$$
\ddd_L\c(\nab\c\Dcal^{-2}\Roc)=\nab\c\Dcal^{-2}\c\ddd_L(\Roc)+C(\Roc)
$$
and, see proposition \ref{prop:main-comm-result},
$$
C(\Roc)=\ddd_L P'+E',\qquad \NN_1(P'), \,\,\, \|E'\|_{\PP^\ep}\les \De_0^2.
$$
Hence, setting $P=P'+ \nab\c\Dcal^{-2}(\Roc)$ and   recalling
 the estimate, see   proposition \ref{prop:L2C(Roc)}, 
$\NN_1\big(\nab\c\Dcal^{-2}(\Roc)\big)\les\RR_0+\De_0^2$  we can write,
\bea
\nab\ze&=&\ddd_L (P) +\nab\c\Dcal^{-1}M+
\nab\c\Dcal^{-2}(\frac{1}{r}R_0)+ E,\label{eq:strchi2}\\
 \NN_1(P)&\les& \RR_0+\De_0^2,\quad
\|E\|_{\PP^0}\les
\RR_0+\De_0^2
\eea
Going back to the transport equation  \eqref{eq:schematictransportM} for $M$
we replace $\nab\hot\ze$ according to \eqref{eq:strchi2}, and derive
$$
\frac{d}{ds} M+\frac{3}{2}\trch M=\Gd\c\big(\ddd_L (P) +\nab\c\Dcal^{-1}M+ E\big)
+\f12 \trch\rhoc+\Gd\c\nab\trch
+\Gd\c \Gd\c \Bd.
$$
Now observe that, according to the auxilliary assumption {\bf BA5},
$\Gd\c \Bd $ is a weak  error type term and can be incorporated in $E$,
\be{eq:schematictransportMagain}
\frac{d}{ds} M+\frac{3}{2}\trch M=\Gd\c\big(\ddd_L (P) +\nab\c\Dcal^{-1}M+\nab\trch+E\big)
+\f12 \trch\rhoc
\end{equation}
It remains to deal with the term $\trch\rhoc$.
 Returning to \eqref{eq:schematictranspbboc}
we write,
\beaa
(\rho,\si)&=&\dcalll^{-1}\ddd_L\bboc +\dcalll^{-1}\big(\frac{1}{r} R_0)+\dcalll^{-1}\big(\Gd\c(R_0+\nab \Gd+\Gd\c
\Bd)\big)\\
&=&\ddd_L(\dcalll^{-1}\bboc)+[\ddd_L, \dcalll^{-1}](\Roc)+\dcalll^{-1}\big(\frac{1}{r}
R_0)\\
&+&\dcalll^{-1}\big(\Gd\c(R_0+\nab \Gd+\Gd\c
\Bd)\big)
\eeaa
In view of proposition
\ref{prop:comm[DL]first} the commutator $[\ddd_L, \dcalll^{-1}](\Roc)$
is a strong error terms and so is $\dcalll^{-1}\big(\Gd\c(R_0+\nab \Gd+\Gd\c
\Bd)$ in view of proposition \eqref{prop:errortype}. Clearly,
$\dcalll^{-1}\big(\frac{1}{r}
R_0)\les \RR_0$.

Hence, 
$$
(\rho,\si)=\ddd_L(\dcalll^{-1}\bboc) +e',\qquad \|e'\|_{\PP^0}\les \RR_0+\De_0^2
$$
Observe that, in view of proposition \ref{prop:RocD-1L2commestim},
  $\NN_1(\dcalll^{-1}\bboc)\les \RR_0+\De_0^2$. On the other hand
$\rho$ and $\rhoc$ differ by a terms of the form $\Gd\c\Bd$ which 
is a weak error term.
 We have thus established
the following representation
formula,
\be{eq:structLPrho}
\rhoc=\ddd_L p'+e', \qquad \NN_1(P'), \|e'\|_{\PP^0} \les \RR_0+\De_0^2
\end{equation}

Going back to \eqref{eq:schematictransportMagain} we deduce, schematically,
$$
\frac{d}{ds} M+\frac{3}{2}\trch M=(\Gd+\frac{1}{r})\c\ddd_L (P)+\Gd\c \big(
\nab\c\Dcal^{-1}M+\nab\trch+E\big)+\frac{1}{r}E
$$
with $E$ verifying the  estimate  $\|e'\|_{\PP^0}\les \RR_0+\De_0^2$. 
We have just proved the following:
\begin{proposition} The mass term $\mu=M$
verifies a transport equation
of the form,
\be{eq:schematictransportMagain-2}
\frac{d}{ds} M+\frac{3}{2}\trch M=(\Gd+\frac{1}{r})\c\ddd_L (P)+
\Gd\c \big( \nab\c\Dcal^{-1}M+\nab\trch+E\big)+\frac{1}{r} E
\end{equation}
 with $P$ and $E$ verifying,
\be{eq:goodestmP'}
\NN_1(P)\les \RR_0+\De_0^2,\qquad \|E\|_{\PP^0}\les\RR_0+ \De_0^2
\end{equation}
\end{proposition}
We   use \eqref{eq:schematictransportMagain-2} to  estimate $M$ with the help of 
proposition
\ref{prop:mainlemmaapplication}. 
 Indeed setting
\beaa
F_1&=& \Gd+\frac{1}{r},\quad F_2=\frac{1}{r},\qquad F_3=\Gd\\
W&=&\nab\c \Dcal^{-1}M+\nab\trch + E
\eeaa
we can rewrite  \eqref{eq:schematictransportMagain-2} in the form, 
$$\ddd_L M+\frac{3}{2}\trch M=F_1\c\ddd_LP+  F_2\c E+F_3\c W.$$
Thus, making use of the estimates for $\|\nab\trch\|_{\BB^0}\les \II_0+\De_0^2$
already derived in the previous subsection,
\beaa
\|M\|_{\BB^0}&\les&\|M(0)\|_{B^0_{2,1}(S_0)} +
\big(\NN_1(F_1)+\|F_1\|_{L_x^\infty L_t^2}\big)\c\NN_1(P)\\
&+&\big(\NN_1(F_2)+\|F_2\|_{L_x^\infty L_t^2}\big)\|E\|_{\PP^0}
+\big(\NN_1(F_3)+\|F_3\|_{L_x^\infty L_t^2}\big)\|H\|_{\PP^0}\\
 &\les&\|M(0)\|_{B^0_{2,1}(S_0)}+\RR_0+\De_0^2+\De_0\c
\|W\|_{\PP^0}\\
\|W\|_{\PP^0}&\les&\|M\|_{\PP^0}+\|\nab\trch\|_{\PP^0}+\RR_0+\De_0^2\les\|M\|_{\BB^0}+\II_0+\RR_0+\De_0^2
\eeaa
Hence,
$
\|M\|_{\BB^0}\les\II_0+\De_0\c\|M\|_{\BB^0}+\RR_0+\De_0^2
$
and therefore,
\be{eq:final-estim-mu}
\|M\|_{\BB^0}\les\II_0+\RR_0+\De_0^2.
\end{equation}

As in the previous subsection we need to derive
also an estimate for $\|M\|_{L_x^2L_t^\infty}$.
We  do this 
 with the help of  lemma \ref{le:transportformulaL2} 
applied to the transport equation \eqref{eq:schematictransportM}.
\beaa
\|M\|_{L_x^2L_t^\infty}&\les &\|M(0)\|_{L^2(S_0)}+\RR_0+
\|\chih\c \nab\ze\|_{L_x^2L_t^1}\\&+&\|G\c\nab\trch\|_{L_x^2L_t^1}
+\|G\c G\c B\|_{L_x^2L_t^1}\\
&\les&\II_0+\RR_0+\De_0^2
\eeaa
Hence,
$$
\|M\|_{L_x^2L_t^\infty}\les \II_0+\RR_0+\De_0^2.
$$

{\bf Step 2:}\quad {\sl Estimates for $\NN_1(\ze)$.}
To estimate $\|\ze\|_{L_t^2L_x^2}+\|\nab\ze\|_{L_t^2L_x^2}$
we simply apply the elliptic estimates of proposition \ref{prop:Hodgeestimates}
to the elliptic system \eqref{eq:schematicze}.
Thus,
\beaa
\|\ze\|_{L_t^2L_x^2}+\|\nab\ze\|_{L_t^2L_x^2}&\les&\|R\|_{L_t^2L_x^2}+\|M\|_{L_t^2L_x^2}
+\|\Gd\c \Bd\|_{L_t^2L_x^2}\\
&\les&\II_0+\RR_0+\De_0^2+\|\Gd\|_{L_x^\infty L_t^2}\c\|\Bd\|_{L_x^2L_t^\infty}\\
&\les& \II_0+\RR_0+\De_0^2
\eeaa
Also, from \eqref{eq:zegeod}
$$\nab_L\ze+\trch\ze= R +\Gd\c \Bd
$$
and,
$$\|\ddd_L \ze\|_{L_t^2L_x^2}\les\|R\|_{L_t^2L_x^2}+
\|\ze\|_{L_t^2L_x^2}+\RR_0+\De_0^2\les
\II_0+\RR_0+\De_0^2$$
Thus,
\be{eq:final-estim-NN1ze}
\NN_1(\ze)\les\II_0+\RR_0+\De_0^2
\end{equation}

{\bf Step 3}:\quad {\sl Estimates for $\|\ze\|_{L_x^\infty L_t^2}$}

We take advantage once more  of the decomposition \eqref{eq:strchi2},
\beaa
\nab\ze&=&\ddd_L (P) +\nab\c\Dcal^{-1}M+\nab\c\Dcal^{-2}(\frac{1}{r}R_0)+ E,\\
 \NN_1(P)&\les& \RR_0+\De_0^2,\quad
\|E\|_{\PP^0}\les
\RR_0+\De_0^2
\eeaa
and corollary
 \ref{corr:funny-classical-trace}: 
\beaa
\|\ze\|_{L_x^\infty L_t^2}&\les& \NN_1(\ze)+\NN_1(P)+
\|\nab\c\Dcal^{-1} M+\nab\c\Dcal^{-2}(\frac{1}{r}R_0)+E\|_{\PP^0}\\
&\les&\II_0+\RR_0+\De_0^2
\eeaa

We gather all the results obtained above in the following:
\begin{proposition}  Assuming the boot-strap assumptions {\bf BA1} --{\bf BA4}
as well as the auxilliary bootstrap assumption {\bf BA5}
 we can derive 
 the following estimates for $\mu$, $\ze$:
\bea
\|\mu\|_{L_x^2L_t^\infty}&\les& \II_0+\De_0^2\\
\NN_1(\ze)&\les&\II_0+\RR_0+\De_0^2\label{eq:NN1ze}\\
\|\ze\|_{L_x^\infty L_t^2}&\les&\II_0+\RR_0+\De_0^2\label{eq:needtracele}\\
\|\mu\|_{\BB^0}&\les&\II_0+\RR_0+\De_0^2
\eea

In particular we can choose $\II_0, \RR_0$ sufficiently small 
and check that the bootstrap assumptions {\bf BA1}, {\bf BA2}
 for $\ze$ and $\mu$ 
are verified with  $\De_0$ replaced by $\De_0/2$.
We also deduce  all the estimates
for $\ze$, $\mu$ in theorem \ref{thm:Main}.
\label{prop:estim-zemu}
\end{proposition}

\subsection{ Estimates  for $(\trchb, \,\,\chibh)$}\quad

{\bf Step 1}:\quad{\sl $\BB^0$ estimates for $(\trchb , \chibh)$}

 Recall  the  transport equations \eqref{eq:trchbgeod}, \eqref{eq:chibhgeod}
 and Codazzi equations \eqref{eq:dvichibhgeod}.
\beaa
\frac{d}{ds}\trchb&=&-2\div\ze-\chi\c\chib+2|\ze|^2+2\rho\\
\nab_L\chibh&=&-\nab\hot
\ze-\f12(\trch\c\chibh+\trchb\c\chih)+\ze\hot\ze\\
\div\chibh&=&\f12\nab\trchb-\f12\trchb\c\ze+\ze\c\chibh+\underline{\b}\\
&=&\f12\nab\trchb-\f12\trchb\c\ze-\ze\c\chibh+\bboc
\eeaa
Recall that $\mu=-\div\ze+\f12\chih\c\chibh -\rho+|\ze|^2$.
Thus,
$$\frac{d}{ds}\trchb+\f12 \trch\c\trchb =2 \mu      -2\chih\c\chibh +4\rho$$
Also, since
$\frac{d}{ds}\trch+\f12 \trch^2 =-|\chih|^2,$
\be{eq:transtrchtrchb}\frac{d}{ds}(\trchb+\trch)+\f12 \trch(\trchb+\trch) =2
\mu+4\rho-2\chih\c\chibh-\chih\c\chih
\end{equation}
 On the other hand, recalling \eqref{eq:structLPrho}
$$
\rho=\ddd_L p+e', \qquad \NN_1(p),\,\, \|e'\|_{\PP^0} \les \RR_0+\De_0^2
$$
\beaa
\frac{d}{ds}(\trch+\trchb)+\f12 \trch(\trchb+\trch)&=&2\mu+4 \ddd_L p+ \chih\c \chibh+e'
\eeaa
Observe that according to the auxilliary assumption {\bf BA5}
 we have $\|\chi\c \chib\|_{\PP^0}\les \De_0^2$ while,
in view of the  estimates of the previous section 
$\|\mu\|_{\PP^0}\les \II_0+ \RR_0+\De_0^2$. 
Therefore   introducing  $e= 2\mu+ e'+\chi\c \chib$  and 
  $b=\trch+\trchb$,  we derive: 
\be{eq:schematictrchbtransport}
\frac{d}{ds}b+\f12 \trch b= 4 \ddd_L p+e
\end{equation}
where
$$
 \NN_1(p), \quad \|e\|_{\PP^0} \les\II_0+ \RR_0+\De_0^2.
$$
 We can write \eqref{eq:schematictrchbtransport} in the form,
$$\frac{d}{ds}b+\f12\trch b=F_1\c\ddd_Lp  +F_2\c  e$$
where $F_1=4$, , $F_2=1$, $b=\trch+\trchb$,  and apply  proposition
\ref{prop:mainlemmaapplication}. 
\bea
\|b\|_{\BB^0}&\les&\|b(0)\|_{B^0_{2,1}(S_0)} + \NN_1(p) +\|e\|_{\PP^0}\nn\\
&\les&\II_0+\RR_0+\De_0^2\les \II_0+\RR_0+\De_0^2 \label{eq:estimforb1}
\eea
We can proceed exactly
in the same manner with the 
transport equation for $\chibh$,
\beaa
\nab_L\chibh+\f12\trch\c\chibh&=&-\nab\hot
\ze-\f12\trchb\c\chih+\ze\hot\ze\\
&=&-\nab\hot
\ze-\f12 b\c\chih+\frac{1}{r} \c\chih+  \Gd\c\Gd
\eeaa
or, recalling  \eqref{eq:strchi2},
\beaa
\nab\hot\ze&=&\ddd_L (P) +\nab\c\Dcal^{-1}\mu+
\nab\c\Dcal^{-2}(\frac{1}{r}R_0)+ E',\\
 \NN_1(P)&\les& \RR_0+\De_0^2,\quad
\|E'\|_{\PP^0}\les
\RR_0+\De_0^2
\eeaa
Therefore, writing $\trchb=b-\trch=b-(\trch-\frac{2}{r})+\frac{2}{r}$
and including $b$ among the terms we have denoted by $\Bd$,
we can write,
\beaa
\ddd_L(\chibh)+\f12\trch\c\chibh&=&-\ddd_L P+\nab\c\Dcal^{-1}\mu+ \frac{1}{r}\big( \chih
+\nab\c\Dcal^{-2}(R_0)\big)+    \Gd\c\Bd+E'
\eeaa
Observe that the term $$E=\nab\c\Dcal^{-1}\mu+ \frac{1}{r}\big( \chih
+\nab\c\Dcal^{-2}(R_0)\big)+    \Gd\c\Bd+E'$$
verifies the estimate,
\beaa
\|E\|_{\PP^0}&\les& \|\mu\|_{\PP^0}+\|\chih\|_{\PP^0}+\|\nab\c\Dcal^{-2}(R_0)\|_{\PP^0}
+\|\Gd\c\Bd\|_{\PP^0}+\|E'\|_{\PP^0}\\
&\les&\II_0+\RR_0+\De_0^2
\eeaa
Therefore,
\bea
\ddd_L(\chibh)+\f12\trch\c\chibh&=&-\ddd_L P+E\label{eq:nice-structure-chibh}\\
\NN_1(P),\,\, \|E\|_{\PP^0}&\les&\II_0+\RR_0+\De_0^2\nn
\eea
to which we can apply as before  proposition
\ref{prop:mainlemmaapplication} to derive,
\beaa
\|\chibh\|_{\BB^0}&\les&\|\chibh(0)\|_{B^0_{2,1}(S_0)}+\NN_1(P)+
\|E\|_{\PP^0}\\
&\les&\II_0+\RR_0+\De_0^2
\eeaa
i.e.,
\be{eq:estim-for-chbh}
\|\chibh\|_{\BB^0}\les \II_0+\RR_0+\De_0^2
\end{equation}

{\bf Step 2}:\quad{\sl Estimates in the norm $L_x^2L_t^\infty$}
We shall  prove the following estimate:
\bea
\|\trchb+\frac{2}{r}\|_{L_x^2L_t^\infty}+\|\chibh\|_{L_x^2L_t^\infty}&\les &\II_0+\De_0^2
\label{eq:estimtrchbLx2Ltinfty}
\eea
We consider once more the transport  equation \eqref{eq:transtrchtrchb}
$$\frac{d}{ds}(\trchb+\trch)+\f12 \trch(\trchb+\trch) =2
\mu-4\rho +\Gd\c\Bd,$$
and apply to it lemma \ref{le:transportformulaL2}.
Taking also into account the estimates for $\mu$ derived
in the previous subsection,
\beaa
\|\trch+\trchb\|_{L_x^2L_t^\infty}&\les& \|(\trch+\trchb)(0)\|_{L_x^2(S_0)}
+\|\mu\|_{L_x^2L_t^1}+\RR_0\\
&+&\|\Gd\c\Bd\|_{L_x^2L_t^1}
\les\II_0+\RR_0+\De_0^2
\eeaa
Thus,
\beaa
\|\trchb+\frac{2}{r}\|_{L_x^2L_t^\infty}&\les&
 \|\trch-\frac{2}{r}\|_{L_x^2L_t^\infty}+\II_0+\RR_0+\De_0^2\\
&\les&\II_0+\RR_0+\De_0^2
\eeaa
as desired.
The corresponding estimates\footnote{Alternatively one can use the Codazzi
 equations for $\chibh$.} for $\chibh$ 
can be derived in the same manner from
its  transport equation,
$$\nab_L\chibh  +\f12 \trch\c\chibh  =-\nab\hot
\ze-\f12\trchb\c\chih+\ze\hot\ze.$$
We also observe, taking into account the estimates
we have already derived for $\trch, \chih,\ze,\trchb,\chibh$
\beaa
\|\ddd_L\chibh\|_{L_t^2L_x^2}&\les &\|\trch\c\chibh\|_{L_t^2L_x^2}+\|\nab\ze\|_{L_t^2L_x^2}
+\|\trchb\c\chih\|_{L_t^2L_x^2}+\|\Gd\|_{L_t^4L_x^4}^2\\
&\les& \II_0+\RR_0+\De_0^2+\|\trchb\|_{L_x^2L_t^\infty}\|\chih\|_{L_x^\infty
L_t^2}\les\II_0+\RR_0+\De_0^2
\eeaa
Similarily, from  \eqref{eq:transtrchtrchb},
$$\|\ddd_L(\trchb+\frac{2}{r})\|_{L_t^2L_x^2}\les \II_0+\RR_0+\De_0^2.$$

We gather  the results obtained above in the following:
\begin{proposition}  Assuming the boot-strap assumptions {\bf BA1} --{\bf BA4}
 as well as the auxilliary assumption {\bf BA5} 
 we can derive the following estimates for $\trchb, \chibh$:
\bea
\|\trchb +\frac{2}{r}\|_{\BB^0}&\les& \II_0+\RR_0+\De_0^2 \\
\|\chibh\|_{\BB^0}&\les&\II_0+\RR_0+\De_0^2\\
\|\trchb+\frac{2}{r}\|_{L_x^2L_t^\infty}+
\|\chibh\|_{L_x^2L_t^\infty}&\les &\II_0+\RR_0+\De_0^2\\
\|\ddd_L\big(\trchb+\frac{2}{r}\big)\|_{L_t^2L_x^2}&\les& \II_0+\RR_0+\De_0^2\\
\|\ddd_L\chibh\|_{L_t^2L_x^2}&\les &\II_0+\RR_0+\De_0^2.
\eea
In particular we can choose $\II_0, \RR_0$ sufficiently small 
and check that the bootstrap assumptions {\bf BA3}, {\bf BA4}
 for $\trchb$ and $\chibh$ 
are verified with  $\De_0$ replaced by $\De_0/2$.
We also deduce  all the estimates
for $\trchb$, $\chibh$ in theorem \ref{thm:Main}.
\label{prop:estim- Bad}
\end{proposition}
It only remains to verify the auxilliary assumption {\bf BA5}.
To achieve that we observe that, in view of \eqref{eq:schematictrchbtransport} 
and \eqref{eq:nice-structure-chibh} the terms we have denoted by 
 $\Bd=(\trchb+\frac{2}{r},\, \chibh) $  have the following structure
\beaa
\ddd_L\Bd+\f12 \trch \Bd& =&-\ddd_LP+E\\
\NN_1(P),\,\, \| E\|_{\PP^0}&\les& \II_0+\RR_0+\De_0^2
\eeaa
To estimate $\|\Gd\c \Bd\|_{\PP^0}$
we make use of proposition \ref{prop:mainlemmaapplication2} as follows:
Observe that the tensor  $W=P+\Bd$ verifies the transport equation,
$$\ddd_L W+\f12 \trch W=\f12 \trch  P+E.$$
Clearly $\|\trch P\|_{\PP^0}\les \II_0+\RR_0+\De_0^2$ and thus
the term $\f12 \trch  P$  can be  incorporated in $E$. According
to  proposition \ref{prop:mainlemmaapplication} we have
\beaa
\|\Gd\c W\|_{\PP^0}&\les &
\big( \NN_1(\Gd)+\|\Gd\|_{L_x^\infty L_t^2}\big)\c
\big ( \|W(0)\|_{B^0_{2,1}(S_0)}+\|E\|_{\PP^0}\big)\\
&\les&\De_0\c\big(\II_0+\RR_0+\De_0^2\big)
\eeaa
Therefore,
\beaa
\|\Gd\c \Bd\|_{\PP^0}&\les &\|\Gd\c W\|_{\PP^0}+\|\Gd\c P\|_{\PP^0}\\
&\les&\De_0\c\big(\II_0+\RR_0+\De_0^2\big)+\NN_1(G)\c \NN_1(P)\\
&\les&\De_0\c\big(\II_0+\RR_0+\De_0^2\big)\le \f12 \De_0^2
\eeaa
provided that $\II_0,\RR_0$ are sufficiently small relative to $\De_0$.
This justifies the auxilliary  bootstrap assumption {\bf BA5}
 and shows that the results of propositions 
 \ref{prop:estim-zemu} and \ref{prop:estim- Bad}
are valid without that additional assumption.

Together the  propositions  \ref{prop:mainstructurenabch}, 
\ref{prop:estim-zemu} and \ref{prop:estim- Bad} finish the proof
of our main theorem.

\section{Appendix} The purpose of this appendix is to prove 
 the commutator lemma \ref{le:pcomm-postponed}. This requires 
a quick recall of some of  the results of  \cite{KR2} for surfaces. 
We note that all our surfaces $S=S_s$, defined by the geodesic foliation on $\HH$,
verify the condition {\bf WS} of proposition \ref{prop:weak-regular-cond}. In
particular they
satisfy the  weak regularity condition 
of \cite{KR2}.
\begin{definition}
We define the negative  fractional powers $\La^{-a}=(I-\lap)^{-a/2}$, $a> 0$, according to the formula,
\be{eq:defineLaga}
\La^{-a}f=c_{a}\int_0^\infty \tau^{\frac{a}{2}-1}e^{-\tau}U(\tau)f d\tau
\end{equation}
where the heat flow  $U(\tau)f =e^{\tau \lap} f$ is  defined by,
\be{eq:heatflow}\pr_\tau U(\tau)f -\lap U(\tau)f=0, \,\, U(0)f=f.
\end{equation}
and $c_a$  appropriate constants. In fact the constnts 
are such that the  standard composition formulas hold:
$$\La^{-a}\c\La^{-b}=\La^{-(a+b)},$$  which reflects 
 the semigroup properties of $U(\tau)$.
\end{definition}
In the following proposition we record the main resullts
from \cite{KR2} which we shall need  in the proof of lemma 
\ref{le:pcomm-postponed}.

\begin{proposition}\label{prop:GN-B-H}
Given a surface $S$ verifying {\bf WS}, with $\f12 \le r\le 2$, we have
 the following inequalities:

{\bf i)} \quad The following calculus inequalities hold for an arbitrary tensor
$F$ and any $2\le p<\infty$,
\be{eq:Gagl-Nir1}
\|F\|_{L^p(S)} \les \|\nab F\|_{L^2(S)}^{1-\frac 2p} \|F\|_{L^2(S)}^{\frac 2p}
+ \|F\|_{L^2(S)}
\end{equation}

{\bf ii)}\quad The following B\"ochner-type inequality holds for any scalar 
function $f$,
\be{eq:Bochconseq}
\|\nab^2 f\|_{L^2(S)}\les  \|\lap f\|_{L^2(S)}  +  I_a \c\|\nab f\|_{L^2(S)}
\end{equation}
where \be{eq: defineKga}
 I_{a}=1+K_a^{\frac 1{1-a}} + K_{a}^\f12, \qquad \qquad 
K_a=\|\La^{-a}(K-\frac{1}{r^2})\|_{L^2(S)} \qquad
\end{equation}
{\bf iii)}\quad  The heat flow $U(\tau)f=e^{\tau \lap}f$ defined by 
\eqref{eq:heatflow}
verifies the following properties:
\bea
\|U(\tau) F\|_{L^2(S)}&\les &\|F\|_{L^2(S)}\label{eq:l2heat1scalar}\\
\|\nab U(\tau) F\|_{L^2(S)}&\les& \tau^{-\f12}\|F\|_{L^2(S)}\label{eq:l2heat2scalar}
\\ \|\lap U(\tau) F\|_{L^2(S)}&\les& \tau^{-1}\|F\|_{L^2(S)}
\label{eq:l2heat3scalar}\\
\| U(\tau)\nab F\|_{L^2(S)}&\les& \tau^{-\f12}\|F\|_{L^2(S)}
\label{le:L2heatscalar}
\eea
Also, for $2\le p<\infty$,
\bea
\|U(\tau) f\|_{L^p(S)}\les \big(1+\tau^{-(1-2/p)}\big)\|f\|_{L^2(S)}\label{eq:heat-GN}
\eea
and the dual estimate, for $1<q\le 2$,
\bea
\|U(\tau) f\|_{L^2(S)}\les \big(1+\tau^{(1-2/q)}\big)\|f\|_{L^q(S)}\label{eq:dual-heat-GN}
\eea
\end{proposition}
\begin{proof}:\quad 
See part I of \cite{KR2}.
\end{proof}
\begin{proof} {\bf of lemma \ref{le:pcomm-postponed}}   \quad We have  to prove the following estimate,
\be{eq:}
\|[\La^{-a}, \nab_L]\, f\|_{L^1_t L^2_x} \les I_a^\ep \c \|f\|_{L^2_t L^2_x}
\end{equation}
with $$I_a=1+K_a^{\frac 1{1-a}} + K_{a}^\f12, \qquad \qquad 
K_a=\|\La^{-a}(K-\frac{1}{r^2})\|_{L^2(S)}.$$
Acoording to \eqref{eq:defineLaga} and Duhamel's 
 formula, 
\beaa
[\La^{-a}, \nab_L]\, f&=&c_{a}\int_0^\infty \tau^{\frac{a}{2}-1}e^{-\tau}
[U(\tau), \ddd_L]f d\tau\\
&=&c_{a}\int_0^\infty d\tau\,\, \tau^{\frac{a}{2}-1}e^{-\tau}\int_0^\tau 
U(\tau-\tau') [\lap,\ddd_L]U(\tau') f d\tau'
\eeaa
We now use the commutator formula  \eqref{eq:commlap}
which we write in the symbolic form,
\beaa
[\ddd_L,\lap]U(\tau') f&=&
\nab\big( (\frac{1}{r}+\Gd)\c\nab U(\tau') f\big)+\big(\nab\Gd+\Gd\c\Gd\big)\c\nab
U(\tau')f\\
&=&\nab \phi_1(\tau') +\phi_2(\tau')
\eeaa
Thus,
\beaa
[\La^{-a}, \nab_L]f\, 
&=&\Phi_1+\Phi_2\\
\Phi_1&=&c_{a}\int_0^\infty d\tau\,\, \tau^{\frac{a}{2}-1}e^{-\tau}\int_0^\tau 
U(\tau-\tau') \nab \phi_1(\tau')   d\tau'\\
\Phi_2&=&c_{a}\int_0^\infty d\tau\,\, \tau^{\frac{a}{2}-1}e^{-\tau}\int_0^\tau 
U(\tau-\tau')  \phi_2(\tau')   d\tau'
\eeaa
Using the heat flow  estimates of proposition \ref{prop:GN-B-H},
\beaa
\|\Phi_1\|_{L_t^1L_x^2}&\les& \int_0^\infty
 d\tau\,\, \tau^{\frac{a}{2}-1}e^{-\tau}\int_0^\tau \|
U(\tau-\tau') \nab \phi_1(\tau') \|_{L_t^1L_x^2}  d\tau'\\
&\les&\int_0^\infty
 d\tau\,\, \tau^{\frac{a}{2}-1}e^{-\tau}\int_0^\tau  |\tau-\tau'|^{-\f12}
\|  \phi_1(\tau') \|_{L_t^1L_x^2}  d\tau'
\eeaa
Now,  using  first H\"older for some $p>2$, sufficiently close to $2$, followed
by the calculus, Bochner inequalities, and heat flow estimates  of proposition \ref{prop:GN-B-H}
as well as assumption {\bf BA1},
we derive,  
\beaa
\|  \phi_1(\tau') \|_{L_t^1L_x^2}&\les&\|A\|_{L_t^2L_x^{\frac{2p}{p-2}}}\c \|\nab U(\tau'
)f\|_{L_t^2L_x^p} + \|r^{-1}\|_{L^2_t L^\infty_x} \|\nab U(\tau)\|_{L^2_t L^2_x}\\
&\les& \De_0\|\nab^2 U(\tau'
)f\|_{L_t^2L_x^2}^{1-\frac{2}{p}} \c \|\nab U(\tau')f\|_{L_t^2L_x^2}^{\frac{2}{p}}+
\|\nab U(\tau)\|_{L^2_t L^2_x}\\
&\les& \De_0\|\lap U(\tau') f\|_{L_t^2L_x^2}^{1-\frac 2p} \c \|\nab
U(\tau')f\|_{L_t^2L_x^2}^{\frac{2}{p}}+ (1+\De_0\, 
I_{a}^{1-\frac 2p})
\|\nab U(\tau' )f\|_{L_t^2L_x^2}\\ &\les &
\|f\|_{L_t^2L_x^2}\,\big(    \De_0\,{\tau'}^{-1+\frac 1p}  + \,
(1+\De_0\,I_a^{1-\frac{2}{p}})\,\c {\tau'}^{-\frac 12}\,\big)
\eeaa
Therefore, 
\beaa
\|\Phi_1\|_{L^1_t L^2_x } &\les & \|f\|_{L^2_t L^2_x} \,
\int_0^\infty \tau^{-\frac 34} e^{-\tau} 
\int_0^\tau (\tau-\tau')^{-\f12}\big( \De_0\,{\tau'}^{-1+\frac 1p} +(1+\De_0\,I_a^{1-\frac
2p})\c{\tau'}^{-\frac 12}\big)  \\
&\les & (1+ \De_0\, I_a^{1-\frac 2p})\, \|f\|_{L^2_t L^2_x}
\eeaa
or, with $p=\frac{2}{1-\ep}$ for some  small $\ep>0$,
\be{eq:estim-Phi1}\|\Phi_1\|_{L^1_t L^2_x } \les (1+ \De_0\, I_a^{\ep})\, \|f\|_{L^2_t L^2_x}
\end{equation}
Similarly, with the help of \eqref{eq:dual-heat-GN}, for some $1<r$ close to $1$,
\beaa
\|\Phi_2\|_{L_t^1L_x^2}&\les& \int_0^\infty
 d\tau\,\, \tau^{\frac{a}{2}-1}e^{-\tau}\int_0^\tau \|
U(\tau-\tau')  \phi_2(\tau') \|_{L_t^1L_x^2}  d\tau'\\
&\les&\int_0^\infty
 d\tau\,\, \tau^{\frac{a}{2}-1}e^{-\tau}\int_0^\tau\big( 1+
|\tau-\tau'|\big)^{-\frac{1}{r}+\f12}
\|  \phi_2(\tau') \|_{L_t^1L_x^r}  d\tau'
\eeaa
On the other hand, observing that $p=\frac{2r}{2-r}>2$ can be made arbitrarily close 
to to if $r>1$ is close to $1$, we can estimate $\phi_2(\tau')$ precisely
 as $\phi_1(\tau')$, 
\beaa
\|  \phi_2(\tau') \|_{L_t^1L_x^r}&\les&\|\nab\Gd+\Gd\c\Gd\|_{L_t^2L_x^2}\c
 \|\nab U(\tau') f\|_{L_t^2  L_x^{\frac{2r}{2-r}}}\\
&\les&\De_0\c  \|\nab U(\tau') f\|_{L_t^2 L_x^{p}}\\
&\les&
\De_0   \c\|f\|_{L_t^2L_x^2}\,\big(    \,{\tau'}^{-1+\frac 1p}  + \,
I_a^{1-\frac{2}{p}}\,\c {\tau'}^{-\frac 12}\,\big), \qquad p=\frac{2r}{2-r}
\eeaa
Thus, just as for $\Phi_1$, with $r=\frac{1}{1-\ep}$,
\be{eq:estim-Phi2}
\|\Phi_1\|_{L^1_t L^2_x }\les  \De_0\, I_a^{\ep}\, \|f\|_{L^2_t L^2_x}
\end{equation}
Combining \eqref{eq:estim-Phi1} with \eqref{eq:estim-Phi2} we deduce,
$$\|[\La^{-a}, \ddd_L] f\|_{L^1_t L^2_x }\les 
  (1+\De_0 I_a^{\ep})\, \|f\|_{L^2_t L^2_x}\les 
 I_a^{\ep}\, \|f\|_{L^2_t L^2_x}
$$
as desired.
\end{proof}

\end{document}